\documentclass[reqno,11pt]{amsart}
\usepackage{graphicx}
\usepackage{amssymb}
\usepackage{amsmath}
\usepackage{amsfonts}
\usepackage{amsthm}
\usepackage{color}
\usepackage{hyperref}
\usepackage[mathscr]{eucal}
\usepackage{enumerate}
\usepackage{times}
\usepackage{ulem}                    
\usepackage{enumitem}            

\numberwithin{equation}{section}
\allowdisplaybreaks[1]

\newtheorem{Def}{Definition}[section]
\newtheorem{Thm}[Def]{Theorem}

\newtheorem{Prop}[Def]{Proposition}
\newtheorem{Lemma}[Def]{Lemma}

\newtheorem{Corollary}[Def]{Corollary}

\newcommand{\beq}{\begin{equation}}
\newcommand{\eeq}{\end{equation}}
\newcommand{\Proof}{\begin{proof}}
\newcommand{\QED}{\end{proof} \noindent}

\newcommand{\R}{\mathbb{R}}

\newcommand{\1}{\mbox{\rm 1 \hspace{-1.05 em} 1}}

\newcommand{\N}{\mathcal{N}}
\newcommand{\U}{\mathcal{U}}
\newcommand{\A}{\mathcal{A}}
\newcommand{\B}{\mathcal{B}}
\newcommand{\D}{\mathcal{D}}
\newcommand{\F}{\mathcal{F}}

\newcommand{\visc}{\varepsilon}
\newcommand{\ka}{\kappa}
\newcommand{\sos}{\sigma^2}
\newcommand{\kati}{\tilde{\kappa}}

\title[Dissipative Relativistic Fluid Flow]
{Dissipative relativistic fluid flow:\\ A simple Lorentz invariant causal model capturing entropy shocks in its zero viscosity limit}

\author[A.\ Chaddha]{Adhiraj Chaddha$^{*}$ }
\address[*]{Department of Mathematics\\ City University of Hong Kong \\ SAR Hong Kong}
\email{achaddha2-c@my.cityu.edu.hk}

\author[M.\ Reintjes]{Moritz Reintjes$^{**}$ \\ \\ September 17, 2026}
\address[**]{Department of Mathematics\\ City University of Hong Kong \\ SAR Hong Kong}
\email{moritzreintjes@gmail.com}

\begin{document}

\begin{abstract}
Zero viscosity limits, by identifying the physical (Lax admissible) shocks, are a cornerstone in the design of analytical and numerical schemes for the study of classical shock waves. For relativistic fluid flow, however, the underlying Laplacian based  dissipation mechanism (so-called ``artificial viscosity'') violates Lorentz invariance, the fundamental principle of Special Relativity.   In this paper we show that replacing the Laplacian on conserved quantities by the wave operator on the fluid four-velocity alone, provides a simplest Lorentz invariant description of dissipative relativistic fluid flow. We prove the resulting equations are causal and well-posed in one spatial dimension. We establish non-linear dissipativity and global-in-time existence of smooth solutions with small data in $H^s(\R)\cap L^1(\R)$, based on a Kawashima compensated energy method; (this appears to be the first such result for relativistic dissipation in 1-D). Moreover, we prove shock waves have profiles (a unique viscous travelling wave approximation in $L^2$) if and only if the shock wave is Lax admissible, and we prove that entropy production of travelling wave solutions is positive if and only if they obey the speed of light bound.    
This establishes the dissipative relativistic Euler equations introduced in this paper as a viable model for the study of relativistic shock waves in the zero viscosity limit, both in analytical and numerical approaches, consistent with the laws of Relativity. 
\end{abstract}

\maketitle 

\setcounter{tocdepth}{1}
\tableofcontents

\section{Introduction}  \label{Sec_intro}

Shock waves are central to the study of Fluid Dynamics, as smooth solutions of the compressible Euler equations are known to break down in finite time to form shock discontinuities whenever the flow is sufficiently compressive \cite{BuckmasterShkollerVicol_1, BuckmasterDrivasShkollerVicol, Dafermos, Lax1, Sideris, Smoller}.  Zero viscosity limits play a fundamental role in the theory of shock waves, both, for generating the physically correct shock wave solutions and for studying shock waves as mollified limits \cite{Dafermos, Kruzkov}. For instance, in the celebrated work of Tartar and Diperna \cite{Tartar, Diperna} shock wave solutions of the compressible Euler equations in 1-D, satisfying the correct entropy (Lax admissibility) condition, are obtained as zero viscosity limits of smooth solutions to the compressible Euler equations with a Laplacian based dissipative pressure term, (see also \cite{ChenPerepelitsa, Dafermos, HuangYang, SzepessyXin} and references therein). Conversely, in numerical schemes it is often crucial to study shock waves as smooth solutions at a mollified level in order to avoid steep gradients at shock discontinuities and potential unphysical oscillations or total collapse of the simulation \cite{VonNeumannRichtmyer, LaxWendroff, Harten, Bemfica_et-al}.  This is most often accomplished by inserting the spatial Euclidean Laplace operator on the fluid variables into the Euler equation as a simplest underlying dissipation mechanism (so-called ``artificial viscosity''), which correctly captures the physics of compressible fluid flow in the small viscosity limit, without turning to the more complicated Navier-Stokes viscosity.

Relativistic fluids play an equally important role in Physics, with applications ranging from Cosmology to the study of stellar structures and black hole formation \cite{Choquet, Weinberg}. Shock waves form in relativistic fluids whenever the flow is sufficiently compressive \cite{Christodoulou, PanSmoller}. The existence theory of Lax \cite{Lax} and of Glimm \cite{Glimm} for constructing shock waves was extended to the {\it relativistic} Euler equations by Smoller and Temple in \cite{SmollerTemple}, and extended to General Relativity in \cite{GroahTemple}. Using classical Laplacian based artificial viscosity, Chen and Schrecker \cite{ChenSchrecker} as well as LeFloch and Yamazaki \cite{LeFlochYamazaki}, succeeded to construct relativistic shock waves by the vanishing viscosity method.  However, the relativistic Euler equations with Laplacian based dissipation violate Lorentz invariance\footnote{The Laplace operator is Galilei invariant, the basic invariance of Classical Mechanics, which is preserved by the classical Euler equations with classical artificial viscosity; but the Laplacian violates Lorentz invariance.}---the principle at the foundation of Special Relativity  ensuring the speed of light bound.  Since it is generally advantageous in approximation schemes to preserve the underlying physical principles, (in particular Lorentz invariance in Relativity), the following question arises: \textit{Does there exist a dissipation mechanism for the relativistic Euler equations, consistent with the laws of Special Relativity (Lorentz invariance and causality), and whose zero viscosity limits yield the correct Lax admissible entropy shock waves in $1-D$, while being simple enough to implement in analytical and numerical approximation schemes?}

A sophisticated Lorentz invariant formulation of dissipative relativistic fluid flow, as an extension of classical Navier-Stokes dissipation to Relativity, has been given and analysed by Freist\"uhler and Temple in their profound works \cite{FreistuehlerTemple_1, FreistuehlerTemple_2, FreistuehlerTemple_3, FreistuehlerTemple_4}, (see also \cite{Freistuehler, Freistuehler2, FreistuehlerReintjesTemple, Pellhammer, SroczinskiZumbrun}), superseding the pioneering but less physical models of Landau \cite{LandauLifshitz}, Eckart \cite{Eckart} and Israel-Stewart \cite{IsraelStewart}. Alternative models, equally sophisticated and compelling, have been given by Bemfica, Disconzi and Noronha in \cite{Disconzi1, Disconzi2, Disconzi3, Disconzi4, Disconzi5}, by Kovtun in \cite{Kovtun}, (collective known as the BDNK-model), and by Mueller-Ruggeri in \cite{MuellerRuggeri} (the Theory of Extended Thermodynamics). A profound global-in-time existence theory for dissipative second order hyperbolic systems, applicable to the equations in \cite{FreistuehlerTemple_1, FreistuehlerTemple_2, FreistuehlerTemple_3}, has been developed by Sroczinski in \cite{Sroczinski1, Sroczinski2, Sroczinski3, Sroczinski4}, see also \cite{FreistuehlerSroczinski2, FreistuehlerSroczinski, FreistuehlerReintjesSroczinski}, building on the pioneering work of Kawashima \cite{Kawashima_diss}. Numerical studies of the BDNK-model have been done in \cite{Pretorius1, Pretorius2, Pretorius3, Pretorius4}.  By capturing the underlying relativistic mechanics and thermodynamics to a high degree of accuracy, the above models may indeed be too sophisticated for some analytical or numerical approaches merely aimed at modelling relativistic shock waves near the zero viscosity limit.\footnote{Simpler proposals, closer in spirit to artificial viscosity, have been given in \cite[Thm. 3]{FreistuehlerSroczinski} based on Gudonov variables, (see also \cite{Baerlin}).}  

In this paper, inspired by the Freist\"uhler-Temple dissipation model, we introduce a new, arguably simplest model of Lorentz invariant dissipation for the relativistic Euler equations, based on the wave operator on the fluid four velocity alone, and we establish its consistency for the study of relativistic shock waves as its zero viscosity limits. 

To introduce our model of relativistic artificial viscosity, consider a fluid with mass/energy-density $\rho$, pressure $p(\rho)$ and fluid $4$-velocity $u^\mu$. The relativistic Euler equations governing compressible fluid flow without viscosity are
\beq\label{rel_Euler}
\partial_{\nu}T^{\mu \nu} = 0,
\eeq 
where $T^{\mu \nu}$ is the energy-momentum tensor of a perfect fluid,
\begin{equation}\label{perfect_fluid}
T^{\mu \nu}= (p+\rho)u^{\mu}u^{\nu} + p \eta^{\mu \nu},
\end{equation}
$\eta={\rm diag}(-1,1,1,1)$ is the Minkowski metric of Special Relativity, and the fluid $4$-velocity is assumed to satisfy the speed of light bound in the form of the constraint $u^\sigma u_\sigma =-1$.\footnote{We use natural units with the speed of light set to $c=1$. We use the Einstein summation convention of summing over repeated upper and lower indices from $0$ to $3$, e.g., $\partial_{\nu}T^{\mu \nu} = \sum_{\nu =0}^3 \partial_{\nu}T^{\mu \nu}$ is the spacetime divergence. We raise and lower indices with the Minkowski metric $\eta_{\mu\nu}$ and its inverse $\eta^{\mu\nu}$, e.g., $u_\mu = \eta_{\mu\nu} u^\nu$.}  In this paper we introduce and analyze the equation
\begin{equation} \label{rel_art_visc}
\partial_{\nu}T^{\mu \nu} = \visc\; \square u^{\mu},
\end{equation}
where $\square  = -\partial_{tt}^2 + \Delta$ is the D'Alembert operator and $\Delta = \partial_{xx}^2 +\partial_{yy}^2 + \partial_{zz}^2$ the Laplacian, both taken component-wise on $u^\mu$, and $\visc >0$ is a constant representing viscosity. We assume throughout the paper that $p=p(\rho)$ is a barotropic equation of state obeying the speed of light bound $0<\frac{dp}{d\rho}<1$. We refer to \eqref{rel_art_visc} for simplicity as the \textit{dissipative relativistic Euler equations}. In this paper we establish the basic consistency of \eqref{rel_art_visc} for the study of ralativistic shock waves: 

\begin{itemize}[leftmargin=.6cm]
\item[(I)] To begin, in Proposition \ref{Thm_decay}, we show that Fourier Laplace modes of \eqref{rel_art_visc} all decay and establish Lorentz invariance of \eqref{rel_art_visc}. Moreover, we show that \eqref{rel_art_visc} reduce to a classical Navier Stokes type equation in the Newtonian limit. 
\item[(II)] In our first theorem, Theorem \ref{Thm_profiles}, we prove that shock profiles exists if and only if the shocks satisfy the Lax admissibility condition. In the latter case, we prove there exists a unique profile approximating the shock discontinuity in $L^2$ in space and uniformly in time. Moreover, we establish positivity of entropy production on travelling waves. 
\item[(III)] The dissipative relativistic Euler equations \eqref{rel_art_visc} are not second order hyperbolic. They are a mixed first-second order system, due to the constraint $u^\sigma u_\sigma =-1$, (c.f., \eqref{linearized_RAV} below). By expressing in one spatial dimension the system as a symmetrizable hyperbolic first order system, we establish in Theorem \ref{Thm_existence_causality} well-posedness of the Cauchy problem for \eqref{rel_art_visc}, including continuous dependence on data. Moreover, in Theorem \ref{Thm_existence_causality} (iii), we prove that \eqref{rel_art_visc} is \textit{causal}. That is, propagation of information is bounded by the speed of light, a central requirement of Relativity violated by the Eckart model \cite{Eckart} and Laplacian based artifical viscosity. 
\item[(IV)] In Theorem \ref{Thm_global} we finally establish non-linear dissipativity of \eqref{rel_art_visc} in one spatial dimension, and prove global-in-time existence of  solutions with small data in $H^s(\R)\cap L^1(\R)$ near a constant density state. This is based on a Kawashima compensated energy argument, adapting ideas in  \cite{KawashimaShizuta} to our setting.    To the authors knowledge, this is the first global-in-time existence result for smooth solutions of the Cauchy problem of a causal relativistic dissipative fluid system, in one spatial dimension, the setting where classical shock wave theory is well-established and which is the goal of this paper.\footnote{The existence theory in \cite{Sroczinski1, Sroczinski2, Sroczinski3, Sroczinski4} applies in spatial dimension $d\geq 3$, and requires full rank of coefficient matrices. It hence does not apply to the mixed-order system \eqref{rel_art_visc}. Finite time blow-up of solutions in the Israel-Stewart model is shown in \cite{Disconzi6, Disconzi7}.}
\end{itemize}

By design, equation \eqref{rel_art_visc} correctly encode the principles of Special Relativity--Lorentz invariance and causality. Interestingly, even though the laws of Thermodynamics are not build into the dissipative relativistic Euler equation by design, the second law holds near shock waves in the sense that entropy production of travelling wave solutions is positive if and only if the wave obeys the speed of light bound, (see Theorem \ref{Thm_profiles}). 

Taken together, this establishes the dissipative relativistic Euler equations as a promising simplest model for the study of shock waves in the zero viscosity limit, simple enough to be implemented in analytical and numerical approaches, and consistent with the laws of Special Relativity.

\section{Statement of Results} \label{Sec_results}            

\subsection{Fourier-Laplace modes and Newtonian limit}
Our first proposition establishes decay of Fourier-Laplace modes of the linearization of  \eqref{rel_art_visc}, as a consistency check of dissipativity. For this, we derive in Section \ref{Sec_diss} the linearization of \eqref{rel_art_visc} at a steady state solution $(\rho_0,u^\mu_0)$ to be given by
\beq \label{linearized_RAV}
\begin{cases}
\partial_t \rho + \ka\; \nabla \cdot \vec{u}  = 0 \cr
\ka\; \partial_t u_j + \sos \; \partial_j \rho = \visc\; \Box u_j ,
\end{cases}   
\eeq
$j=1,...,3$, and where $\vec{u} \equiv (u^1,u^2,u^3)$ are the spatial components of the fluid $4$-velocity $u^\mu$, $\nabla=(\partial_{x}, \partial_{y},\partial_{z})$ is the spatial gradient, $\sigma \equiv \sqrt{p'(\rho_0)}\equiv \sqrt{\frac{dp}{d\rho}|_{\rho_0}} >0$ is the speed of sound, $\ka \equiv p(\rho_0) + \rho_0 >0$, and we assumed $u^\mu_0= (1,0,0,0)$ and $\rho_0>0$ is constant. A Fourier-Laplace (FL-) mode solution of \eqref{linearized_RAV} is of the form
\beq \label{FL-modes}
\begin{pmatrix} \rho \cr \vec{u} \end{pmatrix} = \Psi e^{\lambda t + i \vec{\xi} \vec{x}},
\eeq
where $\Psi \in \mathbb{C}^4$ is a fixed vector, and $\lambda \in \mathbb{C}$ and $\vec{\xi} \in \R^3$ are the FL-modes determined by the dispersion relation of \eqref{linearized_RAV}. Our first result establishes decay of FL-mode solutions, Lorentz invariance, and shows that \eqref{rel_art_visc} reduce to a classical Navier-Stokes type equations in the Newtonian limit.  

\begin{Prop} \label{Thm_decay}
\textbf{(i)} Equations \eqref{rel_art_visc} is Lorentz invariant.\\ 
\textbf{(ii)} Assume $p(\rho_0) + \rho_0 >0$ and  $0< p'(\rho_0) <1$. Then, for any $\vec{\xi} \in \R^3$ there exists some $\lambda \in \mathbb{C}$ such that the corresponding Fourier-Laplace mode \eqref{FL-modes} solves \eqref{linearized_RAV}. Moreover, if $\vec{\xi} \neq 0$, then $\lambda \neq 0$, and any Fourier-Laplace mode solution \eqref{FL-modes} of \eqref{linearized_RAV} with $\vec{\xi} \neq 0$ decays in the sense that ${\rm Re}(\lambda) <0$.\footnote{Since FL-modes are not Lorentz invariant, decay of FL-modes for linearizations at constant fluid $4$-velocities $u^\mu_0$ with non-zero spatial components cannot be inferred from Proposition \ref{Thm_decay}.} \\
\textbf{(iii)} Equation \eqref{rel_art_visc} converges (formally) under the classical limit ($c\to\infty$) to
\beq \label{class_art_visc}
\begin{cases}
\partial_t \rho + \nabla \cdot(\rho \vec{v}) = 0   
\cr
\rho \left( \partial_t \vec{v} + \vec{v}\cdot \nabla \vec{v} \right) + \nabla p = \visc \; \Delta \vec{v} , 
\end{cases}
\eeq
where $\rho(t,\vec{x})$ is the mass-density, $p(t,\vec{x})$ the pressure, $\vec{v}(t,\vec{x})\in \R^3$ the fluid velocity, $\nabla$ the spatial gradient and $\Delta \equiv \partial_{xx}^2 + \partial_{yy}^2 + \partial_{zz}^2$ the Laplacian. 
\end{Prop}

\subsection{Shock profiles and entropy production}
Our first main result establishes that shock profiles of the dissipative relativistic Euler equations in one spatial dimension exist and select the correct Lax admissible shock waves.  More precisely, substituting the traveling wave ansatz $\rho(t,x) = \rho(\zeta)$ and $u^\mu(t,x)= u^\mu(\zeta)$, for $\zeta = \frac{x-st}{\visc}$, $u^2=0=u^3$, $s\in \R$, into the dissipative relativistic Euler equations \eqref{rel_art_visc}, and expressing the normalization condition $u^\sigma u_\sigma=-1$ (imposing the speed of light bound) in the form $u^0 = \sqrt{1+u^2}$ for $u\equiv u^1$, we first solve in Section \ref{Sec_profiles} for the density as a function of $u$ explicitly, and write the dissipative relativistic Euler equations \eqref{rel_art_visc} equivalently as       
\begin{eqnarray}
\dot{u} &=& \tfrac{1}{1-s^2} \Big( T^{1\nu} - T^{1\nu}_R \Big) \zeta_\nu,\label{ODE} \\
\rho(u) &=& \tfrac{1}{s-v(u)} \Big( v(u) T^{1\nu}_R - T^{0\nu}_R \Big) \zeta_\nu , \label{density-function}   
\end{eqnarray}
where $v(u) \equiv \frac{u}{\sqrt{1+u^2}}$ is the classical fluid velocity, $\zeta_\nu \equiv (-s,1)$ and $\dot{u} \equiv \frac{du}{d\zeta}$. Equation \eqref{density-function} gives the density as an explicit function of $u$, and thereby reduces \eqref{rel_art_visc} to the scalar ODE  \eqref{ODE}, which leads to a significant simplification in the analysis of shock profiles of the dissipative relativistic Euler equations. 

Now, to state our result on shock profiles, we take $s\in \R$ to be the shock speed of two constant states $(\rho_L,u_L)$ and $(\rho_R,u_R)$, which satisfy the Rankine Hugoniot (RH) conditions       
\beq \label{RH-condt}
[T^{\mu\nu}] \zeta_\nu \equiv (T^{\mu\nu}_L - T^{\mu\nu}_R)  \zeta_\nu =0,
\eeq
where $\zeta_\nu = (-s,1)$, and $T^{\mu\nu}_{L/R}$ denotes the energy momentum tensor \eqref{perfect_fluid} evaluated at $(\rho_{L/R},u_{L/R})$ respectively. We show in Section \ref{Sec_profiles} that local solutions of \eqref{ODE} with initial data in between $u_L$ and $u_R$ extend to global solutions $u(\zeta)$, $\zeta \in (-\infty,\infty)$, which converge asymptotically to $u_L$ and $u_R$, because $u_L$ and $u_R$ are the only rest points of \eqref{ODE}. Now, by convention $u_L$ is the state on the left hand side of the shock surface $\{x-st=0\}$, and $u_R$ is the state on the right hand side of the shock surface, and thus the physically relevant solutions of \eqref{ODE} for the problem of shock profiles have asymptotics 
\beq \label{asymptotics}
\lim\limits_{\zeta \to - \infty} (\rho(\zeta),u^\mu(\zeta)) = (\rho_L,u^\mu_L)
\hspace{.5cm} \text{and} \hspace{.5cm} 
\lim\limits_{\zeta \to \infty} (\rho(\zeta),u^\mu(\zeta)) = (\rho_R,u^\mu_R).
\eeq
Theorem \ref{Thm_profiles} below establishes that the asymptotics \eqref{asymptotics} on solutions of \eqref{ODE} are equivalent to the Lax admissibility conditions for either 1- or 2-shocks, that is,
\begin{eqnarray} 
&\lambda_1(u_R) < s_1 < \lambda_1(u_L) \ \ \ \ \ \text{and}  \ \ \ \ \   s_1 < \lambda_2(u_R),  & \label{Lax-condition_1}  \\
\text{or} \ \ \ \ \ & \lambda_2(u_R)  < s_2 < \lambda_2(u_L) \ \ \ \ \ \text{and}  \ \ \ \ \   \lambda_1(u_L) <  s_2,  & \label{Lax-condition_2} 
\end{eqnarray}
where $s_i$ denotes the shock speed of the $i$-th shock curve, $i=1,2$, and 
\beq 
\lambda_1(u) \equiv \frac{v-\sigma}{1-v\sigma}
\hspace{.5cm} \text{and} \hspace{.5cm} 
\lambda_2(u) \equiv \frac{v+\sigma}{1+v\sigma} 
\eeq
are the characteristic speeds of the relativistic Euler equations \eqref{rel_Euler} in terms of the speed of sound $\sigma \equiv \sqrt{\frac{dp}{d\rho}}$.  Concerning consistency with the second law of thermodynamics, Theorem \ref{Thm_profiles} (III) shows that entropy production of the dissipative relativistic Euler equations is positive for travelling wave solutions obeying the speed of light bound. In summary, we prove in Section \ref{Sec_profiles} the following theorem.

\begin{Thm} \label{Thm_profiles}  
Assume a pressure law with $0<p'(\rho)<1$, $p''(\rho)\geq 0$ and $p(\rho)>0$ for $\rho>0$.  \\
\textbf{(I)} Assume $(\rho_L,u_L) \neq (\rho_R,u_R)$ are constant states with $\rho_L>0$, $\rho_R>0$, which satisfy the RH conditions \eqref{RH-condt} with shock speed $s$.  The following is equivalent:
\begin{itemize}[leftmargin=.25in]  \setlength\itemsep{.3em}    
\item[(i)] There exists a traveling wave solution $\rho(t,x) = \rho(\zeta)$ and $u^\mu(t,x)= u^\mu(\zeta)$, for $\zeta \equiv \frac{x-st}{\visc}$, of the dissipative relativistic Euler equations \eqref{rel_art_visc} with asymptotics \eqref{asymptotics}.
\item[(ii)] The first Lax condition \eqref{Lax-condition_1} holds with $s=s_1$ if $u_L$ and $u_R$ are separated by a 1-shock, and the second Lax condition \eqref{Lax-condition_2} holds with $s=s_2$ if $u_L$ and $u_R$ are separated by a 2-shock.                        
\end{itemize}      
\textbf{(II)} For any Lax-admissible shock wave solution of \eqref{rel_Euler} there exists a unique family of travelling wave solutions of \eqref{rel_art_visc} with $\visc$-independent center which converge to the Lax shock as $\visc \to 0$, $L^2$ in space and uniformly in time. This family is given by $\rho(t,x) = \rho(\zeta)$ and $u^\mu(t,x)= u^\mu(\zeta)$, with $\rho(\zeta)$ and $u(\zeta)$ solving \eqref{ODE}, and with $\zeta \equiv \frac{x-st-x_0}{\visc}$, where $x_0\in \R$ is the position of the shock discontinuity at $t=0$. \\
\textbf{(III)} Assume the first law of Thermodynamics, and consider traveling wave solutions to \eqref{rel_art_visc} in one spatial dimension of form $u^{\mu}=u^{\mu}(\xi)$ and $\rho=\rho(\xi)$ with normalization $u^\sigma u_\sigma =-1$, and $\xi\equiv x-ct$ for some constant $c\in \R$. Then, in the region where $\frac{du^\mu}{d\xi} \neq 0$, any such travelling wave has positive entropy production, if and only if $|c| < 1$; and entropy production is zero if and only if $|c|=1$. 
\end{Thm}

\subsection{Local well-posedness and causality}
We next address the question of causality and well-posedness of the dissipative relativistic Euler equations. Based on the constraint $u^\sigma u_\sigma=-1$, equations \eqref{rel_art_visc} are a mixed first-second order system, not a second order hyperbolic system of full rank, and standard well-posedness results do not apply to \eqref{rel_art_visc} as stated. In Section \ref{Sec_Cauchy} we rewrite the dissipative relativistic Euler equations \eqref{rel_art_visc}, subject to the constraint $u^\sigma u_\sigma=-1$, as a first order symmetrizable hyperbolic system in one spatial dimension, and identify an explicit symmetrizer. Based on Kato's existence theory \cite{Kato}, we then obtain well-posedness, (local-in-time existence and continuous dependence on data), in one spatial dimension. Interestingly, symmetric hyperbolicity and well-posedness is achieved without requiring positivity (or even non-negativity) of the density, an interesting property in light of the physical vacuum boundary problem \cite{HadzicShkollerSpeck}. 

Moreover, we prove that the dissipative relativistic Euler equations are causal. That is, propagation of information is bounded by the speed of light bound. This is a central pillar of Special Relativity, violated by the Eckart model \cite{Eckart}, (see the discussion in \cite{FreistuehlerTemple_1}), and by Laplacian based dissipation. 

We assume throughout that $p(\rho)$ is a smooth (at least $C^{s+2}$) function of $\rho$, defined on some open interval $I \subset \R$. We consider \eqref{rel_art_visc} with one-dimensional initial data
\beq \label{Cauchy-data}
\rho(0,x)=\rho_{0}(x), \qquad  u(0,x)=u_{0}(x), \qquad \partial_t u(0,x)= w_{0}(x).
\eeq 
Since $\rho_0\in H^s(\R)$ would imply $\rho_0$ vanish asymptotically, we introduce a fixed reference state $\bar{\rho} \in I$ and require $\rho_0 - \bar{\rho} \in H^s(\R)$.

\begin{Thm} \label{Thm_existence_causality}
Let $s \geq 2$ be an integer, $\D_s \subset H^s(\mathbb{R},\R)$ open and bounded, $\bar{\rho} \in I$ a constant density, and assume $\overline{\{\bar\rho + \rho_0(x) \mid \rho_0\in\mathcal{D}_s,\ x\in\mathbb R\}} \subset I$.\\
\noindent \textbf{(i)} Assume initial data $(\rho_0, u_0, w_0)$ with $\rho_0\in I$, $\rho_0 - \bar{\rho} \in \D_s$, $u_0 \in\D_{s+1}$, $w_0 \in\D_s$. Then there is a unique solution $(\rho,u)$ of \eqref{rel_art_visc}, subject to the constraint $u^\sigma u_\sigma=-1$, satisfying the initial data $\rho(0,x)=\rho_0(x)$, $u(0,x)=u_0(x)$, $\partial_t u(0,x)=w_0(x)$, defined on $[0, T]$ for some $T>0$, such that
\beq 
\begin{aligned}  \label{Cauchy_regularity_local}
\rho - \bar{\rho} &\in \; C([0,T]; \D_s) \cap C^1([0,T]; H^{s-1}(\R)), \cr
u &\in \; C([0,T]; \D_{s+1}) \cap C^1([0,T]; \D_{s})  \cap C^2([0,T]; H^{s-1}(\R)), 
\end{aligned}
\eeq 
and $T$ can be chosen common to all initial data sufficiently close to $(\rho_0,u_0,w_0)$. Moreover, if $\inf\limits_{x\in \R}\rho_0(x)>0$ then $\inf\limits_{x\in \R}\rho(t,x)>0$ for all $t>0$ sufficiently small. 

\noindent \textbf{(ii)} Solutions of \eqref{rel_art_visc} depend continuously on the initial data in $H^s \times H^{s+1}$, in the following sense: Consider $(\rho_0 - \bar{\rho},u_0,w_0)$ and a sequence $(\rho_0^n - \bar{\rho}, u_0^n, w_0^n)$, both in $\D_s\times\D_{s+1} \times \D_s$, $n\in \mathbb{N}$, and let $(\rho,u)$ and $(\rho^n, u^n)$ be the unique solutions of \eqref{rel_art_visc} with data $(\rho, u, \partial_t u)(0)=(\rho_0,u_0,w_0)$ and $(\rho^n, u^n, \partial_t u^n)(0)=(\rho^n_0, u^n_0, w^n_0)$. Assume the solution $(\rho,u)$ exists on a closed interval $[0, T_0]$. If 
$$
\lim_{n\to\infty}  \Big( \|\rho_0^n-\rho_0\|_{H^s} + \|u_0^n - u_0\|_{H^{s+1}} + \|w_0^n - w_0\|_{H^s} \Big)=0,
$$ 
then the solutions $(\rho^n, u^n)$ also exist on $[0, T_0]$ for sufficiently large $n$, and 
$$
\lim_{n\to\infty}  \Big(\|\rho^n-\rho\|_{H^s} + \|u^n - u\|_{H^{s+1}} + \|\partial_t u^n - \partial_t u\|_{H^s} \Big)=0,
$$ 
uniformly in $t \in [0,T_0]$.

\noindent \textbf{(iii)}  All solutions of \eqref{rel_art_visc} with regularity \eqref{Cauchy_regularity_local} are causal throughout their existence. That is, their backward Cauchy developments are contained inside backward light-cones.
\end{Thm}

\subsection{Global-in-time existence and non-linear dissipation}
Our third main result asserts that the local solutions furnished by Theorem  \ref{Thm_existence_causality} extent to global solutions for sufficiently small perturbations of a constant rest state. For the proof, we introduce a nonlinear change of variables, to write the symmetrizable hyperbolic $4\times 4$ system underlying Theorem  \ref{Thm_existence_causality} in balance-law form with a rank-$1$ damped forcing term. By adapting the ideas of Kawashima-Shizuta \cite{KawashimaShizuta} to our setting, we introduce a compensating cross term to transfer this damping to the remaining components, from which we obtain an energy estimate controlling the high-frequency part of solutions; (this builds upon the pioneering work of Kawashima and Shizuta \cite{KawashimaShizuta}, Kawashima \cite{Kawashima_diss}, and Kawashima-Umeda-Shizuta \cite{KawashimaUmedaShizuta}, as well as the later advancements \cite{HanouzetNatalini, Yong, BianchiniHanouzetNatalini, FreistuehlerSroczinski, FreistuehlerReintjesSroczinski} and the more general work \cite{Sroczinski1, Sroczinski2, Sroczinski3, Sroczinski4, FreistuehlerSroczinski2} partially based on para-differential operators).  At low frequencies, the linearized flow exhibits spectral decay which leads in combination with Duhamel's formula to decay at rate $O(t^{-\frac14})$.  The resulting low-frequency decay estimate and high-frequency compensated energy estimate suffice to close a small-data bootstrap argument providing a uniform Sobolev bound. This then allows the local solution to be continued for all time. For solutions to decay, positivity of the reference density state $\bar{\rho}$ and of the initial density $\rho_0$ is crucial.

\begin{Thm} \label{Thm_global} 
Let $s \geq 2$ be an integer, $I\subset(0,\infty)$ an open interval and $\bar{\rho}~\in~I$ a constant positive density state. Assume a pressure law $p\in C^{s+2}(I)$ with $0<p'(\bar\rho)<1$ and $\rho+p(\rho)>0$ for $\rho\in I$. Define $W_{0} \equiv  (\rho_0-\bar\rho, w_0, \partial_x u_0, u_0)$ for initial data \eqref{Cauchy-data}. Then there exists $\delta_0>0$ such that, if $W_0 \in L^1(\mathbb R;\mathbb R^4)\cap H^s(\mathbb R;\mathbb R^4)$ satisfies
$$
\delta \equiv  \|W_0\|_{L^1(\R)} + \|W_0\|_{H^s(\R)} \leq\delta_0,
$$
then the local solution furnished by Theorem \ref{Thm_existence_causality} extends uniquely to a global solution of \eqref{rel_art_visc}, subject to $u_\mu u^\mu=-1$, satisfying
\begin{align*}
\rho-\bar\rho &\in C\bigl([0,\infty);H^s(\mathbb R)\bigr) \cap C^1\bigl([0,\infty);H^{s-1}(\mathbb R)\bigr), \cr 
u &\in C\bigl([0,\infty);H^{s+1}(\mathbb R)\bigr) \cap C^1\bigl([0,\infty);H^s(\mathbb R)\bigr) \cap  C^2\bigl([0,\infty);H^{s-1}(\mathbb R)\bigr).
\end{align*}
This solution is dissipative in the sense that there exists a constant $C>0$ such that
\beq \label{non-lin_dissipation}
\|\rho(t)-\bar\rho\|_{H^s(\R)} + \|\partial_tu(t)\|_{H^s(\R)} + \|u(t)\|_{H^{s+1}(\R)} \leq C\delta(1+t)^{-1/4},
\eeq 
for all $t\geq 0$. After possibly decreasing $\delta_0$, the density remains in $I$ and satisfies
\[
\rho(t,x)\geq\frac{\bar\rho}{2}>0, \qquad \forall\; (t,x)\in[0,\infty)\times\mathbb R.
\]
The constants $C>0$ and $\delta_0 >0$ depend only on $s$, $\visc$, $\bar\rho$, and the pressure law $p(\cdot)$ on a neighborhhod of $\bar{\rho}$.
\end{Thm}

We prove Proposition \ref{Thm_decay} on the decay of Fourier-Laplace modes and Lorentz invariance of \eqref{rel_art_visc} in Section \ref{Sec_Lorentz_FL-decay}. Theorem \ref{Thm_profiles} is proven in Section \ref{Sec_profiles}. The proof of Theorem \ref{Thm_existence_causality} is given in Section \ref{Sec_Cauchy}, and we prove our global existence result of Theorem \ref{Thm_global} in Section \ref{Sec_gl-ex}. In the appendix we compute the Newtonian limit of \eqref{rel_art_visc} and give an extension to General Relativity.

\section{Lorentz Invariance, Decay of FL Modes and the Newtonian Limit}  \label{Sec_Lorentz_FL-decay}

\subsection{Lorentz Invariance} \label{Sec_LorCauLim}

Lorentz invariance is a direct consequence of the invariance of $\partial_\nu T^{\mu\nu}$ and the wave operator. We prove this here for completeness.

\begin{Lemma} \label{Lemma_Lorentz-invariance}
The dissipative relativistic Euler equations \eqref{rel_art_visc} is Lorentz invariant.
\end{Lemma}

\Proof 
This follows directly from $u^\mu$ and $\partial_\nu T^{\mu\nu}$  both transforming as vectors, and $\Box$ being invariant under Lorentz transformations. In more detail, let $\Lambda^\mu_\alpha$ denote a Lorentz transformation, mapping a coordinate system $x^\mu$ to coordinates $\bar{x}^\alpha$ by $x^\mu = \Lambda^\mu_\alpha \bar{x}^\alpha + a^\mu$, where $a^\mu \in \R^4$ is the vector of a constant spacetime translation, and $\Lambda: \R^4 \to \R^4$ is a constant  $(4\times 4)$-matrix which preserves the Minkowski metric $(\eta_{\mu\nu}) = (\eta_{\alpha\beta}) = {\rm diag}(-1,1,1,1)$ under tensor transformation,
\beq \label{Lorentz_transfo_general}
\Lambda^{\mu}_{\alpha} \Lambda^{\nu}_{\beta} \eta_{\mu\nu} = \eta_{\alpha\beta}.
\eeq 
Now, the fluid density and pressure transform as scalars, e.g., $\rho(x) = \rho(x(\bar{x})) \equiv \bar{\rho}(\bar{x})$, while the fluid $4$-velocity transforms as a vector, $u^\mu = \Lambda^\mu_\alpha \bar{u}^\alpha$. It follows that the energy-momentum tensor \eqref{perfect_fluid} of a fluid transforms as $T^{\mu\nu} = \Lambda^\mu_\alpha \Lambda^\nu_\beta \bar{T}^{\alpha\beta}$, and thus 
\beq \label{Lorentz_transfo_fluid-div}
\partial_\nu T^{\mu\nu} 
= (\Lambda^{-1})_\nu^\gamma \partial_\gamma \big(\Lambda^\mu_\alpha \Lambda^\nu_\beta \bar{T}^{\alpha\beta} \big) 
= \Lambda^\mu_\alpha \partial_\beta  \bar{T}^{\alpha\beta} 
\eeq
transforms as a vector. On the other hand, since
\beq \label{Lorentz}
\Box = \eta^{\nu\sigma} \partial_\nu \partial_\sigma 
=  \eta^{\nu\sigma} (\Lambda^{-1})^\beta_\nu (\Lambda^{-1})^\gamma_\sigma \partial_\beta \partial_\gamma 
= \eta^{\beta\gamma} \partial_\beta \partial_\gamma = \bar{\Box}
\eeq
transforms as a scalar, $\Box u^\mu = \Lambda^\mu_\alpha \bar{\Box} \bar{u}^\alpha$ transforms as a vector. Combining this with \eqref{Lorentz_transfo_fluid-div}, we conclude that \eqref{rel_art_visc} transforms as
\beq
0= \partial_\nu T^{\mu\nu} - \visc \Box u^\mu = \Lambda^\mu_\alpha \Big(  \partial_\beta  \bar{T}^{\alpha\beta}  - \visc \bar{\Box} \bar{u}^\alpha \Big),
\eeq
which is the sought-after Lorentz invariance.  
\QED

\subsection{Decay of FL-modes - Proof of Proposition \ref{Thm_decay}} \label{Sec_diss}

We now prove existence of Fourier-Laplace (FL-) mode solutions \eqref{FL-modes} of form $\Psi e^{\lambda t + i \vec{\xi} \vec{x}}$ with $\lambda \in \mathbb{C}$ and $\vec{\xi} \in \R^3$ to the linearized dissipative relativistic Euler equations \eqref{linearized_RAV}, and prove decay in the sense that ${\rm Re}(\lambda)<0$. For this we first derive the linearization \eqref{linearized_RAV} of \eqref{rel_art_visc} at a steady state solution $(\rho_0,u^\mu_0)$ on a formal level, assuming sufficiently smooth solutions.

\begin{Lemma} \label{Lemma_linearization}
Assume a constant density $\rho_0 >0$ and fluid $4$-velocity $u^\mu_0 = (1,0,0,0)$. Assume $(\bar{\rho},\bar{u})$ is a smooth solution of \eqref{rel_art_visc} of form $\bar{\rho} = \rho_0 + \epsilon \rho$ and $\bar{u}^\mu \equiv u_0^\mu + \epsilon u^\mu$ with $u^\mu = (0,\vec{u})$ and $\vec{u} \equiv (u^1,u^2,u^3)$, where $\epsilon >0$ is constant. Then, for $\epsilon >0$ sufficiently small, omitting terms of order $O(\epsilon^2)$, $(\rho,u^\mu)$ solves the following linearization of \eqref{rel_art_visc}, 
\beq \label{linearized_RAV_proof}
\begin{cases}
\partial_t \rho + \ka\; \nabla \cdot \vec{u}  = 0 \cr
\ka\; \partial_t u_j + \sos \; \partial_j \rho = \visc\; \Box u_j ,
\end{cases}   
\eeq
for $j=1,..,3$, and where $\nabla=(\partial_{x}, \partial_{y},\partial_{z})$, and $\visc >0$ denotes viscosity. 
\end{Lemma}

\Proof
To begin, set $u^\mu = (0,u^1,u^2,u^3)$, which implies the normalization $\bar{u}^\sigma \bar{u}_\sigma =-1$ to hold within order $\epsilon^2$-errors,
\beq \label{App_eqn1}
\bar{u}^\sigma \bar{u}_\sigma 
=  u_0^\sigma (u_0)_\sigma  + \epsilon^2 u^\sigma u_\sigma  
=  -1 + O(\epsilon^2).
\eeq
By Taylor expansion in $\epsilon$, we find for the barotropic pressure law that
\beq  \label{App_eqn2}
p(\bar{\rho}) = p(\rho_0) + \epsilon \sos \rho + O(\epsilon^2),
\eeq
where $\sos = p'(\rho_0)$.              Thus, the energy momentum tensor $\bar{T}^{\mu\nu}$ of $(\bar{\rho},\bar{u})$, up to order $\epsilon^2$-errors, is given by
\beq \label{App_eqn3}
\bar{T}^{\mu\nu} = T^{\mu\nu}_0  +  \epsilon \Big( \ka ( u_0^\mu u^\nu + u_0^\nu u^\mu )  + (\sos + 1) u_0^\mu u_0^\nu \rho  + \sos  \eta^{\mu\nu} \rho \Big)   + O(\epsilon^2),
\eeq
where $T^{\mu\nu}_0$ denotes the energy momentum tensor of $(\rho_0,u^\mu_0)$ and $\ka = p(\rho_0) + \rho_0$. Taking now the divergence of \eqref{App_eqn3} gives
\beq \label{App_eqn4}
\partial_\nu \bar{T}^{\mu\nu} =  \epsilon \Big( \ka ( u_0^\mu \partial_\nu u^\nu + u_0^\nu \partial_\nu u^\mu )  + (\sos + 1)u_0^\mu u_0^\nu \partial_\nu \rho + \sos  \eta^{\mu\nu} \partial_\nu \rho \Big)   + O(\epsilon^2).
\eeq
Now, substituting \eqref{App_eqn4} for the left hand side in \eqref{rel_art_visc}, dividing by $\epsilon$ and omitting terms of order $\epsilon$, yields the linearization
\beq \label{App_eqn5}
 \ka ( u_0^\mu \partial_\nu u^\nu + u_0^\nu \partial_\nu u^\mu )  + (\sos + 1)u_0^\mu u_0^\nu \partial_\nu \rho + \sos  \eta^{\mu\nu} \partial_\nu \rho  =  \visc \Box u^\mu .
\eeq
The sought-after linearized equations \eqref{linearized_RAV} follow by substitution of $u^\mu_0 = (1,0,0,0)$ and $u^\mu = (0,u^1,u^2,u^3)$ in \eqref{App_eqn5}.
\QED

We now write \eqref{linearized_RAV} equivalently as a matrix equation on FL-modes. For this, substitute the FL-modes $\rho = \Psi^0 e^{\lambda t + i \vec{\xi} \vec{x}}$ and $u^j = \Psi^j e^{\lambda t + i \vec{\xi} \vec{x}}$ into \eqref{linearized_RAV_proof}. It follows then that the linearized dissipative relativistic Euler equations \eqref{linearized_RAV} are equivalent to the matrix equation
\beq \label{FL-modes_matrix-eqn}
M \cdot \Psi =0 , \hspace{1cm} \text{for}    \hspace{1cm}
M \equiv \begin{pmatrix}  \lambda   &   i \kappa \vec{\xi}^T \\  i \sos \vec{\xi}   &   (\kappa \lambda + \visc \lambda^2 + \visc |\vec{\xi}|^2) I_3
 \end{pmatrix},
\eeq
where $\sos \equiv p'(\rho_0)$, $\ka \equiv p(\rho_0) + \rho_0$, $\rho_0>0$  constant and $I_3$ denotes the identity matrix on $\R^3$. In our next lemma, we compute the dispersion relation $\det(M)=0$ associated to \eqref{FL-modes_matrix-eqn}, which gives a necessary and sufficient condition for the existence of FL-mode solutions. For this we now make without loss of generality the simplifying assumption that $\vec{\xi} =(\xi,0,0)$ for some $\xi \in \R$, (this can always be arranged for by a spatial rotation as shown in the end of this section).

\begin{Lemma}
Assume $\vec{\xi} =(\xi,0,0)$. A FL-mode \eqref{FL-modes} solves the linearized dissipative relativistic Euler equations if and only if $\lambda$ and $\xi$ satisfy the dispersion relation 
\beq \label{dispersion_relation}
 \big( \lambda^3 + \kati \lambda^2 + \xi^2 \lambda + \kati \sos \xi^2 \big) 
\big( \lambda^2 + \kati \lambda + \xi^2 \big) =0 ,
\eeq
where $\kati = \frac{\ka}{\visc}$ and $\visc>0$.
\end{Lemma}

\Proof
Substituting $\vec{\xi} = (\xi, 0, 0 )$ into \eqref{FL-modes_matrix-eqn} reduces the matrix $M$ to block diagonal form $M = {\rm diag}(M_1,M_2)$ with
\beq \label{FL-modes_matrix-eqn_special}
M_1 = \begin{pmatrix}  \lambda   &   i \kappa {\xi} \\  i \sos {\xi}   &   \kappa \lambda + \visc \lambda^2 + \visc\xi^2  \end{pmatrix}
\hspace{0.5cm} \text{and} \hspace{0.5cm} 
M_2 = (\kappa \lambda + \visc \lambda^2 + \visc\xi^2) I_2  .
\eeq
Thus, since $\det(M) = \det(M_1) \det(M_2)$,  the dispersion relation $\det(M)=0$ reduces to the sought-after equation \eqref{dispersion_relation}.
\QED

We now keep $\xi \in \R$ fixed and view each factor in \eqref{dispersion_relation} as a polynomial is $\lambda$. We next show that their non-zero roots all satisfy ${\rm Re}(\lambda) <0$, which implies the sought-after decay of FL-modes. Note, the case $\lambda =0$ implies the trivial FL-mode $(\lambda,\xi)=0$ of a steady state, which we can ignore here.

\begin{Lemma} \label{Lemma_decay}
For any $\xi \in \R$, there always exist some $\lambda \in \mathbb{C}$ which solves the dispersion relation \eqref{dispersion_relation}, and if $\xi \neq 0$, then $\lambda \neq 0$. Moreover, any non-zero solution $\lambda$ of \eqref{dispersion_relation} satisfies ${\rm Re}(\lambda) <0$.
\end{Lemma}

\Proof
Note first that, if $\xi=0$, then \eqref{dispersion_relation} implies either $\lambda=0$ (the trivial case of a steady state), or $\lambda = -\kati <0$. We now study the roots of the two polynomials in \eqref{dispersion_relation} separately for the non-trivial case $\xi \neq 0$. For this, set $\lambda = a + ib$, for $a,b \in \R$. 

The first condition \eqref{dispersion_relation} gives rise to, i.e., $\big( \lambda^3 + \kati \lambda^2 + \xi^2 \lambda + \kati \sos \xi^2 \big) = 0$, is equivalent to 
\begin{eqnarray} 
(a^3 -3ab^2) +\kati(a^2 - b^2) + \xi^2 a + \kati \sigma^2 \xi^2 &=& 0  \label{diss_techeqn1} \\
(3a^2b - b^3) + 2\kati a b + b \xi^2 &= &0 .   \label{diss_techeqn2}
\end{eqnarray}
Equation \eqref{diss_techeqn2} holds if either $b=0$ or $b^2=2\kati a + 3a^2 + \xi^2$. In case $b=0$, $a \equiv {\rm Re}(\lambda)$ satisfies 
\beq \label{diss_techeqn3}
a^3 + \kati a^2 + \xi^2 a +  \kati \sigma^2 \xi^2 = 0,
\eeq
which directly implies that $a < 0$, since all coefficients in \eqref{diss_techeqn3} are strictly positive, keeping in mind that $\kati > 0$, $\xi^2 > 0$ and $0<\sigma<1$. In case that $b^2=2\kati a + 3a^2 + \xi^2$, $a\equiv {\rm Re}(\lambda)$ satisfies
\beq \label{diss_techeqn4}
a^3 +  \kati a^2 + \frac{1}{4}(\kati^2 + \xi^2) a + \frac{\kati \xi^2}{8} (1-\sigma^2) = 0 ,
\eeq
which again implies that $a<0$, since again all coefficients are strictly positive.

From the second condition in \eqref{dispersion_relation}, i.e., $\lambda^2 + \kati \lambda + \xi^2 =0$, one can directly compute 
$$
\lambda = - \frac{\kati}{2} \pm \frac{1}{2} \sqrt{\kati^2 - 4 \xi^2}
$$
which implies again ${\rm Re}(\lambda) < 0$. This completes the proof.
\QED

To complete the proof of Proposition \ref{Thm_decay} we now show that a general $\vec{\xi}' = (\xi_1,\xi_2,\xi_3)$ of an FL-mode in a frame $\vec{x}'$ can always be reduced to $\vec{\xi}=(\xi,0,0)$ by a spatial rotation without changing $\lambda$. That is, choosing a $3\times 3$ rotation matrix $R$ such that $R^T \vec{\xi}' = \vec{\xi} \equiv (\xi,0,0)$, the exponent of the FL-modes \eqref{FL-modes} transforms under the spatial rotation $\vec{x}' = R \vec{x}$ as
\beq \nonumber
\lambda t + i (\vec{\xi}')^T \vec{x}' 
= \lambda t + i (\vec{\xi}')^T R \vec{x} 
= \lambda t + i (R^T \vec{\xi}')^T \vec{x} 
= \lambda t + i \vec{\xi}^T \vec{x}
= \lambda t + i \xi x.
\eeq	
Thus, $\lambda$ is the same in both frames $\vec{x}$ and $\vec{x}'$, so that Lemma \ref{Lemma_decay} readily implies the sought after decay of the FL-modes $(\lambda,\vec{\xi}')$ and $(\lambda,\vec{\xi})$.  Note finally that the dissipative relativistic Euler equations are invariant under spatial rotations, which map solutions of \eqref{linearized_RAV_proof} again to solutions. This completes the proof of Proposition \ref{Thm_decay} (ii).

\subsection{Newtonian limit}

We now compute the Newtonian limit (taking the speed of light $c\to \infty$) of the dissipative relativistic Euler equations.    In order to take the limit $c \to \infty$ in a formal sense, we first reinsert the speed of light $c$ into \eqref{rel_art_visc}, which in coordinates $x^\mu = (ct,x,y,z)$ gives 
\begin{equation} \label{Limit_rel_art_visc}
\partial_{\nu}T^{\mu \nu} = \visc\; \square u^{\mu}
\end{equation}
for 
$$
T^{\mu\nu} \equiv  (\rho + \frac{p}{c^2}) u^\mu u^\nu + p \eta^{\mu\nu} 
\ \ \ \ \ \text{and} \ \ \ \ \  
u^\mu  \equiv \gamma(\vec{v})\left( \begin{array}{c} c \cr \vec{v}  \end{array} \right),
$$ 
with $\gamma(v) = (1-\vec{v}^2/c^2)^{-\frac{1}{2}}$, Minkowski metric $\eta = {\rm diag}(-1,1,1,1)$, D'Alembert operator $\square \equiv  \eta^{\mu\nu} \partial_\mu \partial_\nu = - \tfrac{1}{c^2}\partial_{tt}^2 + \partial_{xx}^2 +\partial_{yy}^2 + \partial_{zz}^2$, classical fluid velocity $\vec{v}$, and fluid $4$-velocity $u^\mu$ satisfying $u^\mu u_\mu = -c^2$. 

Contracting \eqref{Limit_rel_art_visc} with $u_\mu$ yields the relativistic conservation of mass-energy equation,
\beq \label{Limit_eqn1}
u^\mu \partial_\mu \rho + \left(  \rho + \frac{p}{c^2} \right) \partial_\mu u^\mu = -\frac{\visc}{c^2} \; u_\mu \Box u^\mu ,
\eeq
where we used that $u_\mu \, \partial_\nu u^\mu= 0$ as a result of the normalization $u^\mu u_\mu = -c^2$. The relativistic conservation of 3-momentum equation is obtained by contraction of \eqref{Limit_rel_art_visc} with   $\Pi^{\mu\nu} \equiv \eta^{\mu\nu} +  c^{-2} u^\mu u^\nu$, (the projection orthogonal to $u^\mu$), as 
\beq \label{Limit_eqn2}
\left(  \rho + \frac{p}{c^2} \right) u^\nu \partial_\nu u^\mu + \Pi^{\mu\nu} \partial_\nu p =  \visc \Pi^{\mu\nu} \Box u_\nu ,
\eeq
using that $\Pi_{\mu\nu} u^\nu = 0$ and, (by $u_\nu \, \partial_\sigma u^\nu= 0$), that $\Pi^{\mu\nu}  \partial_\sigma u_\nu = \partial_\sigma  u^{\mu}$.

We now compute the limit $c\rightarrow \infty$ of \eqref{Limit_eqn1} and \eqref{Limit_eqn2}. First, setting $v^\mu \equiv u^\mu / \gamma(\vec{v})$, we find using the Leibnitz rule that
\beq \label{Limit_eqn3}
\Box u^\mu = \begin{cases} \frac{1}{2c} \Delta |v|^2 + O(c^{-3}) , \ \ \ \mu =0, \cr  \gamma(\vec{v}) \Box v^j + O(c^{-2}) , \ \ \ \mu =j , \end{cases}
\eeq
since $\gamma(\vec{v}) = O(1)$, $\partial_\nu \gamma(\vec{v}) = O(c^{-2})$ and $\partial_{\nu\mu}^2 \gamma(\vec{v}) = O(c^{-2})$. Thus, since $\Box = - \tfrac{1}{c^2}\partial_{tt}^2 + \Delta$ and $v^0 =c$, we obtain 
\beq \label{Limit_eqn4}
\begin{aligned}
\lim\limits_{c\to \infty} \Box u^j  &= \Delta v^j  \ \ \ \ \text{for} \ \mu =j \in \{1,2,3\}, \cr 
\lim\limits_{c\to \infty} u_\mu \Box u^\mu &= - |\nabla v|^2.
\end{aligned}
\eeq
To now compute the classical limit of \eqref{Limit_eqn1}, observe that
\beq \label{Limit_eqn5}
\partial_\sigma u^\sigma = \gamma(\vec{v}) \nabla \cdot \vec{v} + O(c^{-1})
\ \ \ \ \ \text{and} \ \ \ \  \
u^\sigma \partial_\sigma \rho = \gamma(\vec{v}) \left( \partial_{t} \rho + \vec{v} \cdot \nabla \rho  \right).
\eeq 
Substitution of \eqref{Limit_eqn5} into \eqref{Limit_eqn1} yields
\beq \nonumber
\gamma(\vec{v}) \left( \partial_{t} \rho + \vec{v} \cdot \nabla \rho  \right)
  + \left( \rho  + \frac{p}{c^2} \right) \left( \gamma(v) \nabla \cdot \vec{v} + O(c^{-1}) \right) = -\frac{\visc}{c^2} \; u_\sigma \Box u^\sigma ,
\eeq
and taking the limit $c\rightarrow \infty$ we obtain the sought-after conservation of mass equation,
\beq \label{Limit_eqn6}
 \partial_{t} \rho + \vec{v} \cdot \nabla \rho  +  \rho   \nabla \cdot \vec{v}  = 0,
\eeq
interpreting $\rho$ as classical mass density.
To compute the classical limit of \eqref{Limit_eqn2}, we note that 
\beq       \nonumber
\lim_{c\rightarrow \infty} u^\sigma  \partial_\sigma u^\nu  = \left( \begin{array}{c} 0 \cr  \partial_t \vec{v} + \vec{v} \cdot \nabla \vec{v} \end{array} \right)   
\ \ \ \ \ \text{and} \ \ \ \ 
\lim_{c\rightarrow \infty} \Pi_{\mu\nu} = \left( \begin{array}{cc} 0 & 0 \cr 0 & \text{id}_3  \end{array}  \right),
\eeq
where $\text{id}_3$ denotes the identity on $\R^3$. Substituting the above identities into \eqref{Limit_eqn2}, we obtain in the limit $c\to \infty$ the sought after conservation of momentum equation
\beq \label{Limit_eqn7}
 \rho  \partial_t \vec{v} + \rho \vec{v} \cdot \nabla \vec{v}  + \nabla p = \visc
 \; \Delta \vec{v} .
\eeq
This completes the proof of Proposition \ref{Thm_decay} (iii).

\section{Shock Profiles and Entropy Production} \label{Sec_profiles}

\subsection{Shock Profiles}
To construct traveling wave solutions it suffices to assume one dimensional fluid flow in the direction of the $x^1$-axis, that is, we assume $u^2=0=u^3$ and set $u \equiv u^1$, which implies $u^0 = \sqrt{1+u^2}$ by the normalization condition $u^\sigma u_\sigma=-1$. We assume further that the speed of sound $\sigma \equiv \sqrt{p'(\rho)}$ is bounded by the speed of light,  $0< \sigma < 1$.  Let $(\rho_L,u_L)$ and $(\rho_R,u_R)$ be  two constant states, (the states left and right of position of the shock discontinuity on the $x$-axis),  which satisfy the Rankine Hugoniot conditions \eqref{RH-condt}.   We make the traveling wave ansatz
\beq \label{traveling_wave}
\rho(t,x) \equiv \rho(\zeta)  
\hspace{.5cm} \text{and}  \hspace{.5cm}
u(t,x)= u(\zeta),
\eeq
for $\zeta \equiv \frac{x-st}{\visc}$ and $s\in \R$ fixed, and assume  
\begin{eqnarray}  
\lim\limits_{\zeta \to +\infty} u(\zeta) &=& u_R
\hspace{.5cm} \text{and}  \hspace{.5cm}
\lim\limits_{\zeta \to +\infty}  \rho(\zeta) = \rho_R,   \label{profiles_limits_1}   \\
\lim\limits_{\zeta \to -\infty}  u(\zeta) &=&  u_L
\hspace{.5cm} \text{and}  \hspace{.5cm}
\lim\limits_{\zeta \to -\infty}  \rho(\zeta) = \rho_L.    \label{profiles_limits_2}
\end{eqnarray}

Now, substituting the traveling wave ansatz \eqref{traveling_wave} into the dissipative relativistic Euler equations \eqref{rel_art_visc} yields a first order ODE. That is, a direct computation then shows that $\partial_\nu T^{\mu\nu}  = \tfrac{1}{\visc} \dot{T}^{\mu\nu} \zeta_\nu$ and $\Box u^\mu = \tfrac{1-s^2}{\visc^2} \ddot{u}^\mu(\zeta)$, where $\zeta_\nu = (-s,1)$ and a dot denotes differentiation with respect to $\zeta$. The dissipative relativistic Euler equations are thus equivalent to               
\beq  \label{ODE0}
\dot{T}^{\mu\nu} \zeta_\nu = (1-s^2 ) \ddot{u}^\mu(\zeta) ,
\eeq
and integration in $\zeta$ gives
\beq \label{profiles_techeqn1}
\dot{u}^\mu(\zeta) = \tfrac{1}{1-s^2} T^{\mu\nu} \zeta_\nu + \theta^\mu
\eeq
for integration constants $\theta^\mu = (\theta^0,\theta^1)$.  Imposing the limits \eqref{profiles_limits_1} in \eqref{profiles_techeqn1} determines the integration constants to be $\theta^\mu = -\tfrac{1}{1-s^2} T^{\mu\nu}_R \zeta_\nu$, since in particular $\lim\limits_{\zeta \to +\infty} \dot{u}(\zeta) =0$ by the first condition in \eqref{profiles_limits_1}.  Similarly, imposing \eqref{profiles_limits_2} together with $\lim\limits_{\zeta \to -\infty} \dot{u}(\zeta) =0$ in \eqref{profiles_techeqn1} implies, by the Rankine Hugoniot condition \eqref{RH-condt}, that $s$ is the shock speed, 
\beq \label{shockspeed}
s = \frac{[T^{01}]}{[T^{00}]} = \frac{[T^{11}]}{[T^{10}]} .
\eeq
We conclude that the dissipative relativistic Euler equations \eqref{rel_art_visc} for the traveling waves \eqref{traveling_wave} can be written equivalently as\footnote{In \eqref{profiles_RAV-Euler}, one can replace $T_R^{\mu\nu}$ by $T_L^{\mu\nu}$, since $T_R^{\mu\nu} \zeta_\nu = T_L^{\mu\nu} \zeta_\nu$ by the RH-condition \eqref{RH-condt}.} 
\beq \label{profiles_RAV-Euler}
\dot{u}^\mu(\zeta) = \tfrac{1}{1-s^2} \big( T^{\mu\nu} - T_R^{\mu\nu}\big) \zeta_\nu .
\eeq

We next show that, by the normalization condition $u^\sigma u_\sigma=-1$, equation \eqref{profiles_RAV-Euler} reduces to the scalar first order ODE \eqref{ODE} coupled to the explicit expression \eqref{density-function} for the density as a function of $u$ alone.       

\begin{Lemma}  \label{Lemma_ODE-decoupling} 
Let $\zeta_\nu = (-s,1)$, define the classical velocity $v(u) \equiv \frac{u}{\sqrt{1+u^2}}$, and assume $v(u) \neq s$. Then \eqref{profiles_RAV-Euler} is equivalent to     
\begin{eqnarray}
\rho(u) &=& \tfrac{1}{s-v(u)} \Big( v(u) T^{1\nu}_R - T^{0\nu}_R \Big) \zeta_\nu , \label{density-function_proof}   \\
\dot{u} &=& \tfrac{1}{1-s^2} \Phi(u),
\hspace{.5cm} \text{for} \hspace{.5cm}
\Phi(u) \equiv  \Big( T^{1\nu} - T^{1\nu}_R \Big) \zeta_\nu .  \label{ODE_proof} 
\end{eqnarray}
\end{Lemma}

\Proof
Substituting $u^0 = \sqrt{1+u^2}$ and its derivative $\dot{u}^0= v(u) \dot{u}$ into the first component of \eqref{profiles_RAV-Euler}, gives
\beq \label{ODE-decoupling_techeqn1}
v(u) \dot{u} = \tfrac{1}{1-s^2} \Big( T^{0\nu} - T^{0\nu}_R \Big) \zeta_\nu.
\eeq
Substitution of the second component of \eqref{profiles_RAV-Euler} for $\dot{u}$ in \eqref{ODE-decoupling_techeqn1} yields
\beq
v(u)  \Big( T^{1\nu} - T^{1\nu}_R \Big) \zeta_\nu  
=  \Big( T^{0\nu} - T^{0\nu}_R \Big) \zeta_\nu,
\eeq
which we write equivalently as
\beq  \label{ODE-decoupling_techeqn2}
\Big( v(u) T^{1\nu} - T^{0\nu} \Big) \zeta_\nu 
= \Big( v(u) T^{1\nu}_R - T^{0\nu}_R \Big) \zeta_\nu . 
\eeq
We now show that the expression on the left hand side of \eqref{ODE-decoupling_techeqn2} is equal to $\rho(s-v(u))$. For this, using $\zeta_\nu =(-s,1)$, observe that
\begin{eqnarray} 
\Big( v T^{1\nu} - T^{0\nu} \Big) \zeta_\nu 
&=& v T^{11} - (1+sv) T^{01} + s T^{00}    .
\end{eqnarray}
Substituting $T^{11} = (p+ \rho)u^2 +p$, $T^{01} = (p+ \rho)u \sqrt{1+u^2}$ and $T^{00} = (p+ \rho)(1+u^2) -p$, and separating terms multiplied by $(p+\rho)$, we find that 
\begin{eqnarray} \label{ODE-decoupling_techeqn3}
\Big( v T^{1\nu} - T^{0\nu} \Big) \zeta_\nu 
&=& -(s-v)p + (p+\rho) \Psi      ,
\end{eqnarray}
where 
\begin{eqnarray}  \label{ODE-decoupling_techeqn4}
\Psi   &\equiv &   v u^2 - (1+sv) u \sqrt{1+u^2} + s (1+u^2)  \cr
&=&  v u^2 -  u \sqrt{1+u^2}   + s  \cr
&=&    -v + s,
\end{eqnarray}
as can be verified by direct computation using $v = \frac{u}{\sqrt{1+u^2}}$. Thus, substitution of \eqref{ODE-decoupling_techeqn4} back into \eqref{ODE-decoupling_techeqn3} gives us as claimed
\beq  \label{ODE-decoupling_techeqn5}
\Big( v T^{1\nu} - T^{0\nu} \Big) \zeta_\nu = \rho (s-v).
\eeq
Substituting  \eqref{ODE-decoupling_techeqn5} into \eqref{ODE-decoupling_techeqn2}, and solving for $\rho$, gives us the sought-after expression \eqref{density-function_proof}.  Now, since $\rho$ is a function of $u$ alone, the second component of \eqref{profiles_RAV-Euler} turns directly into the sought-after scalar ODE \eqref{ODE_proof}. This completes the proof of Lemma \ref{Lemma_ODE-decoupling}.
\QED

We first derive the Lax conditions \eqref{Lax-condition_1} - \eqref{Lax-condition_2} in the case that $s=0$, and then obtain the Lax condition for general $s$ by Lorentz transformation. For this, note frist that the Rankine Hugoniot condition for $s=0$ imply $[T^{01}] =0$ for $T^{01} \equiv (p+\rho)u \sqrt{1+u^2}$, from which we infer that $u_L$ and $u_R$ have the same sign and thus are both non-zero whenever a shock discontinuity is present; the assumption $v(u)\neq s$ of Lemma \ref{Lemma_ODE-decoupling} is thereby met. Now, since \eqref{ODE_proof} is a scalar ODE in $u$, the analysis of shock profiles is reduced to a one dimensional fixed point problem for \eqref{ODE_proof}. We now show that \eqref{ODE_proof} has exactly two fixed points for $u$ in between $u_L$ and $u_R$, namely $u=u_L$ and $u=u_R$.

\begin{Lemma}  \label{Lemma_fixedpoints}
For $u$ in between $u_L$ and $u_R$, $\Phi(u)=0$ if and only if either $u=u_L$ or $u=u_R$. 
\end{Lemma}

\Proof
Assuming $s=0$, we find that 
\beq \label{fixedpoints_eqn1}
\Phi(u) =  \Big( T^{1\nu} - T^{1\nu}_R \Big) \zeta_\nu  = T^{11} - T^{11}_R ,
\eeq 
where $T^{11} = (p+\rho) u^2 + p$. From this we directly find that $\Phi(u)=0$ if $T^{11} = T^{11}_R$ or $T^{11}=T^{11}_L$ by the Rankine Hugoniot condition for $s=0$. Now, a direct computation shows that the density function \eqref{density-function_proof} satisfies $\rho(u_L) = \rho_L$ and $\rho(u_R) = \rho_R$. This implies that $\Phi(u_L) =0$ and $\Phi(u_R)=0$.

We next prove that no more fixed points exists between $u_L$ and $u_R$, by studying critical points. For this, we first express $\Phi$ in terms of the classical fluid 3-velocity $v\equiv \frac{u}{\sqrt{1+u^2}}$. For this we substitute the equivalent expression $u^2 = \frac{v^2}{1-v^2}$ into $T^{11} = (p+\rho) u^2 + p$ in \eqref{fixedpoints_eqn1}, which gives us
 \beq \label{fixedpoints_eqn2}
\Phi(v) = \frac{\Psi(v)}{1-v^2}, 
\eeq 
for 
\begin{eqnarray}
\Psi(v) &\equiv &  (p+\rho)v^2 + (p-T^{11}_R) (1-v^2)  \cr
&=& p(v) + v T^{01}_R - T^{11}_R,
\end{eqnarray}
where we used in the last equality that $\rho(v) = - T^{11}_R + \frac{1}{v} T^{01}_R$ by \eqref{density-function_proof}. Since $0 < v^2 <1$ by the speed of light bound, it follows that the roots of $\Phi(v)$ and $\Psi(v)$ are identical. We now show that $\Psi(v)$ has only one critical point, and thus at most two roots. A direct computation gives us
\beq \label{fixedpoints_eqn3}
\Psi'(v) = p'(v) + T^{01}_R = \frac{1}{v^2} (v^2 - \sigma^2) T^{01}_R,
\eeq
since $p'(v) = \sigma^2 \rho'(v)$ for $\sigma^2 \equiv \frac{dp}{d\rho}>0$, and since $\rho'(v)= - \frac{T^{01}_R}{v^2}$. From \eqref{fixedpoints_eqn3}, since $\frac{d^2p}{d\rho^2} \geq 0$ by assumption, it follows that $\Psi(v)$ has only one critical point on the interval between $u_L$ and $u_R$, as can be shown by noting that $G(v)\equiv v^2 - \sigma^2(v)$ is monotone, (since by \eqref{density-function_proof} $G'(v) = 2v + \frac{T_R^{01}}{v^2} p''(\rho(v))$ which has the same sign as $v$). On the other hand, since the Rankine Hugoniot condition $T^{01}_L = T^{01}_R$ for $T^{01} = (p +\rho) u \sqrt{1+u^2}$ imply that $u_L$ and $u_R$ have the same sign, we conclude that $\Psi(v)$ has at most one critical point in between $v_L$ and $v_R$. In summary this shows that $\Phi(v)$ has exactly two fixed points, namely $v_L$ and $v_R$. This proves the lemma.
\QED

Now, since \eqref{ODE_proof} has only $u_L$ and $u_R$ as fixed points, it follows that any local solution with initial data lying in between $u_L$ and $u_R$ automatically extends to a global solution $u(\zeta)$ on $(-\infty,\infty)$ such that $\lim\limits_{\zeta \to -\infty} u(\zeta) = u_L$ and $\lim\limits_{\zeta \to \infty} u(\zeta) = u_R$, based on the parameter direction fixed in \eqref{profiles_limits_1} and \eqref{profiles_limits_2}. We conclude that $u_L$ is an unstable rest point, while $u_R$ is asymptotically stable. Since 
$$
\Psi''(v)=\frac{(T_R^{01})^2}{v^4}p''(\rho(v)) +\frac{2 T_R^{01}}{v^3}p'(\rho(v)) >0,
$$
the stability of $u_R$ and instability of $u_L$ implies the inequality       
\beq \label{Phi_ineq1}
\Phi'(u_R) \ < \ 0 \ < \ \Phi'(u_L).
\eeq
From \eqref{Phi_ineq1} we now derive the Lax admissibility conditions \eqref{Lax-condition_1} - \eqref{Lax-condition_2} and thereby complete the proof of Theorem \ref{Thm_profiles}.      We begin by computing $\Phi'(u)$ in \eqref{Phi_ineq1}.

\begin{Lemma} \label{Lemma_Phi'}
Assume $s=0$ and set $\sigma \equiv \sqrt{\frac{dp}{d\rho}}$, then  
\beq \label{Phi'}
\Phi'(u) =  \frac{v^2 - \sigma^2}{v^2} (p + \rho) u,
\eeq
when evaluated either at $(u,\rho)=(u_L,\rho_L)$ or $(u,\rho)=(u_R,\rho_R)$.
\end{Lemma}

\Proof
For $s=0$, \eqref{density-function_proof} reduces to $\rho(u) = - T^{11}_R + \tfrac{1}{v(u)} T^{01}_R$, from which we get
\beq \label{Phi'_techeqn1}
\rho'(u) = -\tfrac{v'(u)}{v^2(u)} T^{01}_R  =  - \frac{1}{u^2 \sqrt{1+u^2}} T^{01}_R.
\eeq 
From \eqref{ODE_proof}, we find for $s=0$ that $\Phi(u) \equiv T^{11} - T^{11}_R$. From $T^{11} = (p+\rho)u^2 + p$ and \eqref{Phi'_techeqn1} we obtain
\begin{eqnarray} \label{Phi'_techeqn2}
\Phi'(u) &=& \frac{d T^{11}  }{du} 
\ = \  2(p+\rho) u  + \big(u^2 + u^2 \sigma^2 + \sigma^2 \big) \rho'(u)   \cr
&=&  \frac{1}{\sqrt{1+u^2}}  \Big( 2 T^{01} -  \frac{(1+ \sigma^2)u^2 + \sigma^2}{u^2} T_R^{01}  \Big),
\end{eqnarray}
recalling that $T^{01}= (p+\rho)u \sqrt{1+u^2}$ and $\sigma^2 = \frac{dp}{d\rho}$. Thus evaluation at $u_R$ gives
\begin{eqnarray}  \label{Phi'_techeqn3}
\Phi'(u_R) & = &  T^{01}_R \; \frac{u^2 - (1+ u^2) \sigma^2}{u^2\sqrt{1+u^2}}\bigg|_{u_R}   
= (p+\rho) \; \frac{u^2 - (1+ u^2) \sigma^2}{u}  \bigg|_{u_R} \cr
&=&  (p+\rho)\;\frac{v^2 - \sigma^2}{v \sqrt{1-v^2}}  \bigg|_{u_R} 
=  (p + \rho) u \;\frac{v^2 - \sigma^2}{v^2} \bigg|_{u_R},
\end{eqnarray}
which is the sought-after expression \eqref{Phi'} at $u=u_R$, while \eqref{Phi'} at $u=u_L$ follows by \eqref{Phi'_techeqn3}, since $T^{01}_R=T_L^{01}$ in light of the Rankine Hugoniot conditions. This completes the proof of Lemma \ref{Lemma_Phi'}.
\QED

Now, by \eqref{Phi'}, we write \eqref{Phi_ineq1} equivalently as
\beq \label{Phi_ineq2}
\frac{v_R^2 - \sigma_R^2}{v_R^2} (p_R + \rho_R) u_R
 \ < \ 0 \ < \ 
\frac{v_L^2 - \sigma_L^2}{v_L^2} (p_L + \rho_L) u_L ,
\eeq
and since $v^2>0$ and $(p+\rho)_{L/R}>0$, condition \eqref{Phi_ineq2} is equivalent to
\beq \label{Phi_ineq3}
(v_R^2 - \sigma_R^2)  u_R \ < \ 0 \ < \ (v_L^2 - \sigma_L^2) u_L .
\eeq
Now, by the Rankine Hugoniot condition for $s=0$, we have $[T^{01}] =0$ for $T^{01} \equiv (p+\rho)u \sqrt{1+u^2}$, which, since $(p+\rho)_{L/R}>0$ and $\sqrt{1+u^2}>0$, implies that $u_L$ and $u_R$ always have the same sign. Thus, considering first the case that $u_L, u_R >0$, (corresponding to 1-shocks), it follows that $v_L+\sigma_L>0$ and $v_R + \sigma_R>0$, and \eqref{Phi_ineq3} is equivalent to 
\beq \label{Phi_ineq4}
v_R - \sigma_R  \ < \ 0 \ < \ v_L - \sigma_L 
\hspace{1cm} \text{and} \hspace{1cm}
0 < v_R + \sigma_R.
\eeq 
Now, since $1 \pm v \sigma >0$ by the speed of light bounds $|v| < 1$ and $0\leq \sigma <1$, dividing \eqref{Phi_ineq4} by $1 \mp v\sigma$ yields directly the Lax condition for 1-shocks,
\begin{eqnarray} \label{Lax_condt_1}
\lambda_1(u_R) < 0 < \lambda_1(u_L) 	
\hspace{1cm} \text{and} \hspace{1cm}
0 < \lambda_2(u_R), 
\end{eqnarray}
where the characteristic speeds (as computed in \cite{SmollerTemple}) are
\beq \label{char_speeds}
\lambda_1 = \frac{v-\sigma}{1 - v \sigma} 
\hspace{.5cm} \text{and} \hspace{.5cm}
\lambda_2 = \frac{v+\sigma}{1 + v \sigma} .
\eeq
Similarly, turning now to the case of 2-shocks by assuming $u_L, u_R <0$, it follows that $v_L-\sigma_L<0$ and $v_R - \sigma_R<0$. Thus \eqref{Phi_ineq3} is equivalent to 
\beq \label{Phi_ineq5}
v_R + \sigma_R  \ < \ 0 \ < \ v_L + \sigma_L 
\hspace{1cm} \text{and} \hspace{1cm}
v_L - \sigma_L < 0,
\eeq 
and dividing \eqref{Phi_ineq5} by $1 \pm v \sigma$ gives the Lax condition for the 2-shocks,
\beq \label{Lax_condt_2}
\lambda_2(u_R) < 0 < \lambda_2(u_L)
\hspace{1cm} \text{and} \hspace{1cm}
\lambda_1(u_L) < 0.
\eeq
The Lax conditions for general $s\neq 0$ now follows by a Lorentz boost: 

\begin{Lemma} \label{Lemma_Lax-condt}
Under a Lorentz boost to a frame in which the 1-shock moves with velocity $s_1$, the Lax condition \eqref{Lax_condt_1} turns into 
\begin{eqnarray} \label{Lax_condt_1_Lemma}
\lambda_1(u_R)  < s_1 < \lambda_1(u_L) \ \ \ \ \ \text{and}  \ \ \ \ \   s_1 < \lambda_2(u_R).   
\end{eqnarray}
Similarly, under a Lorentz boost to a frame in which the 2-shock moves with velocity $s_2$, the Lax condition \eqref{Lax_condt_2} turns into 
\begin{eqnarray} \label{Lax_condt_2_Lemma}
\lambda_2(u_R)  < s_2 < \lambda_2(u_L) \ \ \ \ \ \text{and}  \ \ \ \ \   \lambda_1(u_L) <  s_2.
\end{eqnarray}
\end{Lemma}

\Proof
Instead of Lorentz transforming \eqref{Lax_condt_1} and \eqref{Lax_condt_2} directly, we transform equations \eqref{Phi_ineq4} and \eqref{Phi_ineq5}. For this, by the velocity addition formula, observe that the fluid 3-velocities $v$ in the frame moving with the shock wave (in which the shock speed is $s=0$) is given by 
\beq \label{vel_addition}
v = \frac{\bar{v} - s}{1-s\bar{v}}
\eeq 
where $\bar{v}$ is the fluid 3-velocity in the frame in which the shock wave moves with the constant speed $s$. Note that, since the density and the pressure are Lorentz scalars, the speed of sound transforms as a scalar as well. Thus, equation \eqref{Phi_ineq4} transforms as
\beq \label{Phi_ineq4_transo} \nonumber
\frac{\bar{v}_R - s_1}{1-s_1\bar{v}_R} - \sigma_R  \ < \ 0 \ < \  \frac{\bar{v}_L - s_1}{1-s_1\bar{v}_L} - \sigma_L 
\hspace{1cm} \text{and} \hspace{1cm}
0 < \frac{\bar{v}_R - s_1}{1-s_1\bar{v}_R} + \sigma_R,
\eeq 
which can be written equivalently as \eqref{Lax_condt_1_Lemma}, denoting for simplicity $\bar{v}_{L/R}$ as $v_{L/R}$. In a similar fashion, replacing $v_{L/R}$ in \eqref{Phi_ineq5} by \eqref{vel_addition}, a direct computation 
leads to the Lax condition \eqref{Lax_condt_2_Lemma} for the 2-shock.
\QED

\subsection{Proof of Theorem \ref{Thm_profiles} - Part I}
Under the assumption that $(\rho_L,u_L)$ and $(\rho_R,u_R)$ are constant states which satisfy the RH conditions \eqref{RH-condt} with shock speed $s$, Theorem \ref{Thm_profiles} asserts that the following statements are equivalent: (i) There exists a traveling wave solution $\rho(t,x) = \rho(\zeta)$ and $u^\mu(t,x)= u^\mu(\zeta)$, $\zeta = \frac{x-st}{\visc}$, of the dissipative relativistic Euler equations \eqref{rel_art_visc} with asymptotics \eqref{asymptotics}. (ii) The first Lax condition \eqref{Lax-condition_1} holds with $s=s_1$ if $u_L$ and $u_R$ are separated by a 1-shock, and the second Lax condition \eqref{Lax-condition_2} holds with $s=s_2$ if $u_L$ and $u_R <0$ are separated by a 2-shock.                        

That (i) implies (ii) follows directly from the above constructions and Lemma \ref{Lemma_Lax-condt}. That is, by (i) it follows that $u^\mu(t,x)= u^\mu(\zeta)$ solves the ODE \eqref{ODE_proof} with $\rho(t,x) = \rho(u)$ given by \eqref{density-function_proof}, such that the right hand side of the ODE  \eqref{ODE_proof} satisfies the inequality \eqref{Phi_ineq1}. Lemma \ref{Lemma_Lax-condt} then implies the sought-after Lax conditions \eqref{Lax-condition_1} and \eqref{Lax-condition_2}.

Vice versa, that (ii) implies (i) follows since the Lax conditions \eqref{Lax-condition_1} and \eqref{Lax-condition_2} imply the inequality \eqref{Phi_ineq1}. Thus, since by Lemma \ref{Lemma_fixedpoints} $u_L$ and $u_R$ are the only rest points of \eqref{ODE_proof}, and since by \eqref{Phi_ineq1} $u_R$ is stable and $u_L$ is unstable, it follows that any local solution  with data in between $u_L$ and $u_R$ (which exist by the Picard Lindel\"off theorem) extends to a global solution with the sought-after asymptotics \eqref{asymptotics}. The sought after traveling wave solutions is then obtained by reversing the steps between equations \eqref{ODE0} to  \eqref{ODE_proof}, and substituting $\zeta = \frac{x-st}{\visc}$.  This completes the proof of Theorem \ref{Thm_profiles}.   
\hfill $\Box$

\subsection{Proof of Theorem \ref{Thm_profiles} - Part II}   \label{Sec_L2-approx}
We consider states $(\rho_L,u_L)$ and $(\rho_R,u_R)$ connected by a single Lax admissible shock curve. That is, let $\tilde{u}(t,x) = u_0(x-st)$ and $\tilde{\rho}(t,x) = \rho_0(x-st)$ be a Lax-admissible shock wave solution to \eqref{rel_Euler} subject to the normalization $u^\sigma u_\sigma =-1$, (cf.  \cite{SmollerTemple}), for Riemann initial data 
\begin{equation} \label{Riemann-problem}
u_0(x) =\begin{cases} 
      u_L, & x < x_0, \\
      u_R, & x  \geq x_0,
   \end{cases}
   \quad\mathrm{and}\quad 
\rho_0(x) =\begin{cases} 
      \rho_L, & x < x_0, \\
      \rho_R, & x  \geq x_0,
   \end{cases}
\end{equation}
$(t,x) \in [0,\infty) \times \R$, where $s\in \R$ is the shock speed and $x_0\in \R$ is the shock position which, in contrast to Theorem \ref{Thm_profiles}, we now need to specify.  We proved in the previous section that each such shock wave solution has a profile, i.e., a travelling wave solution 
\beq \label{L2-approx_eqn1}
\big(\rho_\visc(t,x), u_\visc(t,x) \big) = \big(\rho(\zeta), u(\zeta)\big),
\eeq  
where $(\rho(\zeta),u(\zeta))$ is a solution of the ODE system \eqref{density-function_proof} - \eqref{ODE_proof}, connecting the limit states $(\rho_L,u_L)$ and $(\rho_R,u_R)$ as $\zeta \to \pm \infty$.  However, these profiles are not unique,    until locked in with the shock position.      Namely, \eqref{L2-approx_eqn1} is a travelling wave solution of \eqref{rel_art_visc} for any parameter $\zeta \equiv \frac{x-st-k}{\visc}$, for any center $k \in \R$ independent of $\visc$.\footnote{This can be verified following the steps in the derivation of \eqref{density-function_proof} - \eqref{ODE_proof}, and solving the ODE \eqref{ODE_proof} by separation of variables in combination with the inverse function theorem. This show the general solution of the ODE \eqref{ODE} is of form $u(\zeta) = f(\zeta+k)$ for some function $f$.}     To establish Theorem \ref{Thm_profiles} (Part II), we prove in the following proposition that the traveling wave solution \eqref{L2-approx_eqn1} centred at the shock position $x_0$ (by setting $k=x_0$), is the unique such wave approximating the shock in the zero viscosity limit. 

\begin{Prop} \label{L2_approx_Prop}  
There exists a unique travelling wave solution to \eqref{rel_art_visc} of form \eqref{L2-approx_eqn1} with an $\visc$-independent center, namely the solution with parameterization $\zeta = \frac{x-st-x_0}{\visc}$, ($x_0 \in \R$ independent of $\visc$), such that
\beq \label{L2_approx_Prop_maineqn}
\lim_{\visc \to 0} \Big( \|u_\visc(t,\cdot) - \tilde{u}(t,\cdot)\|_{L^2(\R)} + \|\rho_\visc(t,\cdot) - \tilde{\rho}(t,\cdot)\|_{L^2(\R)} \Big) = 0 ,
\eeq 
uniformly in $t \in [0,\infty)$. 
\end{Prop} 

Theorem \ref{Thm_profiles} (Part II) is then a direct consequence of Proposition \ref{L2_approx_Prop}. To prove Proposition \ref{L2_approx_Prop}, recall that the density of a travelling wave solution is given by \eqref{density-function_proof} as a function of $u(\zeta)$, i.e., $\rho = \rho(u)$, while $u(\zeta)$ is a smooth solution of the first order ODE \eqref{ODE_proof}, $\dot{u}(\zeta)=\Phi_s\big(u(\zeta)\big)$ for $\Phi_s\equiv \frac{\Phi}{1-s^2}$. To prove Proposition \ref{L2_approx_Prop}, we first establish the following lemma, which implies existence of a travelling wave approximating the shock component $\tilde{u}$, the approximating density wave then follows from $\rho = \rho(u)$ being a function of $u(\zeta)$, by \eqref{density-function_proof}.

\begin{Lemma}\label{L2_approx_Lemma1}
The travelling wave $u_\visc(t,x) \equiv u(\tilde\zeta)$, for $\tilde\zeta \equiv \frac{x-st-x_0}{\visc}$ and $u(\cdot)$ the solution of \eqref{ODE_proof}, satisfies
\begin{equation} \label{L2_approx_Lemma1_eqn1}
\lim_{\visc \to 0}  \|u_\visc(t,\cdot) - \tilde{u}(t,\cdot)\|_{L^2(\R)} =0,
\end{equation}
uniformly in $t\in [0,\infty)$.      
\end{Lemma}

\begin{proof}
Writing the shock wave as
\begin{equation} \label{L2_approx_Lemma1_techeqn1}
\tilde u(x,t) = g(\tilde\zeta), \qquad \text{for} \quad  g(\tilde\zeta) \equiv u_L\,\chi_{\tilde\zeta\le 0}+u_R\,\chi_{\tilde\zeta\geq 0} ,
\end{equation}
and changing variables from $x$ to $\tilde\zeta \equiv \frac{x-st-x_0}{\visc}$, (so $dx=\visc\,d\tilde\zeta$), we obtain
\begin{align} \label{L2_approx_Lemma1_techeqn2}
\|u_\visc(\cdot,t)-\tilde{u}(\cdot,t)\|_{L^2(\mathbb{R})}^2
& = \visc \int_{\R}|u(\tilde\zeta)-g(\tilde\zeta)|^2\,d\tilde\zeta   \equiv \visc \mathcal{I}.
\end{align}
Since $ \mathcal{I}$ is independent of time $t$, the proof is complete once we show that
\begin{equation} \label{L2_approx_Lemma1_techeqn3}
 \mathcal{I} = \int_{-\infty}^{0} |u(\tilde\zeta)-u_L|^2\,d\tilde\zeta +\int_{0}^{\infty} |u(\tilde\zeta)-u_R|^2\,d\tilde\zeta < \infty .
\end{equation}
For this, we prove the following assertion. \vspace{.1cm}

\noindent \textit{Claim:} {\it There exists constants $K>0$, $c>0$ and $\zeta_\infty >0$ such that }
\beq \label{L2_approx_Lemma1_eqn2}
\begin{aligned}
|u(\zeta)-u_R| \leq \sqrt{K} e^{-c\zeta}, \qquad \forall\,\zeta \ge \zeta_\infty, \cr 
|u(\zeta)-u_L| \leq \sqrt{K} e^{c\zeta}, \qquad \forall\,\zeta \leq -\zeta_\infty.
\end{aligned}
\eeq 

\noindent \textit{Proof of Claim:} We only prove the first estimate in \eqref{L2_approx_Lemma1_eqn2}, the second one then follows by an analogous argument. Recall that the solution $u(\zeta)$ of the first order ODE \eqref{ODE_proof}, has asymptotics 
\begin{equation}\label{eq:profile-asymptotics}
\lim_{\zeta\to-\infty} u(\zeta)=u_L
\qquad \text{and } \qquad 
\lim_{\zeta\to+\infty} u(\zeta)=u_R,
\end{equation}
and recall that $\Phi(u_R)=0=\Phi(u_L)$ and $\Phi'(u_R) < 0 < \Phi'(u_L)$. Since $\Phi'(u_R)<0$ and since $\Phi'(\cdot)$ is continuous, there exist $\delta>0$ and
$c>0$ such that
\begin{equation}\label{eq:phi-prime-neg}
\Phi'(u)\le -c<0,
\qquad\text{whenever } |u-u_R|<\delta.
\end{equation}
Because $u(\zeta)\to u_R$ as $\zeta\to+\infty$, there exists $\zeta_\infty \in \R$ such that
\begin{equation}\label{eq:tail-stays-in-neighborhood}
|u(\zeta)-u_R| \le \delta ,
\qquad \forall\,\zeta \ge \zeta_\infty.
\end{equation}
Define $w(\zeta)\equiv u(\zeta)-u_R$. Then for $\zeta\ge\zeta_\infty$ we have $|w(\zeta)|\le \delta$ and
\begin{equation}\label{eq:wprime}
\frac{d}{d\zeta} w(\zeta)=\frac{d}{d\zeta}u(\zeta)=\Phi_s\!\big(u(\zeta)\big)=\Phi_s(u_R+w(\zeta)).
\end{equation}
By the mean value theorem, using $\Phi(u_R)=0$, we find that
\begin{equation}\label{eq:mvt}
\Phi(u_R+w)=\Phi(u_R)+\Phi'(u_R+\theta w)\,w
=\Phi'(u_R+\theta w)\,w
\end{equation}
for some $\theta=\theta(\zeta)\in(0,1)$. Hence, for $\zeta\ge\zeta_\infty$,
\begin{equation}\label{eq:wprime-linear}
\frac{d}{d\zeta}w(\zeta)=\Phi_s'(u_R+\theta w)\,w.
\end{equation}
Since $w(\zeta)\neq 0$ for all $\zeta <\infty$, we have
\begin{align}
\frac{d}{d\zeta}| w(\zeta) |
&=sgn(w(\zeta))\,\frac{d}{d\zeta}w(\zeta)
=\Phi_s'(u_R+\theta w)\,|w(\zeta)| .
\end{align}
By \eqref{eq:tail-stays-in-neighborhood}, we have
$|(u_R+\theta w)-u_R|\le |w|\le \delta$, so \eqref{eq:phi-prime-neg} yields
\begin{equation}\label{eq:v-ineq}
\frac{d}{d\zeta} |w(\zeta)|\le -c \,|w(\zeta)|,
\qquad \forall\,\zeta\ge\zeta_\infty.
\end{equation}
Integrating \eqref{eq:v-ineq} from $\zeta_\infty$ to $\zeta\ge\zeta_\infty$ gives
\[
|w(\zeta)|\le |w(\zeta_\infty)|\,e^{-c(\zeta-\zeta_\infty)},
\]
which establishes \eqref{L2_approx_Lemma1_eqn2} and proves the claim.  

Using the estimates \eqref{L2_approx_Lemma1_eqn2} in \eqref{L2_approx_Lemma1_techeqn3} directly shows that the integral for $|\zeta| > \zeta_\infty$ is finite, while finiteness over the compact set $|\zeta| \leq \zeta_\infty$ follows by continuity of the integrands. This completes the proof of Lemma \ref{L2_approx_Lemma1}. 
\end{proof}

\begin{Lemma}\label{L2_approx_Lemma2}
Assume $u_\visc(t,x) \equiv u(\tilde\zeta)$ as in Lemma \ref{L2_approx_Lemma1}, converging to $\tilde{u}$ as in \eqref{L2_approx_Lemma1_eqn1}. Then $\rho_{\visc}(t,x) \equiv \rho(u_{\visc}(t,x))$ with $\rho(\cdot)$ defined by \eqref{density-function_proof}, converges to the shock density $\tilde\rho$ in the sense that 
\beq \label{L2_approx_Lemma2_eqn1}
\lim_{\visc \to 0}\|\rho_\visc(t,\cdot) -\tilde\rho(t,\cdot)\|_{L^2(\R)} = 0,
\eeq  
uniformly in $t \in [0,\infty)$. 
\end{Lemma}

\Proof
From \eqref{density-function_proof}, we have that $\rho(u) = \tfrac{1}{s-v(u)} \Big( v(u) T^{1\nu}_R - T^{0\nu}_R \Big) \zeta_\nu$ is a $C^1$ function in $u$, provided that $v(u) \neq s$, where $v \equiv \frac{u}{\sqrt{1+u^2}}$.
In the frame where $s=0$, we show below Lemma \ref{Lemma_Phi'} that $u_L$ and $u_R$ have the same sign and, (since $u_{L}$ and $u_R$ are rest points of the ODE \eqref{ODE_proof}), that $u(\zeta)$ lies in the interval $[u_l,u_R]$ for all $\zeta \in \R$. This implies that $u(\zeta) \neq s$, and hence $v(u) \neq s$, for all $\zeta \in \R$, and by Lorentz transformation also in frames where $s\neq 0$.  
We conclude that $\rho(\cdot)$, defined by \eqref{density-function_proof}, is a $C^1$ function.
From this we obtain the Lipschitz estimate $|\rho(u_1) - \rho(u_2)| \leq L |u_1 - u_2|$, for $L\equiv \max_{u\in [u_L,u_R]} |\rho'(u)|$, which in turn yields the estimate 
$$
\|\rho(u_{\visc}(t,\cdot)) - \rho(\tilde{u}(t,\cdot))\|_{L^2(\R)} 
\leq L \|u_{\visc}(t,\cdot) - \tilde{u}(t,\cdot)\|_{L^2(\R)} ,
$$
which implies the sought-after convergence \eqref{L2_approx_Lemma2_eqn1}.
\QED

\Proof[Proof of Proposition \ref{L2_approx_Prop}:] 
Lemma \ref{L2_approx_Lemma1} together with Lemma \ref{L2_approx_Lemma2} imply the existence of a travelling wave solution approximating the shock wave $(\tilde\rho,\tilde{u})$ in the sense of \eqref{L2_approx_Prop_maineqn}, as asserted by Proposition \ref{L2_approx_Prop_maineqn}. That the travelling wave with parameterization $\tilde\zeta\equiv \frac{x-st + x_0}{\visc}$ is the only such travelling wave approximating the shock wave is a direct consequence of the next lemma, which once proven completes the proof of Proposition \ref{L2_approx_Prop}.
\QED

\begin{Lemma} \label{L2_approx_Lemma3}
Let $t_0\in[0,\infty)$ be fixed. 
For $k\in\R$, define the travelling wave
\begin{equation}\label{eq:visc-approx-k}
u_\visc^{k}(t,x)\equiv u\!\left(\frac{x-st-k}{\visc}\right),
\end{equation}
where as before $u$ is a solution of the ODE \eqref{ODE_proof} with asymptotics  \eqref{profiles_limits_1} - \eqref{profiles_limits_2}. Then for all $k\neq x_0$ there exists a constant $c_k>0$ such that
\begin{equation}\label{eq:lower-bound}
\|u_\visc^{k}(t_0,\cdot)-\tilde u(t_0,\cdot)\|_{L^2(\R)}^2 \ge c_k 
\end{equation}
for all $\visc>0$ sufficiently small. 
\end{Lemma}

\Proof
Let $t_0\in[0,\infty)$ be fixed, and set $\eta\equiv\frac{1}{8}|u_L-u_R|$. Assume without loss of generality $k>x_0$.  By the asymptotics \eqref{profiles_limits_1} - \eqref{profiles_limits_2} of $u(\zeta)$, there exists some $\zeta_\infty>0$ 
such that
\begin{equation}\label{eq:m-choice}
|u(\zeta)-u_L|< \eta,  \quad \forall \ \zeta \le -\zeta_\infty,
\end{equation}
and $|u(\zeta)-u_R|<\eta$ for all $\zeta\ge \zeta_\infty$. Consider now points $x$ left of the viscous center $s t_0+k$ by at least $\delta \equiv \frac{|k-x_0|}{4}$, that is, consider points $x \in (-\infty, st_0+k-\delta)$. Choosing $\visc \leq \frac{\delta}{\zeta_\infty}$, we obtain the bound $\frac{x-st_0-k}{\visc}\le - \frac{\delta}{\visc} \le -\zeta_\infty$ for all $x \in (-\infty, st_0+k-\delta)$. Thus, by \eqref{eq:m-choice}, we obtain for all $x \in (-\infty, st_0+k-\delta)$ that 
\begin{equation}\label{eq:uk-close-uL}
|u_\visc^{k}(t_0,x)-u_L|
=\left|u\left(\tfrac{x-st_0-k}{\visc}\right)-u_L\right|<\eta .
\end{equation}
Now, consider points restricted to the interval $J\equiv \left[st_0+x_0+\delta,\; st_0+k-\delta\right]$. Since $k=x_0 + 4 \delta$, $J$ has length $|J| = 2\delta$. Moreover, since any $x\in J$ satisfies $x>st_0+x_0$, the shock wave structure \eqref{Riemann-problem} implies $\tilde u(t_0,x)=u_R $ for all $x\in J$. On the other hand, since any $x\in J$ satisfies $x\le s t_0+k-\delta$, \eqref{eq:uk-close-uL} implies that  $|u_\visc^{\,k}(t_0,x)-u_L|<\eta$ for all $x\in J$. Thus, for all $x\in J$, we have
\begin{align*}
|u_\visc^{k}(t_0,x)-\tilde u(t_0,x)|
&\ge |u_R-u_L|-|u_\visc^{k}(t_0,x)-u_L|
\ge 7 \eta ,
\end{align*}
which directly implies
\begin{align*}
\|u_\visc^{k}(t_0,\cdot)-\tilde u(t_0,\cdot)\|_{L^2(\R)}^2
\geq \|u_\visc^{k}(t_0,\cdot)-\tilde u(t_0,\cdot)\|_{L^2(J)}^2
\geq 2 \eta^2 \delta \equiv c_k .
\end{align*}
This establishes the sought-after bound \eqref{eq:lower-bound} and proves the lemma. 
\QED

\subsection{Positivity of Entropy Production} \label{Sec_entropy}

Assume $(\rho,u^\mu)$ is a solution of the dissipative relativistic Euler equation \eqref{rel_art_visc} with $u^\sigma u_\sigma =-1$, and $\rho \in C^1([0,T]\times \R)$, $u^\mu \in C^2([0,T]\times \R)$ for some $T>0$. Assume the first law of thermodynamics in non-rest frame form       
\beq \label{first-law}
D\rho = \frac{\rho + p}{n} Dn + nT Ds,
\eeq 
where $s$ is the specific entropy, $T>0$ is the temperature, and $n>0$ is particle number density assumed to be conserved in the sense that $\partial_\nu(n u^\nu)=0$, and $D \equiv u^\nu \partial_\nu$ is the material derivative accounting for a shift away form the rest frame, (see Chapter 22 in \cite{MisnerThorneWheeler}).

A direct computation, starting with $T^{\mu\nu} = (p+\rho)u^{\mu}u^{\nu} + p \eta^{\mu \nu}$ and using that $u^\sigma u_\sigma =-1$ implies $u_\sigma \partial_\nu u^\sigma =0$, gives 
\begin{align} \label{thermodynamics-1}
 u_\mu \partial_{\nu}T^{\nu \mu} 
& = u_\mu \partial_{\nu} \Big( (\rho + p)u^\mu u^\nu \Big) +  u^\nu \partial_{\nu} p \cr 
& = - u^\nu \partial_{\nu}\rho - (\rho +p)\partial_{\nu}u^\nu .
\end{align}
Thus, substituting the first law \eqref{first-law}, we obtain the identity
\begin{align} \label{thermodynamics-2}
 u_\mu \partial_{\nu}T^{\nu \mu} 
& = - \Big( \frac{\rho + p}{n} u^\nu\partial_{\nu} n + nT  u^\nu\partial_{\nu} s \Big) - (\rho +p)\partial_{\nu}u^\nu  \cr 
& = - nT u^\nu \partial_{\nu} s   - \frac{\rho + p}{n} \partial_{\nu}(u^\nu n) \cr 
& = - nT u^\nu \partial_{\nu} s,
\end{align}
where the last equality follows by particle number conservation, $\partial_{\nu}(u^\nu n)=0$. We conclude that the dissipative relativistic Euler equation \eqref{rel_art_visc} implies that 
\begin{equation}\label{Entropy-rel}
 TnDs =  - \varepsilon u_{\mu}\Box u^{\mu} ,
\end{equation}
where $Ds \equiv u^\nu \partial_{\nu} s$ is the entropy production. We now obtain the following lemma, which directly implies Theorem \ref{Thm_profiles} (III).

\begin{Lemma} \label{Lemma_entropy-production}
Any traveling wave solution of form $u^{\mu}=u^{\mu}(\xi)$, $\rho=\rho(\xi)$ of \eqref{rel_art_visc} in one spatial dimension, for $\xi \equiv x-ct$, $c\in \R$, restricted to the region where $\frac{du^\mu}{d\xi} \neq 0$, has positive entropy production if and only if $|c|<1$. Moreover, it has zero entropy production if and only if $|c|=1$.
\end{Lemma}

\Proof
We first compute the right hand side of \eqref{Entropy-rel}. Using that $\Box(u_{\nu} u^{\nu}) = 2u_{\nu} \Box u^{\nu} + 2(\partial_{\alpha}u_{\nu})(\partial^{\alpha}u^{\nu})$ and $\Box(u_{\nu} u^{\nu})=0$ by the normalization condition, we find that in one spatial dimension
\begin{align}
- u_{\nu} \Box u^{\nu} 
& = (\partial_{\alpha}u_{\nu})(\partial^{\alpha}u^{\nu}) 
  = \gamma^4 \big( (\partial_x v)^2 - (\partial_t v)^2\big),
\end{align}
where $v = \frac{u}{\sqrt{1+u^2}}$ and $\gamma^2 = \frac{1}{1-v^2}$. Thus, \eqref{Entropy-rel} is equivalent to  
\beq \label{entropy-sign-pde}
TnDs = \varepsilon \gamma^4 \big( (\partial_x v)^2 - (\partial_t v)^2 \big),
\eeq 
and for a traveling wave solution of form $v=v(\xi)$ for $\xi=x-ct$, \eqref{entropy-sign-pde} reduces to 
$$
TnDs = \varepsilon \gamma^4(1-c^2) \left(\frac{dv}{d\xi}\right)^2.
$$ 
We conclude, since $T>0$ and $n>0$, and since $\frac{dv}{d\xi}\neq 0$ in the region where the travelling wave is non-constant, that the entropy production is positive if and only if $|c|<1$, and constant if and only if $|c|=1$, just as claimed.
\QED

\section{Well-posedness and Causality}  \label{Sec_Cauchy}

The goal of this section is to prove Theorem \ref{Thm_existence_causality}, concerning well-posedness and causality of the dissipative relativistic Euler equations \eqref{rel_art_visc}. We begin by writing \eqref{rel_art_visc} in one spatial dimension as a first order system in the unknown variables $\rho$, $u\equiv u^1$, $w\equiv \partial_t u$ and $z\equiv \partial_x u$, substituting again the constraint $u^0 = \sqrt{1+ u^2}$ to implement the normalization $u_\sigma u^\sigma =-1$.

\begin{Lemma}   \label{Lemma_Cauchy_first-order-sys}
Let $u^0 = \sqrt{1+ u^2}$ for $u\equiv u^1$. Then equation \eqref{rel_art_visc} in the unknowns $(\rho,u)$ is equivalent to the first order system
\beq   \label{Cauchy_first-order-sys}
A^0(\U) \partial_t \U + A^1(\U) \partial_x \U = \mathcal{L}(\U),
\eeq 
in the unknowns $\mathcal{U}^T \equiv (\rho, w, z, u)$, for $w\equiv \partial_t u$ and $z\equiv \partial_x u$, with coefficient matrices
\begin{align} \label{Cauchy_def_A0A1}
A^0(\U) \equiv 
\begin{pmatrix}
T^{00}_\rho & \visc v(u) & 0 & 0 \\
T^{01}_\rho & \visc & 0 & 0 \\
0 & 0 & 1 & 0 \\
0 & 0 & 0 & 1 
\end{pmatrix},
\hspace{.5cm}
A^1(\U) \equiv 
\begin{pmatrix}
T^{01}_\rho & 0 & - \visc v(u) & 0 \\
T^{11}_\rho & 0 & -\visc & 0 \\
0 & -1 & 0 & 0 \\
0 & 0 & 0 & 0 
\end{pmatrix},
\end{align}
and source term
\begin{align} \label{Cauchy_def_L}
\mathcal{L}(\U) \equiv 
\begin{pmatrix}
\visc \frac{z^2 - w^2}{(1+ u^2)^{\frac{3}{2}}} - T^{00}_{u} w - T^{01}_u z \\
- T^{01}_u w - T^{11}_u z \\
0 \\
w 
\end{pmatrix}.
\end{align}
where $T^{\mu\nu}_\rho \equiv \partial_\rho T^{\mu\nu}$ and $T^{\mu\nu}_u \equiv \partial_u T^{\mu\nu}$. Moreover, the compatibility condition $z=u_x$ is propagated by \eqref{Cauchy_first-order-sys}.
\end{Lemma}

\Proof 
We start with \eqref{rel_art_visc} and impose the normalization by writing
\[
u^\mu=(u^0,u^1)=\bigl(\sqrt{1+u^2},\,u\bigr),
\]
where $u=u^1$ is the unknown.
Thus, we get
\begin{equation}\label{eq:chainrule-matrix}
\begin{pmatrix}
T^{00}_{\rho} & T^{00}_{u}\\[2pt]
T^{01}_{\rho} & T^{01}_{u}
\end{pmatrix}
\binom{\rho_t}{u_t}
+
\begin{pmatrix}
T^{01}_{\rho} & T^{01}_{u}\\[2pt]
T^{11}_{\rho} & T^{11}_{u}
\end{pmatrix}
\binom{\rho_x}{u_x}
=
\varepsilon
\binom{\Box\sqrt{1+u^2}}{\Box u}.
\end{equation}
Let, $\psi(u)\equiv\sqrt{1+u^2}$. A direct computation gives
\begin{equation}\label{eq:boxpsi}
\Box\psi(u)=\psi''(u)\,(u_x^2-u_t^2)+\psi'(u)\,\Box u.
\end{equation}
Define the quadratic term
\[
Q(u,\partial u)\equiv\psi''(u)\,(u_x^2-u_t^2),
\]
then \eqref{eq:boxpsi} reads
\[
\Box\sqrt{1+u^2}=Q(u,\partial u)+\psi'(u)\,\Box u.
\]
Thus the right-hand side of \eqref{eq:chainrule-matrix} becomes
\[
\varepsilon\binom{\Box\sqrt{1+u^2}}{\Box u}
=
\varepsilon\binom{Q(u,\partial u)}{0}
+
\varepsilon\binom{\psi'(u)}{1}\,\Box u.
\]
Setting now $w\equiv u_t$ and $z\equiv u_x$, along with the compatibility relations
\begin{equation}\label{eq:compatibility}
u_t=w,\qquad z_t-w_x=0,
\end{equation}
we obtain that 
$$
\Box u =  -\partial_{tt}^2 u + \partial_{xx}^2 u = - \partial_t w + \partial_x z ,
$$
and we rewrite \eqref{eq:chainrule-matrix} and \eqref{eq:compatibility} as the first order system,
\begin{equation}\label{eq:firstorder}
A^0(\U)\,\partial_t \U + A^1(\U)\,\partial_x \U = \mathcal{L}(\U),
\end{equation}
where $A^0(\U)$ and $A^1(\U)$ are given by \eqref{Cauchy_def_A0A1}, and
\beq \label{Cauchy_def_L(U)}
\mathcal{L}(U)=
\begin{pmatrix}
\varepsilon Q(u,w,z) - T^{00}_{u}\,w - T^{01}_{u}\,z\\
- T^{01}_{u}\,w - T^{11}_{u}\,z\\
0\\
w
\end{pmatrix},
\eeq 
where $Q(u,w,z)=\psi''(u)\,(z^2-w^2)$. Noting that $\psi'(u) = \frac{u}{\sqrt{1+u^2}}\equiv v(u)$ and $\psi''(u) = (1+u^2)^{-\frac{3}{2}}$ gives us the result. 

Finally, observe that the last two rows in \eqref{Cauchy_first-order-sys}, imply
\begin{equation}\nonumber
z_t-w_x=0,\qquad u_t=w,\qquad \partial_t(z-u_x)=0.
\end{equation}
Hence the compatibility condition $z=u_x$ is propagated by the evolution, (that is, if $z=u_x$ holds at $t=0$, then $z=u_x$ holds for all $t>0$ for which a solution $\mathcal{U}$ of \eqref{Cauchy_first-order-sys} exist).
\QED

The eigenvalues of the generalized eigenvalue problem 
\beq \label{Cauchy_gen_eigenvalue_problem}
\Big( \lambda A^0(\U) - A^1(\U) \Big) \vec{V} = 0,
\eeq  
are 
\beq  \label{Cauchy_eigenvalues}
\lambda_1(\U) = 1 , \hspace{.5cm}
\lambda_2(\U) = 0 , \hspace{.5cm}
\lambda_3(\U) = -1, \hspace{.5cm}
\lambda_4(\U) = v(u) ,
\eeq
with associated linearly independent eigenvectors
\beq \label{Cauchy_eigenvectors}
\vec{V}_1 = \begin{pmatrix} 0 \\ 1 \\ -1 \\ 0  \end{pmatrix}, \hspace{.5cm}
\vec{V}_2 = \begin{pmatrix} 0 \\ 0 \\ 0 \\ 1  \end{pmatrix}, \hspace{.5cm}
\vec{V}_3 = \begin{pmatrix} 0 \\ 1 \\ 1 \\ 0  \end{pmatrix} , \hspace{.5cm}
\vec{V}_4 = \begin{pmatrix} \visc (1-v(u)^2) \\ - p'(\rho) v(u) \\ p'(\rho) \\ 0  \end{pmatrix}.
\eeq 
Note that the above eigenvectors are linearly independent for all values of $u \in \R$, since $v(u) = \frac{u}{\sqrt{1+u^2}} \in (-1,1)$, and all values of $p'(\rho)$. 

We are now prepared to write \eqref{Cauchy_first-order-sys} as a symmetric hyperbolic system. For this, we first multiply \eqref{Cauchy_first-order-sys} by $(A^0)^{-1}$, to obtain
\beq   \label{Cauchy_first-order-sys_B}
\partial_t \U + \mathcal{B}(\U) \partial_x \U = \big(A^0(\U)\big)^{-1}\mathcal{L}(\U),
\eeq 
where $\mathcal{B}(\U) \equiv \big(A^0(\U)\big)^{-1}A^1(\U)$. A direct computation gives
\begin{equation}\label{eq:sec8-A0-inverse}
(A^0(\U))^{-1}=
\begin{pmatrix}
1&-v&0&0\\ -\epsilon^{-1}T^{01}_{\rho}&\epsilon^{-1}T^{00}_{\rho}&0&0\\ 0&0&1&0\\ 0&0&0&1
\end{pmatrix},
\end{equation}
and using the identities $T^{01}_{\rho}-vT^{11}_{\rho}=v$ and $T^{00}_{\rho}T^{11}_{\rho}-(T^{01}_{\rho})^2=p'(\rho)$ as well as $vT^{01}_{\rho}-T^{00}_{\rho}=-1$, 
we obtain
\begin{equation} \label{eq:sec8-B-matrix}
\mathcal B(\U) = \begin{pmatrix} v(u)&0&0&0\\ \epsilon^{-1}p'(\rho)&0&-1&0\\ 0&-1&0&0\\ 0&0&0&0 \end{pmatrix}.
\end{equation}
The eigenvalues of $\mathcal{B}(\U)$ are recorded in \eqref{Cauchy_eigenvalues} and the associated eigenvectors in \eqref{Cauchy_eigenvectors}. These eigenvalues are distinct whenever $u\neq 0$, while the above eigenvectors are always linearly independent, also when $u=0$. 

\begin{Lemma} \label{Cauchy_Lemma_symmetrizer}
A symmetrizer for system \eqref{Cauchy_first-order-sys_B} is 
\beq
\mathcal{S}(\U) = \big(\mathcal{R}^{-1}\big)^T \mathcal{R}^{-1}, 
\eeq 
where $\mathcal{R} \equiv (\vec{V}_1,\vec{V}_2,\vec{V}_3,\vec{V}_4)$ is the matrix diagonalizing $\mathcal{B}$. The resulting symmetric hyperbolic system, equivalent to \eqref{Cauchy_first-order-sys_B}, is
\beq
\mathcal{S}(\U)\partial_t \U + \mathcal{S}(\U)\mathcal{B}(\U) \partial_x \U = \mathcal{F}(\U),
\eeq   
for $\mathcal{F}(\U) \equiv \mathcal{S}(\U)\big(A^0(\U)\big)^{-1}\mathcal{L}(\U)$.
\end{Lemma}

\Proof
We now show that $\mathcal{S}(\U)$ is symmetric and positive definite for all $\U \in \R^4$, and that $\mathcal{S}(\U)\mathcal{B}(\U)$ is a symmetric matrix for all $\U \in \R^4$.

Clearly
\begin{eqnarray}
\mathcal{S}(\U)^T = \Big(\big(\mathcal{R}^{-1}\big)^T \mathcal{R}^{-1} \Big)^T = \big(\mathcal{R}^{-1}\big)^T \mathcal{R}^{-1} = \mathcal{S}(\U),
\end{eqnarray}
while for any $\vec{\xi}\in \R^4$
$$
\langle \mathcal{S}(\U) \vec{\xi} ,\vec{\xi} \rangle = 
\langle \mathcal{R}^{-1}\vec{\xi} , \mathcal{R}^{-1}\vec{\xi} \rangle \geq 0, 
$$
which proves the positive definiteness of $\mathcal{S}(\U)$. 
 
To show that $\mathcal{S}(\U)\mathcal{B}(\U)$ is symmetric, we first note that $\mathcal{B}= \mathcal{R} \Lambda \mathcal{R}^{-1}$ for $\Lambda = {\rm diag}(\lambda_1,\lambda_2,\lambda_3,\lambda_4)$, we then find that
$$
\mathcal{S}\mathcal{B} = \big(\mathcal{R}^{-1}\big)^T  \Lambda \mathcal{R}^{-1}
$$
and thus
\begin{eqnarray}
\Big(\mathcal{S}\mathcal{B}\Big)^T 
= \Big( \big(\mathcal{R}^{-1}\big)^T  \Lambda \mathcal{R}^{-1} \Big)^T  
=   \mathcal{S}\mathcal{B},
\end{eqnarray}
which proves the symmetry.
\QED 

Now, multiplying \eqref{Cauchy_first-order-sys_B} with the symmetrizer $\mathcal{S}$, we obtain
\beq   \label{Cauchy_first-order_sym-hyp}
\mathcal{S}(\U) \partial_t \U + \mathcal{A}(\U)  \partial_x \U = \mathcal{F}(\U),
\eeq 
where $\mathcal{A}(\U) \equiv \mathcal{S}(\U) \mathcal{B}(\U)$ and $\mathcal{F}(\U) \equiv \mathcal{S}(\U)\big(A^0(\U)\big)^{-1}\mathcal{L}(\U)$. System \eqref{Cauchy_first-order_sym-hyp} is a first order system of symmetric hyperbolic PDE's as addressed by Kato's existence theory in \cite{Kato}, from which we now deduce well-posedness.

\subsection{Well-posedness}  
To deduce existence of solutions from Theorem 4.2 in \cite{Kato}, we need to establish the following two technical lemmas, establishing uniform boundedness and Lipschitz continuity with respect to $H^s$-norms of the above coefficients $\mathcal{S}$, $\A$ and $\F$, as functions of $\U$, when restricted to some bounded open subset $\D$ of $H^s (\mathbb{R},\R^4)$. The main observation underlying these lemmas is that $\mathcal{S}(\U)$, $\mathcal{A}(\U)$ and $\mathcal{F}(\U)$ are smooth functions in $\U$. This can be easily verified by noting that $\mathcal{R}(\U)$, $A^0(\U)$ and $\mathcal{L}(\U)$ are smooth functions in $\U=(\rho,w,u,z)$, and likewise that 
\beq \nonumber 
\mathcal{R}^{-1} = \begin{pmatrix}
 \frac{p'(\rho)}{2(1-v)\visc} & \frac{1}{2} & - \frac{1}{2} & 0 \\
0 & 0 & 0 & 1 \\
 \frac{-p'(\rho)}{2(1+v)\visc} & \frac{1}{2} & \frac{1}{2} & 0 \\
 \frac{1}{(1-v^2)\visc} & 0 & 0 & 0 \\
 \end{pmatrix}
\hspace{.3cm} \text{and} \hspace{.3cm}
(A^0)^{-1} = \begin{pmatrix}
1 & - v & 0 & 0 \\
- \frac{T^{01}_\rho}{\visc} &  \frac{T^{00}_\rho}{\visc} & 0 & 0 \\
0 & 0 & 1 & 0 \\
0 & 0 & 0 & 1  
\end{pmatrix}
\eeq 
are smooth in $\U$, since $v(u)= \frac{u}{\sqrt{1+u^2}}$ satisfies $|v(u)|<1$.

\begin{Lemma}  \label{Cauchy_Lemma_uniform-bounds}
Let $s \geq 2$ be an integer. Let $\D \subset H^s(\mathbb{R},\R^4)$ be open and bounded.\\
\textbf{(i)}  $\mathcal{S}(\U)$ and $\mathcal{A}(\U)$ are uniformly bounded on $\D$, in the sense that there exist some constant $\lambda$ independent of $\U$ such that $\|\A(\U)\|_{H^s_{ul}(\R)}\leq \lambda$, for all $\U \in \D$, where
\beq \label{Cauchy_estimate_1} 
\|\A(\U)\|_{H^s_{ul}(\R)}^2  \equiv  \sup_{x\in \R} \sum_{k=0}^s \int_{|x-y|<1} \left\|\frac{d^k}{dy^k}\A(\U(y))\right\|_{\R^{4\times 4}}^2 dy 
\eeq 
and $\|\cdot \|_{\R^{4\times 4}}$ denotes the Hilbert Schmidt norm (or any other matrix norm). \\
\textbf{(ii)}   $\mathcal{F}(\U)$ is uniformly bounded on $\D$, in the sense that there exist some constant $\lambda$ independent of $\U$ such that $\|\mathcal{F}(\U)\|_{H^s(\R)} \leq \lambda$, for all $\U \in \D$. 
\end{Lemma}

\Proof 
\textit{(i):} 
Since $\D$ is assumed to be a bounded subset of $H^s (\mathbb{R},\R^4)$, there exists some $\lambda_0>0$ such that 
\beq \label{Cauchy_Hs-bound_D}
\|\U\|_{H^s} \leq \lambda_0, \qquad \forall\; \U \in \D. 
\eeq  
Thus, Morrey's inequality implies that 
\beq \label{Cauchy_L-infty-bound}
\|\U\|_{L^\infty}\leq C_M \lambda_0, \qquad \forall\; \U \in \D,
\eeq  
for some constant $C_M>0$. Moreover, by smoothness of $\A(\U)$ in $\U$, the mapping $\U \mapsto \|D^\alpha\A(\cdot)\|_{\R^{4\times 4}}$ is a continuous function from $\R^4$ to $\R$, where $D^\alpha \A$ is standard multi-index notation for derivatives of $\A$ with respect to $\U=(\U^1,..., \U^4)$, $|\alpha|\leq s$. Thus, since by \eqref{Cauchy_L-infty-bound} any $\U\in \D$ lies in the compact ball of radius $R\equiv C_M \lambda_0$ in $\R^4$, the function $\U \mapsto \|D^\alpha\A(\cdot)\|_{\R^{4\times 4}}$ attains its supremum, that is,
\beq \label{Cauchy_estimate_2} 
\sup_{\U\in D} \|\A(\U)\|_{\text{Eucl}(s)} \leq K, 
\eeq 
for some constant $K>1$ (taken larger than $1$ for convenience), with respect to the Euclidean norm over multi-index derivatives of $\A$ in $\U$ of order $|\alpha|\leq s\in \N$,
\beq \label{Cauchy_estimate_Eucl-s-norm} 
\|\A(\U)\|_{\text{Eucl}(s)}^2 \equiv \sum_{|\alpha|\leq s} \|D^\alpha\A(\U)\|_{\R^{4\times4}}^2 .
\eeq 
Computing now the derivatives in \eqref{Cauchy_estimate_1} using the F\'aa di Bruno formula (the ``higher order chain rule''), and applying the bound \eqref{Cauchy_estimate_2} to bound all derivatives on $\A$ with respect to $\U$ as well as Morrey's inequality to bound products by the $H^s$-norm, we obtain for some constant $C>0$ that
\beq \label{Cauchy_estimate_2b}
\|\A(\U)\|_{H^s_{ul}(\R)}  \leq C \sum_{k=0}^s\|\U\|^k_{H^s(\R)} \leq \lambda,
\eeq 
where the last inequality follows from the $H^s$-bound \eqref{Cauchy_Hs-bound_D} on $\D$. 

To illustrate the last step leading to \eqref{Cauchy_estimate_2b}, consider the case $s=2$,
\begin{eqnarray} \label{Cauchy_estimate_3}
\|\A(\U)\|_{H^2_{ul}(\R)}^2  &=&   \sup_{x\in \R} \sum_{k=0}^2 \int_{|x-y|<1} \left\|\frac{d^k}{dy^k}\A(\U(y))\right\|_{\R^{4\times 4}}^2 dy 
\end{eqnarray} 
and estimate each of the three terms in the integrand separately as follows:
\begin{eqnarray} \nonumber 
\left\|\A(\U(y))\right\|_{\R^{4\times 4}} 
&\leq& \sup_{\U\in D} \|\A(\U)\|_{\R^{4\times4}} 
\overset{\eqref{Cauchy_estimate_2}}{\leq}  K      ,
\end{eqnarray}
and, denoting matrix coefficients by $\A^i_j(\U)$, we get for the first derivative term
\begin{eqnarray} \nonumber 
\left\|\frac{d}{dy}\A(\U(y))\right\|_{\R^{4\times 4}}^2 
&= & \sum_{i,j=1}^4 \left|\frac{d}{dy}\A^i_j(\U(y))\right|^2
= \sum_{i,j,l=1}^4 \left|\frac{\partial\A^i_j}{\partial \U^l}\frac{d\U^l}{dy}\right|^2  \cr 
&\leq &  \left\|\A(\U)\right\|_{\text{Eucl}(1)}^2\cdot \left\|\frac{d\U}{dy}\right\|_{\R^{4}}^2
\overset{\eqref{Cauchy_estimate_2}}{\leq}  K^2 \left\|\frac{d\U}{dy}\right\|_{\R^{4}}^2 ,  
\end{eqnarray} 
while we bound the second derivative term as   
\begin{eqnarray} \nonumber 
\left\|\frac{d^2}{dy^2}\A(\U(y))\right\|_{\R^{4\times 4}}^2 
&=& \sum_{i,j=1}^4 \left|\frac{d^2}{dy^2}\A^i_j(\U(y))\right|^2  \cr 
&\leq & \sum_{i,j,l,k=1}^4  \left|\frac{\partial^2 \A^i_j}{\partial \U^l \partial \U^k } \frac{d\U^l}{dy} \frac{d\U^k}{dy}\right|^2 +  \sum_{i,j,l=1}^4 \left|\frac{\partial \A^i_j}{\partial \U^l} \frac{d^2\U^l}{dy^2} \right|^2  \cr 
&\leq & \left\|\A(\U)\right\|_{\text{Eucl}(2)}^2 \cdot \left( \sum_{k=1,2} \left\|\frac{d^k\U}{dy^k}\right\|_{\R^{4}}^{2(3-k)} \right) ,
\end{eqnarray} 
where $\left\|\A(\U)\right\|_{\text{Eucl}(2)} \leq K$ by \eqref{Cauchy_estimate_2}. Now, substituting the above three estimates back to replace the bound the integrand in \eqref{Cauchy_estimate_3}, and using the bound \eqref{Cauchy_estimate_2} on $\|\A(\U)\|_{\text{Eucl}(s)}$, as well as Morrey's inequality to bound powers of $\left\|\frac{d^k\U}{dy^k}\right\|_{\R^{4}}$ by powers of the $H^s$-norm on $\U$, (noting that $H^s(\R)$ is a multiplicative algebra), we finally obtain
\begin{eqnarray} \nonumber 
\|\A(\U)\|_{H^2_{ul}(\R)}^2  
&\leq &  3K^2 \sum_{k=0}^2\|\U\|_{H^2(\R)}^{2k} 
\leq  3K^2  \sum_{k=0}^2 \lambda_0^{2k} \equiv \lambda^2 ,
\end{eqnarray}
where we used the $H^s$-bound \eqref{Cauchy_Hs-bound_D} on $\D$ for the last inequality. The general estimate \eqref{Cauchy_estimate_2b} is obtained in a similar fashion by successively increasing $s$. The uniform bound on $\|\mathcal{S}(\U)\|_{H^s_{ul}(\R)}$ follows by the same argument, replacing $\A(\U)$ by $\mathcal{S}(\U)$. This completes the proof of part (i). \vspace{.2cm}

\noindent \textit{(ii):}
This follows by an argument similar to the one above, because  $\mathcal{F}(\U)$ is smooth in $\U$, as verified above the lemma. In particular, in analogy to \eqref{Cauchy_estimate_2}, we obtain the bound 
\beq \nonumber 
\sup_{\U\in D} \|\F(\U)\|_{\text{Eucl}(s)} \leq K,
\hspace{.8cm} \text{for } \quad 
\|\F(\U)\|_{\text{Eucl}(s)}^2 \equiv \sum_{|\alpha|\leq s} \|D^\alpha\F(\U)\|_{\R^{4}}^2 
\eeq 
for some constant $K>1$. We then get, applying again the Fa\'a di Bruno formula to compute derivatives, the estimate     
\begin{eqnarray} \nonumber 
\|\mathcal{F}(\U)\|_{H^s(\R)}^2   
&\leq & C_M K^2 \sum_{k=0}^s\|\U\|_{H^s(\R)}^{2k}  
\; \leq\;  C_M K^2 \sum_{k=0}^s \lambda_0^{2k} \;\equiv\; \lambda^2,
\end{eqnarray}
where we used Morrey's inequality to bound powers of $L^2$-integrants by powers of $H^s$-norms, and we used the $H^s$-bound \eqref{Cauchy_Hs-bound_D} on $\D$ for the last inequality. This completes the proof.
\QED

\begin{Lemma} \label{Cauchy_Lemma_Lipschitz-bounds}
Let $s \geq 2$ be an integer. Let $\D \subset H^s (\mathbb{R},\R^4)$ be open and bounded.\\
\textbf{(i)} The coefficients $\mathcal{A}(\U)$ and $\mathcal{S}(\U)$ are uniformly Lipschitz continuous on $\D$, in the sense that there exist some constant $\lambda >0$ independent of $\U \in \D$ such that 
\beq \label{Cauchy_Lip-bdd_eqn1}
\|\A(\U_1) - \A(\U_2)\|_{H^{s}_{\mathrm{ul}}(\R)}  \leq \lambda\; \|\U_1-\U_2\|_{H^s(\R)},
\eeq 
with analogous estimates for $\mathcal{S}(\U)$, and $H^{s}_{\mathrm{ul}}$-norm defined in \eqref{Cauchy_estimate_1}.\\
\textbf{(ii)} The source term $\mathcal{F}(\U)$ is uniformly Lipschitz continuous on $\D$, in the sense that there exist some constant $\lambda >0$ independent of $\U \in \D$ such that 
\beq \label{Cauchy_Lip-bdd_eqn2}
\|\F(\U_1) - \F(\U_2)\|_{H^{s}(\R)}  \leq \lambda\; \|\U_1-\U_2\|_{H^s(\R)}.
\eeq 
\end{Lemma}

\Proof  
We only prove the part (i) for $\A(\U)$. The argument for establishing Lipschitz bounds \eqref{Cauchy_Lip-bdd_eqn1} for $\mathcal{S}(\U)$, as well as \eqref{Cauchy_Lip-bdd_eqn2} on $\F(\U)$ in part (ii), is similar and omitted here. To establish now the Lipschitz estimate \eqref{Cauchy_Lip-bdd_eqn1} for $\A(\U)$, we first write the norm \eqref{Cauchy_estimate_1} as
\beq \label{Cauchy_estimate_5}
\|\A(\U)\|_{H^s_{ul}(\R)}^2  \equiv  \sum_{k=0}^s \left\|\frac{d^k}{dy^k}\A(\U)\right\|_{H^{0}_{\mathrm{ul}}}^2.
\eeq
We now establish \eqref{Cauchy_Lip-bdd_eqn1} step-by-step, successively increasing $s\geq 0$ in \eqref{Cauchy_estimate_5}. 

By smoothness of $\A$ in $\U$, we have for any $\U_1,\U_2 \in \D$ that
\begin{eqnarray}  \label{Cauchy_Lipschitz-bound-00} 
\|\A(\U_1) - \A(\U_2)\|_{\R^{4\times 4}} 
&\leq & \sup_{\U \in D}\|\A(\U)\|_{\text{Eucl}(1)} \; \|\U_1 - \U_2\|_{\R^4}  \cr 
&\leq &  K \|\U_1 - \U_2\|_{\R^4}, 
\end{eqnarray}
where the last inequality follows by \eqref{Cauchy_estimate_2}. Integration, as in \eqref{Cauchy_estimate_2b}, then implies  
\beq \label{Cauchy_Lipschitz-bound-2}
\|\A(\U_1) - \A(\U_2)\|_{H^{0}_{\mathrm{ul}}} 
\leq K \|\U_1 - \U_2\|_{H^{0}_{\mathrm{ul}}} 
\leq K \|\U_1 - \U_2\|_{H^{0}},  
\eeq 
which is the sought-after Lipschitz bound for $s=0$. 

We now address the case $s=1$. Expressing the chain rule as
$$
\frac{d\U}{dy} = \sum_{l=1}^4\frac{\partial\A}{\partial \U^l}\frac{d\U^l}{dy}  \equiv D\A(\U_1)\cdot\frac{d\U_1}{dy},
$$ 
we obtain  
\begin{align} \label{Cauchy_estimate_5.1}
&\left\|\frac{d}{dy} \left(\A(\U_1) - \A(\U_2)\right)\right\|_{H^{0}_{\mathrm{ul}}}
\leq \left\|D\A(\U_1)\frac{d\U_1}{dy}-D\A(\U_2)\frac{d\U_2}{dy}\right\|_{H^{0}_{\mathrm{ul}}}\cr 
&\leq \left\|\big(D\A(\U_1) - D\A(\U_2)\big) \frac{d\U_1}{dy}\right\|_{H^{0}_{\mathrm{ul}}} + \left\| \Big(\frac{d\U_1}{dy} - \frac{d\U_2}{dy}\Big)  D\A(\U_2)\right\|_{H^{0}_{\mathrm{ul}}} .
\end{align}
To bound the first term, note that by smoothness of $\A$ we obtain similar to \eqref{Cauchy_Lipschitz-bound-00} the Lipschitz bound  
\begin{align} \label{Cauchy_Lipschitz-bound-000} 
\big\|\A(\U_1) - \A(\U_2)\big\|_{\text{Eucl}(1)} 
& \leq \sup_{\U \in \D} \|\A(\U)\|_{\text{Eucl}(2)} \cdot \; \big\| \U_1-\U_2 \big\|_{\R^4} \cr 
&\leq K \; \big\| \U_1-\U_2 \big\|_{\R^4} ,
\end{align} 
where the last inequality follows from the Euclidean bound \eqref{Cauchy_estimate_2}. Integration  then allows us to bound the first term in \eqref{Cauchy_estimate_5.1} as
\begin{align} \nonumber 
\left\|\big(D\A(\U_1) - D\A(\U_2)\big) \frac{d\U_1}{dy}\right\|_{H^{0}_{\mathrm{ul}}} 
& \leq   \left\|\big\|\A(\U_1) - \A(\U_2)\big\|_{\text{Eucl}(1)} \left\|\frac{d\U_1}{dy}\right\|_{\R^4}\right\|_{H^{0}_{\mathrm{ul}}}   \cr 
& \leq K C_M  \;  \|\U_1 - \U_2\|_{H^1(\R)} \; \big\|\U_1\big\|_{H^{1}(\R)}  \cr 
& \leq K C_M \lambda_0 \|\U_1 - \U_2\|_{H^1(\R)},
\end{align} 
where we used Morrey's inequality (with constant $C_M>0$) to bound the $L^\infty$-norm on $\U_1-\U_2$ by the $H^1$-norm. We bound the second term in \eqref{Cauchy_estimate_5.1} by
\begin{align} \nonumber 
\Big\| \Big(\frac{d\U_1}{dy} - \frac{d\U_2}{dy}\Big)  D\A(\U_2)\Big\|_{H^{0}_{\mathrm{ul}}}
&\leq \sup_{\U\in \D} \Big\|\A(\U)\Big\|_{\text{Eucl}(1)} \; \Big\| \Big(\frac{d\U_1}{dy} - \frac{d\U_2}{dy}\Big)\Big\|_{H^{0}_{\mathrm{ul}}} \cr 
&\leq K \|\U_1 - \U_2\|_{H^1}.
\end{align}
Combining the above two estimates, we bound \eqref{Cauchy_estimate_5.1} by  
\begin{align} \label{Cauchy_estimate_5b}
\Big\|\frac{d}{dy} & \left(\A(\U_1) - \A(\U_2)\right)\Big\|_{H^{0}_{\mathrm{ul}}}
 \leq  K(C_M \lambda_0 + 1)  \|\U_1 - \U_2\|_{H^1(\R)}.
\end{align}
This in combination with \eqref{Cauchy_Lipschitz-bound-2} then implies the sought-after Lipschitz bound \eqref{Cauchy_Lip-bdd_eqn1} for $s=1$, with $C\equiv K(C_M \lambda_0 + 2)$. 

Continuing in this fashion, employing the Fa\'a di Bruno formula in place of the chain rule, and using that $\U \to \A(\U)$ is smooth, gives rise to uniform Lipschitz estimates in $\U$ (similar to \eqref{Cauchy_Lipschitz-bound-00} and \eqref{Cauchy_Lipschitz-bound-000}) for all derivatives, as long that the derivatives on $\U$ are bounded by the $H^s$-norm. From this one can then deduce the desired Lipschitz bound
\beq \label{Cauchy_Lipschitz-bound-3}
\|\A(\U_1) - \A(\U_2)\|_{H^{s}_{\mathrm{ul}}} \leq C \; \|\U_1-\U_2 \|_{H^s}
\eeq 
for some finite constant $C>0$ depending only on $\sup_{\U\in \D} \|\A(\U)\|_{\text{Eucl}(s+1)}<\infty$, which is finite since $\A(\U)$ is smooth in $\U$ and since $\D$ is bounded in the sense that $\|\U\|_{L^\infty(\R)}<R$ for all $\U\in \D$, see \eqref{Cauchy_L-infty-bound}. This completes the proof of Lemma \ref{Cauchy_Lemma_Lipschitz-bounds}. 
\QED 

We can now establish our basic well-posedness result:

\begin{Thm} \label{Cauchy_Thm_existence}
Let $\D \subset H^s (\mathbb{R},\R^4)$ be open and bounded, for $s \geq 2$ an integer.\\
\noindent \textbf{(i)} Assume initial data $\U_0=(\rho_0, w_0, z_0, u_0) \in \D$. Then there is a unique solution $\U$ of \eqref{Cauchy_first-order_sym-hyp}, defined on $[0, T]$ for some $T>0$, satisfying $\U(0)=\U_0$, and such that
\beq \label{Cauchy_Thm_existence_reg} 
\U \in C([0,T]; \D) \cap C^1([0,T]; H^{s-1}(\R,\R^4)).
\eeq 
Moreover, $T$ can be chosen common to all initial conditions $\U_0$ varying in a small neighborhood of a given point in $\D$.

\noindent \textbf{(ii)} Solutions of \eqref{Cauchy_first-order_sym-hyp} depend continuously on the initial data in the following sense: Let $\U_0, \U_0^n \in \D$, $n\in \mathbb{N}$, and let $\U, \U^n$ be the unique solutions of \eqref{Cauchy_first-order_sym-hyp} with data $\U(0)=\U_0$ and $\U^n(0)=\U_0^n$, (asserted to exist in part (i)). Assume $\U$ exists on a closed interval $[0, T_0]$. If $\lim_{n\to\infty}\| \U_0^n - \U_0 \|_{H^s}=0$, then the solutions $\U^n$ also exist on $[0, T_0]$ for sufficiently large $n$, and $\lim_{n\to\infty}\| \U^n(t) - \U(t) \|_{H^s} = 0$ uniformly in $t \in [0,T_0]$. 
\end{Thm}

\Proof
\noindent \textbf{(i):}
The result follows by Theorem II in \cite{Kato} upon noting that all assumptions in \cite[Thm. II]{Kato} are met, as follows: Lemma \ref{Cauchy_Lemma_uniform-bounds} implies assumption (4.2) in \cite[Thm. II]{Kato}; Lemma \ref{Cauchy_Lemma_Lipschitz-bounds} implies (4.3) - (4.4) in \cite[Thm. II]{Kato} as well as, (since the coefficients in \eqref{Cauchy_first-order_sym-hyp} do not depend explicitly on time),  (4.5) - (4.6) in \cite[Thm. II]{Kato}; while symmetry of $\A$ and positivity of $\S$, as required by (4.7) - (4.8) in \cite[Thm. II]{Kato}, follows from Lemma \ref{Cauchy_Lemma_symmetrizer}. Therefore Theorem II in \cite{Kato} applies and yields solutions $\U$ of \eqref{Cauchy_first-order_sym-hyp} with regularity \eqref{Cauchy_Thm_existence_reg}  for any data $\U_0$ in $\D$, exactly as asserted in part (i) of Theorem \ref{Cauchy_Thm_existence}. 

\textbf{(ii):} Lemma \ref{Cauchy_Lemma_Lipschitz-bounds} implies condition (4.11) of Theorem III (b) in \cite{Kato}, (that is, uniform Lipschitz continuity with respect to $H^s$-norms). Theorem III (b) in \cite{Kato} then implies part (ii) of Theorem \ref{Cauchy_Thm_existence} exactly as stated, completing the proof.
\QED 

The following version of continuous dependence, (an alternative to the one of part (ii) of Theorem \ref{Cauchy_Thm_existence}), holds as well.

\begin{Corollary}
Let $\U_0, \U_0^n \in \D, n=1, 2, \dots$, and let $\U, \U^n$ be the unique solutions of \eqref{Cauchy_first-order_sym-hyp} with $\U(0)=\U_0, \U^n(0)=\U_0^n$. Then there exist positive numbers $C$ and $T' \le T$, depending on $\U_0$, such that if $\|\U_0^n -\U_0\|_s \le C$ and $\| \U_0^n - \U_0 \|_0 \to 0$ as $n \to \infty$, then the solutions $\U^n$ exist on a common interval $[0, T']$ and $\| \U^n(t) - \U(t) \|_{s-1} \to 0$ uniformly in $t$. In particular, $\U^n(t, x) \to \U(t, x)$ uniformly on $[0, T'] \times \R^m$. 
\end{Corollary}

\Proof
This follows directly from (a) in \cite[Thm.III]{Kato}, which requires the same incoming assumptions than \cite[Thm.II]{Kato}, already established in the proof of Theorem \ref{Cauchy_Thm_existence} above.
\QED

Since Theorem \ref{Cauchy_Thm_existence} requires data in $H^s(\R)$, which hence decays at infinity, Theorem \ref{Cauchy_Thm_existence} excludes constant density states and only yields solutions of \eqref{rel_art_visc} for data $\rho_0$ decaying to zero. For the proof of Theorem \ref{Thm_existence_causality} we need the following variation of Theorem \ref{Cauchy_Thm_existence} in order to include constant reference densities $\bar{\rho}\geq 0$, and initial densities that decay to $\bar{\rho} > 0$ at infinity. For this, we introduce the new variable 
$$
W \equiv \U-\overline \U ,
$$
where $\overline \U$ is a constant reference state, e.g., $\overline \U=(\rho_0,0,0,0)^{\mathsf T}$ and $W =(\rho-\rho_0,w,z,u)^{\mathsf T}$. Substitution of $\U= \overline U+W$ into the symmetric-hyperbolic system \eqref{Cauchy_first-order_sym-hyp} gives the equivalent equation
\beq \label{Cauchy_first-order_sym-hyp-W}
\mathcal{S}(\overline \U+W)W_t + \mathcal{A}(\overline \U+W)W_x = \mathcal{F}(\overline \U+W).
\eeq 
We assume that the set $\mathcal D\subset H^s(\mathbb R;\mathbb R^4)$ is such that 
\[
K_{\mathcal D} \equiv \overline{\{\bar \U +W(x) \mid  W\in\mathcal D,\ x\in\mathbb R\} \subset \R^4},
\]
which is by Sobolev embedding a compact set in $\R^4$, is contained in the domain on which $\mathcal{S}$, $\mathcal{A}$ and  $\mathcal{F}$ are defined and smooth and $S$ is positive definite.

\begin{Corollary}    \label{Cauchy_Cor_exist_const-states}
Let $\overline \U \in \mathbb R^4$ be a constant state satisfying $\mathcal{F}(\overline \U)=0$, let $s\geq2$ and let $\mathcal{D}$ be a bounded neighborhood in $H^s(\mathbb R;\mathbb R^4)$ containing the origin, such that $\mathcal{S}$, $\mathcal{A}$ and $\mathcal{F}$ are defined and smooth (and $\mathcal{S}$ positive definite) on $K_\mathcal{D}$.  \\
\textbf{(i)} Then, for every data $W_0\in \mathcal{D}$, \eqref{Cauchy_first-order_sym-hyp-W} has a unique solution  
$$
W\in C([0,T];\mathcal{D}) \cap C^1([0,T];H^{s-1}(\R,\mathbb R^4)),
$$
and $T>0$ may be chosen uniformly for data sufficiently close to $W_0$ in $H^s$. \\
\textbf{(ii)} Solutions of \eqref{Cauchy_first-order_sym-hyp-W} depend continuously on initial data: Let $W_0, W_0^n \in \D$, $n\in \mathbb{N}$, and let $W, W^n$ be the unique solutions of \eqref{Cauchy_first-order_sym-hyp-W} with data $W(0)=W_0$ and $W^n(0)=W_0^n$. Assume $W$ exists on a closed interval $[0, T_0]$. If $\lim_{n\to\infty}\| W_0^n - W_0 \|_{H^s}=0$, then the solutions $W^n$ also exist on $[0, T_0]$ for sufficiently large $n$, and $\lim_{n\to\infty}\| W^n(t) - W(t) \|_{H^s} = 0$ uniformly in $t \in [0,T_0]$. 
\end{Corollary}

\Proof 
Since $F(\overline \U)=0$, the source belongs to $H^s$ whenever $W\in H^s$.  
The coefficient bounds and Lipschitz estimates of Lemmas \ref{Cauchy_Lemma_uniform-bounds} and \ref{Cauchy_Lemma_Lipschitz-bounds}, used in Theorem \ref{Cauchy_Thm_existence}, hold with $\U$ replaced by $\overline \U+W$.  Thus, as in Theorem \ref{Cauchy_Thm_existence}, the conclusion (i) and (ii) follow from Theorem II in \cite{Kato}.
\QED

\subsection{Causality}
We next establish causality of the first order PDE \eqref{Cauchy_first-order_sym-hyp}.

\begin{Thm} \label{Cauchy_Thm_causality}
The solutions of \eqref{Cauchy_first-order_sym-hyp}, asserted to exist in Theorem \ref{Cauchy_Thm_existence} (i), as well as the solution of \eqref{Cauchy_first-order_sym-hyp-W} established in Corollary \ref{Cauchy_Cor_exist_const-states}, are all causal in the sense that their backward Cauchy development is contained inside the interior of the backward light-cone.
\end{Thm}

The main technical step in the proof of Theorem \ref{Cauchy_Thm_causality} is to establish the following energy estimate in $\U$-variables.

\begin{Lemma}
Let $\U$ be a solution of \eqref{Cauchy_first-order_sym-hyp} as in Theorem \ref{Cauchy_Thm_existence_reg} (i). Then, the energy functional 
\begin{equation}\label{Cauchy_energy_def}
E(t)\equiv  \ \frac12\,\int_{I(t)} \U^{T} \mathcal{S}(\U)\U \,dx, 
\end{equation}
defined with respect to some interval $I(t) = [a(t),b(t)]$ such that $(t,a(t))$ and $(t,b(t))$ lie on the backward light cone of some point $(t_0,x_0)$, [i.e., $a(t)\equiv  x_0-(t_0-t)$ and $ b(t)\equiv x_0+(t_0-t)$ for $0\leq t \leq t_0 \leq T$, $x_0 \in \R$], satisfies 
\begin{equation}\label{Cauchy_energy-estimate}
\frac{d}{dt}E(t)\leq C \,E(t) ,
\end{equation}
where $C= C(\U) > 0$ is a constant depending only on $\lambda_0$ (the $H^s$-bound on $\D$).
\end{Lemma}

\Proof
Recall that the symmetric hyperbolic PDE \eqref{Cauchy_first-order_sym-hyp} has coefficients $\mathcal{S}(\U)$ and $\mathcal{A}(\U) = \mathcal{S}(\U)\B(\U)$, where $\mathcal{S}(\U)$ is a positive definite $4\times 4$-matrix and $\B(\U)$ has for all $\U$ eigenvalues lying in $[-1,1]$, given by \eqref{Cauchy_eigenvalues}.  By \eqref{Cauchy_Thm_existence_reg}, solutions $\U(t,\cdot) \in \D$ for all $t\in [0,T]$, and  are thus bounded in $L^\infty$ by \eqref{Cauchy_L-infty-bound}. Thus, since $\mathcal{S}(\U)$ is positive definite and continuous in $\U$ and since $\D$ is a bounded set, it follows that there exists constants $0 < \lambda \leq \Lambda <\infty$ such that
\begin{equation}\label{S_bounds}
\lambda I \le \mathcal{S}(U) \le \Lambda I,  \qquad \forall\; \U \in \D.
\end{equation}
Moreover, by smoothness of $\mathcal{S}(\U)$ and $\A(\U)$ in $\U$, we obtain for the supremums norm $\|\U\|_\infty \equiv {\rm sup}_{t\in [0,T]} \|\U(t,\cdot)\|_{L^\infty(\R)}$ that
\begin{equation}\label{eq:coeff-bounds}
\|\partial_t \mathcal{S}(\U)\|_{\infty} \le C
\qquad \text{ and } \qquad 
\|\partial_x\big(\A(\U)\big)\|_{\infty} \le C,
\end{equation}
for some constant $C>0$ depending on $\|\U\|_{H^2(\R)}$, where we used that the regularity of $\U$ asserted in \eqref{Cauchy_Thm_existence_reg} implies by Morrey's inequality $\|\U_x(t,\cdot)\|_{L^\infty(\R)} \leq C \|\U(t,\cdot)\|_{H^2(\R)} \leq \lambda_0$, (where $\lambda_0>$ is the upper bound in \eqref{Cauchy_Hs-bound_D} on $\D$), and by substituting the PDE for $\U_t$ 
$$
\|\U_t(t,\cdot)\|_{L^\infty(\R)}\leq \|\U_t(t,\cdot)\|_{H^1(\R)} \leq {\rm max}_{t\in [0,T]} C \|\U(t,\cdot)\|_{H^2(\R)} < C \lambda_0.
$$ 

Now, defining the point-wise energy density and flux
\begin{equation}\label{eq:ef-def}
e(t,x) \equiv \frac12\,\U^{T} \mathcal{S}(\U)\,\U,\qquad
f(t,x) \equiv \frac12\,\U^{T} \A(\U)\,\U,
\end{equation}
\eqref{S_bounds} implies the coercivity bounds
\begin{equation}\label{eq:e-bounds}
\frac{\lambda}{2}\,\|\U\|_{\R^4}^2 \le e(t,x) \le \frac{\Lambda}{2}\,\|\U\|_{\R^4}^2.
\end{equation}
Moreover, a straight-forward computation gives
\begin{align}\nonumber 
\partial_x f
& = \; \U^T \A \U_x + \frac12\,\U^T \A_x \U, \cr 
\partial_t e \;
&= \; \U^T \mathcal{S} \U_t + \frac12\,\U^T \mathcal{S}_t \U 
 \; =\;  -\,\U^T \A \,\U_x + \U^T \mathcal{F}   
+ \frac12\,\U^T \mathcal{S}_t \U,
\end{align}
where we substituted the PDE \eqref{Cauchy_first-order-sys}, as $\mathcal{S}\U_t = - \A \U_x + \mathcal{F}$, in the last step. Combining the above two equations gives the balance law
\beq \label{eq:local-balance}
\partial_t e + \partial_x f = R,  
\eeq 
for
\beq \nonumber 
R(\U)\equiv \U^T \mathcal{F}(\U) + \frac12\,\U^T \mathcal{S}(\U)_t \U + \frac12\,\U^T \A(\U)_x \U.
\eeq

Now, consider the point $(x_0,t_0)$ and recall that
\begin{equation}\label{eq:interval}
I(t)\equiv[a(t),b(t)],\qquad
a(t)\equiv x_0-(t_0-t),\qquad b(t)\equiv x_0+(t_0-t),
\end{equation}
for $t\le t_0$; so the closure of the interior backward light cone is $\bigcup_{0 \leq t \leq t_0}I(t)$. For the energy defined in \eqref{Cauchy_energy_def}, the Leibniz' rule for moving boundaries together with \eqref{eq:interval} imply 
\begin{align}
\frac{d}{dt}E(t)
&=\int_{a(t)}^{b(t)} \partial_t e\,dx - e(t,b(t)) - e(t,a(t)),
\end{align}
and substituting \eqref{eq:local-balance} yields
\begin{align} \label{eq:tech1}
\frac{d}{dt}E(t)
&=\int_{a(t)}^{b(t)} \big(-\partial_x f + R(\U)\big)\,dx - e(b) - e(a) \notag\\
&=\big[f(a)-e(a)\big] - \big[f(b)+e(b)\big] + \int_{I(t)} R(\U)\,dx .
\end{align}
To estimate the boundary terms in \eqref{eq:tech1}, it is convenient to introduce
\begin{equation}\label{eq:wB-def}
w \equiv \mathcal{S}^{1/2} \U,\qquad 
\text{ and } \qquad 
\tilde{\B} \equiv \mathcal{S}^{1/2} \B \mathcal{S}^{-1/2},
\end{equation}
where we used that $\mathcal{S}$ is positive definite and symmetric. Then
\begin{equation}\label{eq:e-in-w}
e= \frac12\,w^T w
\hspace{1cm} \text{and} \hspace{1cm}
f = \frac12\,w^T \tilde{\B} w,
\end{equation}
and since the characteristic polynomial of $\tilde{\B}$ and $\B$ differs only by a positive factor, $\tilde{\B}$ has the same eigenvalues as $\B$. Thus, since these eigenvalues are all contained in $[-1,1]$, and since $\tilde{\B}$ is symmetric, the operator norm on $\tilde{\B}$ is bounded by $1$, which implies that
\begin{equation}\label{eq:flux-bound}
|f| = \frac12|w^T \tilde{\B} w| \leq \frac12 \|\tilde{\B}\|_{\mathrm{op}}\,|w^T w| \le \frac12|w^T w| = e.
\end{equation}
In particular,
\begin{equation}\label{eq:boundary-signs}
f(a)-e(a)\le 0,\qquad f(b)+e(b)\ge 0,
\end{equation}
so the boundary contribution in \eqref{eq:tech1} is non-positive,
\[
\big[f(a)-e(a)\big] - \big[f(b)+e(b)\big] \leq 0,
\]
which implies
\begin{equation}\label{eq:E-ineq}
\frac{d}{dt}E(t)\le \int_{I(t)} R(U)\,dx.
\end{equation}

It remains to bound \eqref{eq:E-ineq} by the total energy. From the uniform bounds \eqref{eq:coeff-bounds}, together with the coercivity bounds \eqref{eq:e-bounds}, we obtain
\[
\frac12\int_{I(t)} \U^T \mathcal{S}_t \U\,dx \le \frac12\|\partial_t \mathcal{S}(\U)\|_\infty \int_{I(t)} \|\U\|_{\R^4}^2 dx
\le C\,E(t),
\]
and similarly
\[
\frac12\int_{I(t)} \U^T \A(\U)_x \U \,dx 
\leq \frac12\|\partial_x\A(\U)\|_\infty \int_{I(t)} \|\U\|_{\R^4}^2 dx
\le C\,E(t).
\]
Moreover, since $\mathcal{F}(\U)$ is Lipschitz continuous in $\U$ with $\mathcal{F}(0)=0$, it follows that $\|\mathcal{F}(\U)\|_{\R^4}\le C\|\U\|_{\R^4}$ for all $\U \in \D$, (since $\D$ is bounded by \eqref{Cauchy_L-infty-bound}), and thus 
\begin{align} \nonumber 
\int_{I(t)} \U^T \mathcal{F}(\U)  \,dx
&\leq \int_{I(t)} \|\U\|_{\R^4}\,\|\mathcal{F}(\U)\|_{\R^4} \,dx
\le C\, \int_{I(t)}\|\U\|_{\R^4}^2\,dx
\le C\,E(t).
\end{align}
Combining the above three inequalities implies
\begin{equation}\label{eq:R-bound}
\int_{I(t)} R(\U)\,dx \le C\,E(t),
\end{equation}
for some constant $C>0$. Combining finally \eqref{eq:E-ineq} and \eqref{eq:R-bound} yields the sought-after energy estimate \eqref{Cauchy_energy-estimate}.
\QED 

\Proof[Proof of Theorem \ref{Cauchy_Thm_causality}:]
By Gr\"onwall's inequality the energy estimate \eqref{Cauchy_energy-estimate} implies
\begin{equation}\label{eq:E-gronwall}
E(t)\le E(0)\,e^{Ct}, \hspace{1cm} \forall \; t \in [0,t_0],
\end{equation}
where $C=C(\U)>0$ is a time-independent constant. Now, assuming the initial data satisfy $\U(0,x)=0$ for all $x\in I(0)$, we have $E(0)=0$ by our definition of energy in \eqref{Cauchy_energy_def}. As a consequence $E(t)=0$ for all $t\in[0,t_0]$, which by \eqref{Cauchy_energy_def} in turn implies that $\U(t,x)=0$ for all $x\in I(t)$ and all $t\in [0,t_0]$. That is, $\U$ vanishes inside the backward light-cone of $(t_0,x_0)$. This proves causality of solutions; data outside of $I(0)$ cannot influence the solution inside the backward light-cone of $(t_0,x_0)$. 

The causality of solutions $W$ of \eqref{Cauchy_first-order_sym-hyp-W} follows by the  same reasoning, using the adapted energy  
\begin{equation}\label{Cauchy_energy_def-W}
E_W(t)\equiv  \ \frac12\,\int_{I(t)} W^{T} \mathcal{S}(W+\bar{\U})W \,dx, 
\end{equation}
and establishing the energy estimate $\frac{d}{dt}E_W(t)\leq C \,E_W(t)$.
\QED

\subsection{Proof of Theorem \ref{Thm_existence_causality}}  
Set $\overline \U=(\overline\rho,0,0,0)^{\mathsf T}$, and $W=\U-\overline \U=(\rho-\overline\rho,w,z,u)^{\mathsf T}$. Then $\overline U$ is an equilibrium, $\mathcal{F}(\overline \U)=0$. The initial data $\rho_0-\bar{\rho} \in \mathcal{D}_s$, $u_0 \in \mathcal{D}_{s+1}$ and $w_0 \in \mathcal{D}_{s}$ in Theorem \ref{Thm_existence_causality}, gives rise to compatible data $W_0$, (i.e., data satisfying $\partial_x u_0 = z_0$), with $W_0 \in \mathcal{D}$. Now,  Corollary \ref{Cauchy_Cor_exist_const-states} applies and gives a unique local compatible solution
$$
W\in C([0,T];\mathcal{D}) \cap C^1([0,T];H^{s-1}),
$$
satisfying $z=\partial_x u$ everywhere. The first and fourth component of $W$, form hence the solution $(\rho,u)$ asserted in Theorem \ref{Thm_existence_causality} (i). The continuous-dependence assertion in Theorem \ref{Thm_existence_causality} follows now directly from Corollary \ref{Cauchy_Cor_exist_const-states} (ii).

Part (iii) of Theorem \ref{Thm_existence_causality}, concerning causality of the dissipative relativistic Euler equations \eqref{rel_art_visc}, follows directly from Theorem \ref{Cauchy_Thm_causality}. Namely, considering the solution $(\rho,u)$ with initial data $(\rho_0-\bar{\rho},u_0,w_0)$ vanishing on the interval $I(0)$, (the intersection of the backward light-cone of a point $(t_0,x_0)$ with the Cauchy surface at $t=0$), then the associated variable $W = (\rho-\bar{\rho}, w, z, u)$ vanishes on $I(0)$ as well. By \ref{Cauchy_Thm_causality}, it follows that $W$ and hence also $(\rho-\bar{\rho},u)$, vanish everywhere inside the backward light-cone of $(t_0,x_0)$. This proves causality of the dissipative relativistic Euler equations \eqref{rel_art_visc}. Similarly, one can show that two solutions of \ref{rel_art_visc} that agree on $I(0)$ will also agree on $I(t)$ for all $t\in [0,t_0]$, the backward light-cone of $(t_0,x_0)$.
\hfill $\Box$

\section{Dissipativity and Global-in-time Existence}  \label{Sec_gl-ex}

Fix an integer $s\geq2$.  Let $(\rho,u)$ be a local solution of the mixed first--second order system \eqref{rel_art_visc} in one spatial dimension, as furnished by Theorem \ref{Thm_existence_causality}.  We study this solution near a constant rest state $(\bar{\rho},0)$. Set again $\kappa= \bar{\rho}+p(\bar{\rho})$ and $\sigma=\sqrt{p'(\bar{\rho})}$, and assume $\bar{\rho}>0$. Assume there exists an open interval $I$ containing $\bar{\rho}$ such that $p\in C^{s+2}(I)$ and 
\begin{equation}\label{eq:sec8-structural-hypotheses} 
p(\rho)+\rho > 0 \quad \forall \rho\in I, 
\qquad \text{and} \qquad    0< \sigma <1.
\end{equation}

Recall from Section \ref{Sec_Cauchy} that in variables $U=(\rho,w,z,u)^{\mathsf T}$, for $w=\partial_tu$ and $z=\partial_x u$, equation \eqref{rel_art_visc} is equivalent to \eqref{Cauchy_first-order-sys_B}
\begin{equation}\label{gl-ex_PDE-B}
\partial_tU+\mathcal B(U)\partial_xU= \tilde{F}(U),
\qquad \tilde{F}(U)\equiv(A^0(U))^{-1}\mathcal L(U).
\end{equation}
Recall that the compatibility condition $z=u_x$ is propagated by \eqref{gl-ex_PDE-B}, as shown in Lemma \ref{Lemma_Cauchy_first-order-sys}. Introducing now the constant state $\overline U=(\bar{\rho},0,0,0)^{\mathsf T}$, and the new variables
\begin{equation}\nonumber
W \equiv U-\overline U =(\varrho,w,z,u)^{\mathsf T},
\end{equation}
for $\varrho \equiv \rho-\bar{\rho}$, it follows that \eqref{gl-ex_PDE-B} can be written equivalently as 
\begin{equation}  \label{gl-ex_PDE-W}
W_t+B(W)W_x=F(W),
\end{equation}  
for $B(W)\equiv\mathcal B(\overline U+W)$ and $F(W)\equiv \tilde F(\overline U+W)$. Corollary \ref{Cauchy_Cor_exist_const-states}, applied to compatible $H^s$ initial data (i.e, $z=u_x$ holds at $t=0$), gives a local solution
\begin{equation}\nonumber
W\in C([0,T];H^{s}(\mathbb R;\mathbb R^4)) \cap C^1([0,T];H^{s-1}(\mathbb R;\mathbb R^4)),
\end{equation}
which is compatible in the sense that $z=u_x$ for all $t\in [0,T]$. We can now state our  global-in-time existence theorem for \eqref{gl-ex_PDE-W}, the main result of this section.

\begin{Thm} \label{Thm_gl-ex_W}
Let $s\geq2$ be an integer, and assume \eqref{eq:sec8-structural-hypotheses} and $\bar{\rho}>0$.  Let
\begin{equation}\nonumber
W_0=(\varrho_0,w_0,z_0,u_0)^{\mathsf T}
\in H^{s}(\mathbb R;\mathbb R^4)\cap L^1(\mathbb R;\mathbb R^4)
\end{equation}
be compatible initial data, (that is, $z_0=\partial_x u_0$). There exists $\delta_0>0$ such that, if
\begin{equation}\nonumber
\delta\equiv\|W_0\|_{L^1}+\|W_0\|_{H^{s}} \leq \delta_0,
\end{equation}
then the local solution of \eqref{gl-ex_PDE-W} extends uniquely to a global solution 
\begin{equation}\nonumber
W\in C\big([0,\infty);H^{s}\big) \cap C^1\big([0,\infty);H^{s-1}\big),
\end{equation}
which satisfies $z=u_x$ and $\rho(t) >0$ for all $t\geq 0$.  Moreover, the solution is dissipative in the sense that
\begin{equation}\label{eq:sec8-decay}
\|W(t)\|_{H^s}\leq C\delta(1+t)^{-1/4},\qquad \forall \; t\geq0.
\end{equation}
The constants $C$ and $\delta_0$ depend only on $s$, $\visc$, $\bar{\rho}$ and the pressure law on a fixed neighborhood of $\bar{\rho}$.
\end{Thm}

The main analytical challenge in the proof of Theorem \ref{Thm_gl-ex_W} is to establish Proposition \ref{Prop_gl-ex_bootstrap} below, whose proof we postpone to Subsections \ref{Sec_gl-ex_var} - \ref{Sec_gl-ex_energy}.  

\begin{Prop}\label{Prop_gl-ex_bootstrap}
Let $s\geq2$ and assume $W(t)$ is a compatible solution on $[0,T]$. Then there exist constants $\eta>0$ and $C>0$, independent of $T$, such that every compatible solution $W(t)$ of \eqref{gl-ex_PDE-W} on $[0,T]$, which satisfies 
\beq \label{gl-ex_bootstrap_X_def}
X_s(T)\leq \eta,
\qquad \text{for} \qquad 
X_s(T) \equiv \sup_{0\leq t\leq T}(1+t)^{1/4}\|W(t)\|_{H^s},
\eeq 
also satisfies
\begin{equation}\label{eq:sec8-bootstrap-combined}
X_s(T)\leq C \Big( \delta + X_s(T)^2 \Big).
\end{equation}
\end{Prop}

Assuming for the moment the validity of Proposition \ref{Prop_gl-ex_bootstrap}, (which will be proven in Subsection \ref{Sec_gl-ex_bootstrapProp}), we now complete the proof of Theorem \ref{Thm_gl-ex_W} by a bootstrap argument.

\Proof[Proof of Theorem \ref{Thm_gl-ex_W}]
Let $\eta>0$ and $C>0$ be the constants asserted to exist by Proposition \ref{Prop_gl-ex_bootstrap}, which we consider fixed. Without loss of generality, assume $C\geq 1$. Set 
\beq \label{gl-ex_Thm_proof_eqn1}
M \equiv 2C \qquad \text{and} \qquad \delta_0 \equiv \min \Big\{ \frac{\eta}{2C}, \frac{1}{8C^2} \Big\},
\eeq 
which by construction gives $M\delta_0 \leq \eta$ and $CM^2 \delta_0 \leq \frac{C}{2} = \frac{M}{4}$.

By Theorem \ref{Thm_existence_causality} there exists a unique local compatible solution $W(t)$ of \eqref{gl-ex_PDE-W} for each initial data with $\delta\equiv\|W_0\|_{L^1}+\|W_0\|_{H^{s}} \leq \delta_0$. We keep one such solution $W(t)$ fixed. Let $[0,T_{\max})$ denote the maximal interval of existence of $W(t)$, and denote with $X_s(T)$ the bootstrap norm of $W(t)$.  

By definition $X_s(0) \leq 2\delta \leq M \delta$, and since $X_s(T)$ is continuous in $T$, it follows that there exist some $T \in (0,T_{\max})$ such that
\beq \nonumber 
X_s(T)   \leq  M \delta .
\eeq 
Define $T_* \equiv \sup\left\{ T<T_{\max}: X_s(T) < M\delta \right\} >0$. 
Since $ M \delta \leq M \delta_0 \leq \eta$ by \eqref{gl-ex_Thm_proof_eqn1}, Proposition \ref{Prop_gl-ex_bootstrap} implies
$$
X_s(T)   \leq  C \delta + C X_s(T)^2 \qquad \forall\; T \in [0,T_*).
$$
Using now that our incoming bound \eqref{gl-ex_Thm_proof_eqn1} implies $C X_s(T)^2 \leq \delta  C M^2 \delta_0 \leq \frac{M}{4} \delta$, we obtain
\beq \label{gl-ex_Thm_proof_eqn2}
X_s(T)   \leq  C \delta + \frac{M}{4} \delta  \ = \ \frac{3}{4} M \delta,
\eeq 
for all $T \in [0,T_*)$. Thus our bootstrap assumption is maintained and by continuity of $X_s(T)$ it follows that $T_{*} =T_{\max}$. Namely, if $T_*<T_{\max}$, then continuity and \eqref{gl-ex_Thm_proof_eqn2} would imply $X_s(T)<M\delta$ on an interval slightly larger than $[0,T_*)$, contradicting the definition of $T_*$; hence $T_*=T_{\max}$.

Finally, to prove global-in-time existence, recalling $\frac{3}{4} M \delta = \frac{3}{2} C \delta$, observe that \eqref{gl-ex_Thm_proof_eqn2} implies the uniform bound
\beq \label{gl-ex_Thm_proof_eqn3}
\sup_{0\leq t < T_{\max}} \|W(t)\|_{H^s} \leq C \delta .
\eeq 
Thus, choosing $\delta>0$ sufficiently small, $W(t)$ stays within the neighborhood $\mathcal{D}$ in $H^s$ to which Corollary \ref{Cauchy_Cor_exist_const-states} applies, such that $\rho \geq \frac12 \bar{\rho}$ with $\rho \in I$. The continuation criterion for local solution of symmetric-hyperbolic equations therefore excludes $T_{\max}<\infty$, and implies $T_{\max}=\infty$, which is the sought-after global-in-time existence,  (see Proposition 1.5 in \cite[Ch.16]{Taylor3} or Theorem II in \cite{Kato}). The decay estimate \eqref{eq:sec8-decay} follows finally from \eqref{gl-ex_Thm_proof_eqn2} for all $t\geq 0$. Moreover, since by \eqref{gl-ex_Thm_proof_eqn3} $\|\rho(t)-\bar{\rho}\|_{L^\infty} \leq C\delta$, and since $\bar{\rho}>0$, it follows that choosing $\delta>0$ sufficiently small implies $\rho(t) > 0$ for all $t\geq 0$.   This completes the proof under the assumption that Proposition \ref{Prop_gl-ex_bootstrap} holds.
\QED 

It remains to prove Proposition \ref{Prop_gl-ex_bootstrap}, which is the subject of Subsections \ref{Sec_gl-ex_var} - \ref{Sec_gl-ex_energy}. Our proof is partially based on Duhamel's principle, which we implement in the next subsection.

\subsection{A nonlinear variable change and Duhamel's principle} \label{Sec_gl-ex_var}

For Duhamel's solution formula to entail the desired temporal decay, we need to transform equation \eqref{gl-ex_PDE-W} from $W$ to new variables $Y$, defined by the transformation 
\begin{equation} \label{gl-ex_Y-var_T(W)}
Y = \mathcal T(W) \equiv 
\begin{pmatrix}
T^{00}(\bar{\rho}+\varrho,u)+\epsilon v(u)w-\bar{\rho}\\
T^{01}(\bar{\rho}+\varrho,u)+\epsilon w\\
z\\u
\end{pmatrix},
\end{equation}
where $W=(\varrho,w,z,u)^{\mathsf T}$ and $\rho=\bar{\rho}+\varrho$. We denote the components of $Y$ by $Y \equiv (e,m,z,u)^{\mathsf T}$. The resulting equation is a balance law in conservation form as detailed in our first lemma.

\begin{Lemma} \label{gl-ex_Lemma_diffeo} 
\textbf{(i)} In a sufficiently small neighborhood $\Omega_0$ of the equilibrium $W=0$, the transformation $\mathcal{T}$ is a smooth local diffeomorphism, and there exists a constant $C\geq1$ such that for every integer $0\leq k\leq s$,\footnote{We write the norm equivalence \eqref{gl-ex_eqn_diffeo-norm} as $\|W\|_{H^k} \asymp \|Y\|_{H^k}$.}
\begin{equation} \label{gl-ex_eqn_diffeo-norm}
C^{-1}\|W\|_{H^k}\leq \|Y\|_{H^k} \leq C\|W\|_{H^k}.
\end{equation}
The variables $\rho$ and $w$ can be reconstructed from $Y$ as 
\beq \label{gl-ex_diffeo-rho-w}
\rho(Y) \equiv \bar{\rho}+e-v(u)m \qquad \text{and} \qquad w(Y) \equiv \frac{1}{\epsilon} \big(m-T^{01}(\rho(Y),u)\big).
\eeq 

\noindent \textbf{(ii)} 
On compatible solutions $W$, the system \eqref{gl-ex_PDE-W} is equivalent to the balance law
\begin{equation}\label{gl-ex_PDE-Y}
Y_t+\partial_xG(Y)=H(Y)
\end{equation}
in $Y$ variables, where 
\begin{equation}\nonumber
G(Y) \equiv  \begin{pmatrix} T^{01}(\rho(Y),u)-\epsilon v(u)z\\ T^{11}(\rho(Y),u)-p(\bar{\rho})-\epsilon z\\-\omega(Y)\\ 0 \end{pmatrix}, \qquad
H(Y) \equiv \begin{pmatrix}0\\0\\0\\\omega(Y)\end{pmatrix}.
\end{equation}
Equation \eqref{gl-ex_PDE-Y} propagates the compatibility condition $z=u_x$ from initial data.  
\end{Lemma}

\Proof
\textbf{(i):} Recall that 
\begin{equation}\label{eq:sec8-T-components}
\begin{aligned}
T^{00}(\rho,u)&=\rho+\bigl(\rho+p\bigr)u^2,\\
T^{01}(\rho,u)&=\bigl(\rho+p\bigr)u\sqrt{1+u^2},\\
T^{11}(\rho,u)&=p+\bigl(\rho+p\bigr)u^2,
\end{aligned}
\end{equation}
where $p=p(\rho)$. Thus, at $(\rho,u)=(\bar{\rho},0)$, we have $T^{00}_{\rho}=1$, $T^{01}_{u}=\kappa$ and $T^{00}_{u}=T^{01}_{\rho}=v=0$. Thus
\begin{equation}\label{eq:sec8-P-linear}
D\mathcal T(0) = \begin{pmatrix} 1&0&0&0\\ 0&\epsilon&0&\kappa\\ 0&0&1&0\\ 0&0&0&1 \end{pmatrix}, 
\qquad \det P=\epsilon>0,
\end{equation}
and the inverse function theorem implies that $\mathcal T$ is a local diffeomorphism at the origin.  The inverse formulas \eqref{gl-ex_diffeo-rho-w} follow from our definition $Y^1 = e \equiv T^{00} +\epsilon v(u)w-\bar{\rho}$ and $Y^2 = m \equiv T^{01}+\epsilon w$, which directly imply the sought-after identities 
\begin{equation}\nonumber
e + \bar{\rho} - v m = T^{00} - v(u)T^{01}=\rho
\qquad \text{and} \qquad 
m-T^{01} = \epsilon w.
\end{equation}
 
To establish the norm equivalence \eqref{gl-ex_eqn_diffeo-norm}, note that by \eqref{gl-ex_Y-var_T(W)} and the inverse formulas \eqref{gl-ex_diffeo-rho-w}, $\mathcal{T}$ and $\mathcal{T}^{-1}$ are both smooth functions. Taylor expansion yields 
\begin{equation}\nonumber
Y\equiv \mathcal{T}(W) =PW+R_2(W), 
\qquad \text{and} \qquad 
W\equiv \mathcal{T}^{-1}(Y)= P^{-1}Y+\widetilde R_2(Y),
\end{equation}
where $P\equiv D\mathcal T(0)$, and $R_2$ and $\widetilde R_2$ are smooth functions that vanish to second order at the origin. Standard estimates similar to those in Lemma \ref{Cauchy_Lemma_uniform-bounds}, using Morrey's inequality to control products, yield the desired norm equivalence \eqref{gl-ex_eqn_diffeo-norm}. Note in addition, for later use, that applying H\"older's inequality to the above Taylor expansion gives    
\beq  \label{gl-ex_diffeo_norms_data}
\|Y\|_{L^1} \leq C\bigl(\|W\|_{L^1}+\|W\|_{L^2}^2\bigr)
\quad \text{and} \quad 
\|W\|_{L^1} \leq C\bigl(\|Y\|_{L^1}+\|Y\|_{L^2}^2\bigr),
\eeq
in a potentially smaller neighborhood $\Omega_0$. This completes the proof of part \textbf{(i)}.

\smallskip 
\noindent \textbf{(ii):}
To begin, write \eqref{rel_art_visc} as
\begin{equation} \label{eq:sec8-mixed-system}
\begin{aligned}
\partial_tT^{00}(\rho,u)+\partial_xT^{01}(\rho,u)
  &=\epsilon\bigl(-\partial_t^2+\partial_x^2\bigr)\gamma(u),\\
\partial_tT^{01}(\rho,u)+\partial_xT^{11}(\rho,u)
  &=\epsilon\bigl(-\partial_t^2+\partial_x^2\bigr)u.
\end{aligned}
\end{equation}
Since $w=u_t$, $z=u_x$, $\gamma(u)=\sqrt{1+u^2}$ and $v(u)=\frac{u}{\gamma(u)}$, and since $\gamma'(u)=v(u)$ and $\gamma''(u)=v'(u)=\gamma(u)^{-3}$, we obtain that
\begin{equation}\nonumber
\bigl(-\partial_t^2+\partial_x^2\bigr)\gamma(u)
=v(u)(-w_t+z_x)+v'(u)(z^2-w^2).
\end{equation}
Substitution into the first equation of \eqref{eq:sec8-mixed-system} gives
\begin{equation}\nonumber
T^{00}_t+T^{01}_x+\epsilon v w_t-\epsilon v z_x +\epsilon v'(u)(w^2-z^2)=0,
\end{equation}
and since $\partial_t(\epsilon v w)=\epsilon v w_t+\epsilon v'w^2$ and $\partial_x(-\epsilon v z)=-\epsilon v z_x-\epsilon v'z^2$, we obtain
\beq \nonumber
\partial_t \big( T^{00} + \epsilon v w \big) + \partial_x \big( T^{01} - \epsilon v z  \big) =0.
\eeq 
This is the first conservation law in \eqref{gl-ex_PDE-Y}, keeping in mind that $Y^1 = e = T^{00} +\epsilon v(u)w-\bar{\rho}$. The second equation in \eqref{eq:sec8-mixed-system} is equivalent to
\begin{equation}\nonumber
\partial_t\bigl(T^{01}+\epsilon w\bigr) +\partial_x\bigl(T^{11}-\epsilon z\bigr)=0,
\end{equation}
which, since $Y^2 = m = T^{01}(\bar{\rho}+\varrho,u)+\epsilon w$, gives the second conservation law in \eqref{gl-ex_PDE-Y}. 

Finally, since $Y^3=z$ and $Y^4=u$, the third and fourth components of the balance law \eqref{gl-ex_PDE-Y} yield the conditions $z_t = w_x$ and $u_t=w$, which is identical the third and fourth components in \eqref{gl-ex_PDE-W}, as can be verified by the explicit formulas of the coefficients in \eqref{Cauchy_def_L(U)}, \eqref{eq:sec8-A0-inverse} and \eqref{eq:sec8-B-matrix}. Indeed, $\mathcal{T}$ acts as the identity on $z$ and $u$. Thus, on compatible solutions, \eqref{gl-ex_PDE-W} and \eqref{gl-ex_PDE-Y} are equivalent, and both equaitons propagate the compatibility conditions $z=u_x$ in the sense that if $z=u_x$ holds at $t=0$, then $z=u_x$ will also hold on $[0,T]$, the temporal domain of existence.
\QED

To introduce Duhamel's solution formula, we introduced the linear operator $L$ of \eqref{gl-ex_PDE-Y} and the associated solution operator $S$ of the linearization of \eqref{gl-ex_PDE-Y} as
\begin{equation}\label{eq:sec8-linear-L}
L=-A_*\partial_x+B_*  \qquad \text{and} \qquad S(t)=e^{tL},
\end{equation}
where $A_*\equiv DG(0)$ and $B_*\equiv DH(0)$, for
\small 
\begin{equation}\nonumber
DG(0)= \begin{pmatrix} 0&0&0&\kappa\\ \sigma^2&0&-\epsilon&0\\ 0&-\epsilon^{-1}&0&\kappa\epsilon^{-1}\\ 0&0&0&0 \end{pmatrix}, \qquad
DH(0)= \begin{pmatrix} 0&0&0&0\\ 0&0&0&0\\ 0&0&0&0\\ 0&\epsilon^{-1}&0&-\kappa\epsilon^{-1} \end{pmatrix}.
\end{equation}
\normalsize 
We now write the system \eqref{gl-ex_PDE-Y} equivalently as 
\begin{equation} \label{gl-ex_PDE_LY} 
Y_t = LY + \mathcal{J} \Phi(Y),
\end{equation}
where 
\begin{equation} \label{eq:sec8-Phi-J}
\Phi(Y)=\begin{pmatrix} T^{01} -\kappa u-\epsilon v(u)z \\ T^{11} - p(\bar{\rho}) - \sigma^2 e \\ \frac{1}{\epsilon} \big( \kappa u - T^{01}\big) \end{pmatrix},
\qquad \text{and} \qquad 
\mathcal J\phi \equiv 
\begin{pmatrix}-\partial_x\phi_1\\-\partial_x\phi_2\\
\partial_x\phi_3\\\phi_3\end{pmatrix}
\end{equation}
defines an operator $\R^3$-valued functions $\phi=(\phi_1,\phi_2,\phi_3)^T$ to $\R^4$-valued functions. Note that we have $\Phi(0)=0$ and $D\Phi(0) = 0$.

\begin{Lemma} \label{lem:sec8-exact-Y-equation}
For compatible solution (satisfying $z=u_x$), \eqref{gl-ex_PDE-Y} and \eqref{gl-ex_PDE_LY} are equivalent. Moreover, the compatibility condition is propagated by \eqref{gl-ex_PDE_LY}.
\end{Lemma}

\Proof 
Setting $\Phi(Y) = \big( \Phi_1(Y),\Phi_2(Y),\Phi_3(Y)\big)^T$, we have 
\beq \nonumber
G(Y) = A_*Y + \begin{pmatrix}\Phi_1(Y) \\ \Phi_2(Y) \\  -\Phi_3(Y) \\ 0 \end{pmatrix}
\qquad \text{and} \qquad 
H(Y)= B_*Y  + \begin{pmatrix}0\\0\\0\\ \Phi_3(Y)\end{pmatrix}.
\eeq 
Substituting these identities into \eqref{gl-ex_PDE-Y} gives 
\eqref{gl-ex_PDE_LY}.  The compatibility statement follows directly
from the definition of $\mathcal J$.  
\QED 

Duhamel's formula for \eqref{gl-ex_PDE_LY} is therefore
\begin{equation}\label{eq:sec8-Duhamel}
Y(t)= S(t)Y_0+ \int_0^tS(t-\tau)\mathcal J\Phi(Y(\tau))\,d\tau.
\end{equation}
The remainder of the proof has two distinct parts.  First, the low-frequency
spectrum of $L$ yields the heat-like decay rate $(1+t)^{-1/4}$, while all
frequencies bounded away from zero decay exponentially.  Second, the quadratic
vanishing of $\Phi$ and the compatible form of $\mathcal J$ provide one effective
spatial derivative in every nonlinear Duhamel term.  Duhamel’s formula will be used only after projection onto the low-frequency region. The high-frequency part of the solution will be controlled directly at order $H^s$ by the compensated energy estimate of Subsection \ref{Sec_gl-ex_energy}.

\subsection{Linear estimates}  \label{subsec:sec8-linear-estimates}

The linearization of equation \eqref{gl-ex_PDE-W} in $W$ at the rest state is
\begin{equation}  \label{gl-ex_PDE_lin-W} 
W_t=L_W W, \qquad \text{for} \quad 
L_W= \begin{pmatrix} 0&0&-\kappa&0\\ -\epsilon^{-1}\sigma^2\partial_x&-\epsilon^{-1}\kappa&\partial_x&0\\ 0&\partial_x&0&0\\ 0&1&0&0 \end{pmatrix}.
\end{equation}
The linear operators $L_W$ and $L = -A_*\partial_x+B_*$ in \eqref{eq:sec8-linear-L} are conjugate, on the compatible subspace $z=u_x$, by the constant matrix $P\equiv D\mathcal{T}(0)$ in \eqref{eq:sec8-P-linear}, 
\begin{equation}\label{eq:sec8-linear-conjugacy}
L=P L_W P^{-1},\qquad S(t)=P e^{tL_W}P^{-1}.
\end{equation}
On the compatible subspace the trivial mode (with eigenvalue $0$) of \eqref{gl-ex_PDE_lin-W} is ruled out. So, in order to bypass this trivial mode in our subsequent computations, we pass to two equivalent formulation of \eqref{gl-ex_PDE_lin-W} on compatible solutions, reduced to the relevant $3$-component subspace and converted to Fourier space.

\smallskip \smallskip \noindent \textbf{Formulation 1:}
For compatible solutions $W=(\varrho,w,z,u)^{\mathsf T}$ of \eqref{gl-ex_PDE_lin-W}, set
\begin{equation}\nonumber
V=(\varrho,w,u)^{\mathsf T},
\end{equation}
and keep $z=u_x$ as a (trivial) constraint. Then, the linearization of the original equation \eqref{rel_art_visc}, given in \eqref{linearized_RAV}, can be written equivalently as 
\beq \nonumber
\begin{cases}
\varrho_t+\kappa u_x&=0,\\
w_t+\frac{\kappa}{\epsilon}w+
\frac{\sigma^2}{\epsilon}\varrho_x-u_{xx}&=0,\\
u_t&=w.
\end{cases}
\eeq 
where $\kappa= \bar{\rho}+p(\bar{\rho})$ and $\sigma=\sqrt{p'(\bar{\rho})}$. With the Fourier convention
$\widehat f(\xi)= \frac{1}{\sqrt{2\pi}}\int_{\mathbb R}e^{-ix\xi}f(x)\,dx$, this becomes
\beq \label{gl-ex_PDE_lin-V}
\partial_t\widehat V=M_V(\xi)\widehat V,
\qquad \text{for} \quad 
M_V(\xi)= \begin{pmatrix} 0&0&-i\kappa\xi\\ -i\epsilon^{-1}\sigma^2\xi&-\tilde{\kappa}&-\xi^2\\ 0&1&0 \end{pmatrix},
\eeq 
where $\tilde{\kappa} \equiv \frac{\kappa}{\epsilon}$. We will use the formulation \eqref{gl-ex_PDE_lin-V} in the small frequency regime. 

\smallskip \smallskip \noindent \textbf{Formulation 2:}
Near large frequencies it is more convenient to use the variable
\begin{equation}\nonumber
Q=(\varrho,w,z)^{\mathsf T},
\end{equation}
and recover $u$ by solving $u_t=w$ with compatible data $u(0,x)=u_0(x)$ (i.e., data satisfying $\partial_x u_0 = z_0$).  For $\xi\neq0$, the variables $Q$ and $V$ are related in Fourier space by $\widehat Q=T(\xi)\widehat V$ in terms of the linear transformation $T(\xi)=\operatorname{diag}(1,1,i\xi)$. Thus, introducing 
\beq \nonumber
M_Q(\xi) \equiv T(\xi)M_V(\xi)T(\xi)^{-1} =
\begin{pmatrix} 
0&0&-\kappa\\ -i\epsilon^{-1}\sigma^2\xi&-\tilde{\kappa}&i\xi\\ 0&i\xi&0
\end{pmatrix},
\eeq 
we obtain the linear equation
\begin{equation} \label{gl-ex_PDE_lin-Q}
\partial_t\widehat Q=M_Q(\xi)\widehat Q,
\end{equation}
which is equivalent to \eqref{gl-ex_PDE_lin-V} on the compatible subspace.

\smallskip 
Both systems, \eqref{gl-ex_PDE_lin-Q} and \eqref{gl-ex_PDE_lin-V} have the common characteristic polynomial
\begin{equation}\label{eq:sec8-characteristic-polynomial}
p(\lambda,\xi)=\lambda^3+\tilde{\kappa}\lambda^2+ \xi^2\lambda+\tilde{\kappa}\sigma^2\xi^2.
\end{equation}
We denote the associated eigenvalues $\lambda_\pm(\xi)$ and  $\lambda_f(\xi)$, where $\lambda_\pm(0)=0$ and $\lambda_f(0)=-\tilde{\kappa} $. In the next three lemmas we establish estimates on these eigenvalues at small, large and intermediate frequencies, from which we conclude with the dissipation estimate on the linear solutions in Proposition \ref{prop:sec8-linear-decay}.

\begin{Lemma}[Small frequencies] \label{lem:sec8-small-frequency}
There exists $r_0>0$ such that for all $|\xi|\leq r_0$, 
\begin{align}
\lambda_\pm(\xi)
&=\pm i\sigma\xi- \frac{1-\sigma^2}{2\tilde{\kappa}}\xi^2+O(|\xi|^3), \label{eq:sec8-small-lambda-pm}\\
\lambda_f(\xi)
&=-\tilde{\kappa}+\frac{1-\sigma^2}{\tilde{\kappa}}\xi^2+O(|\xi|^4). \label{eq:sec8-small-lambda-f}
\end{align}
After possibly decreasing $r_0$, there exists constants $c_0>0$ and $\theta_0>0$ such that
\begin{equation}\nonumber
\operatorname{Re}\lambda_\pm(\xi)\leq-c_0\xi^2,
\qquad
\operatorname{Re}\lambda_f(\xi)\leq-\theta_0,
\end{equation}
and such that the spectral projectors of $M_V$ are uniformly bounded for $|\xi|\leq r_0$, and extend smoothly to $\xi=0$.
\end{Lemma}

\Proof 
Since the root $\lambda_f(0)=-\tilde{\kappa}$ of $p(\lambda,0)$ is simple, and since the implicit function theorem applied at $\xi=0$ and $\lambda_f(0)=-\tilde{\kappa}$ gives the existence of a smooth function $\lambda_f(\xi)$ solving the characteristic equation \eqref{eq:sec8-characteristic-polynomial}, we make the ansatz 
$$
\lambda_f(\xi) = -\tilde{\kappa} + \tilde{\nu} \xi^2 + O(|\xi|^4).
$$ 
Substituting this ansatz into \eqref{eq:sec8-characteristic-polynomial} and setting the resulting polynomial to zero for all $\xi$ sufficiently small, gives $\tilde{\kappa}^2\tilde{\nu}-\tilde{\kappa}(1-\sigma^2)=0$, and hence $\tilde{\nu} = \frac{1-\sigma^2}{\tilde{\kappa}}$, which proves \eqref{eq:sec8-small-lambda-f}.

For $\lambda_\pm(\xi)$, the two roots bifurcating from zero, we make the ansatz
$$
\lambda(\xi)= \mu \xi + \nu\xi^2 + O(|\xi|^3)
$$ 
for some $\mu,\nu\in\mathbb{C}$; (the existence of these smooth branches $\lambda(\xi)$, and resulting consistency with our ansatz, is implied by the implicit function theorem applied to $\tilde{p}(\zeta,\xi) \equiv p(\zeta \xi,\xi)/\xi^2)$ at $\zeta=\pm i \sigma$, where $\partial_\zeta \tilde{p}(\zeta,\xi) \neq 0$). Substitution into \eqref{eq:sec8-characteristic-polynomial} gives
\[
p(\lambda(\xi),\xi) = \tilde{\kappa}(\mu^2+\sigma^2)\xi^2 + \bigl(\mu^3+2\tilde{\kappa}\mu\nu+\mu\bigr)\xi^3 +O(|\xi|^4).
\]
For $p(\lambda(\xi),\xi)$ to vanish for all $\xi$ sufficiently close to zero, it follows from the $\xi^2$ coefficient that $\mu = \pm i \sigma$, and from the $\xi^3$ coefficient that $\nu=-\frac{1-\sigma^2}{2\tilde{\kappa}}$, which proves \eqref{eq:sec8-small-lambda-pm}. The real-part estimates now follow directly, keeping in mind that $0<\sigma^2<1$ by assumption.  

The eigenvectors of $M_V(\xi)$, for any nonzero eigenvalue $\lambda$, are all given up to normalization by 
\begin{equation}\nonumber
r(\lambda,\xi)=
\begin{pmatrix}-i\kappa\xi/\lambda\\\lambda\\1\end{pmatrix},
\end{equation}
as can be seen by direct inspection and use of the characteristic polynomial \eqref{eq:sec8-characteristic-polynomial}. Thus, eigenvectors associated to $\lambda_\pm(\xi)$ and $\lambda_f(\xi)$ are given by 
\small 
\beq  \nonumber 
r_+(\xi)=\begin{pmatrix}-\kappa/\sigma+O(\xi)\\ i\sigma\xi+O(\xi^2)\\1\end{pmatrix}, \quad 
r_-(\xi)=\begin{pmatrix}\kappa/\sigma+O(\xi)\\ -i\sigma\xi+O(\xi^2)\\1\end{pmatrix}, \quad
r_f(\xi) =\begin{pmatrix}i\epsilon\xi+O(\xi^3)\\ -\tilde{\kappa}+O(\xi^2)\\1\end{pmatrix}.
\eeq 
\normalsize 
At $\xi=0$, these eigenvectors are evidently linearly independent. Hence these eigenvectors are linearly independent for $\xi$ sufficiently close to $0$. The corresponding projectors therefore are uniformly bounded for $|\xi| \leq r_0$ and $r_0 > 0$ sufficiently small, and  extend smoothly to $\xi=0$.
\QED 

We next establish an analogous result for high frequencies. Even though this result is not directly used in the proof of Theorem \ref{Thm_gl-ex_W}--high frequency terms are controlled by an energy estimate--, we require it for the full linear estimates of Proposition \ref{prop:sec8-linear-decay} which we include for completeness.

\begin{Lemma}[Large frequencies]\label{lem:sec8-large-frequency}
There exists $R_0>0$ such that, for $|\xi|\geq R_0$, the roots of
\eqref{eq:sec8-characteristic-polynomial} satisfy (potentially after relabelling)
\begin{align}
\lambda_\pm(\xi) &=\pm i|\xi|-\frac{\tilde{\kappa}}{2}(1-\sigma^2) +O(|\xi|^{-1}),
\label{eq:sec8-large-lambda-pm}\\
\lambda_f(\xi)&=-\tilde{\kappa}\sigma^2+O(|\xi|^{-1}). \label{eq:sec8-large-lambda-f}
\end{align}
Moreover, after potentially increasing $R_0$, there exists a constant $\theta_1>0$ such
that
\begin{equation}\label{eq:sec8-large-real-parts}
\max_{j\in\{+, -, f\}}\operatorname{Re}\lambda_j(\xi)
\leq-\theta_1, \qquad |\xi|\geq R_0,
\end{equation}
and the spectral projectors of $M_Q(\xi)$ are uniformly bounded for all $|\xi|\geq R_0$.
\end{Lemma}

\Proof 
Substituting the ansatz $\lambda=\mu |\xi| + \nu +O(|\xi|^{-1})$ into \eqref{eq:sec8-characteristic-polynomial}, the coefficients of $|\xi|^3$ and $|\xi|^2$ give, respectively,
\[
\mu(\mu^2+1)=0,
\qquad
(3\mu^2+1)\nu+\tilde{\kappa}\mu^2+\tilde{\kappa}\sigma^2=0.
\]
Solving the first equation with $\mu=\pm i$, the second equation yields $\nu=-\tilde{\kappa}(1-\sigma^2)/2$, which gives $\lambda_\pm(\xi)$ in \eqref{eq:sec8-large-lambda-pm}.  The root $\lambda_f(\xi)$ is obtained by taking $\mu=0$, in which case the second equation above yields $\nu=-\tilde{\kappa}\sigma^2$, which proves \eqref{eq:sec8-large-lambda-f}. The bounds \eqref{eq:sec8-large-real-parts} follows directly from \eqref{eq:sec8-large-lambda-pm} -- \eqref{eq:sec8-large-lambda-f}.

For an eigenvalue $\lambda$ of $M_Q(\xi)$, the corresponding eigenvector (up to a factor) is
\begin{equation}\nonumber
r(\lambda,\xi)= \begin{pmatrix}-\kappa/\lambda\\ \lambda/(i\xi)\\1\end{pmatrix},
\end{equation}
so as $|\xi|\to\infty$ we have
\beq 
r_\pm(\xi) \longrightarrow \begin{pmatrix}0\\ \pm \operatorname{sgn}\xi\\1\end{pmatrix}
\qquad \text{and} \qquad
r_f(\xi)\longrightarrow \begin{pmatrix}\epsilon/\sigma^2\\0\\1\end{pmatrix}.
\eeq 
Thus, the limiting eigenvectors are linearly independent, which implies that the associated spectral projectors are uniformly bounded for $|\xi| \geq R_0$ sufficiently large.
\QED 

We now address the intermediate frequency region, which is crucial for the proof of Proposition \ref{Prop_gl-ex_bootstrap} and Theorem \ref{Thm_gl-ex_W}. Since eigenvalues may fail to be simple in this region, decomposing the flow operator into spectral projectors may no longer yield uniform bounds. For this reason, we estimate the flow operator directly without studying the eigenvalues and their spectral projectors.

\begin{Lemma}[Intermediate frequencies]\label{lem:sec8-middle-frequency}
For every $0<r_0<R_0$ there exist constants $C>0$ and $\theta_m>0$ such that
\begin{equation} \label{eq:sec8-middle-semigroup}
\|e^{tM_Q(\xi)}\|\leq Ce^{-\theta_m t}, \qquad r_0\leq|\xi|\leq R_0,
\end{equation}
where $\|\cdot\|$ denotes a matrix norm (e.g., the Hilbert-Schmidt norm).
\end{Lemma}

\Proof 
By Lemma \ref{Lemma_decay}, every root of
\[
p(\lambda,\xi) = \lambda^3+\tilde{\kappa}\lambda^2+\xi^2\lambda
+\tilde{\kappa}\sigma^2\xi^2
\]
has strictly negative real part whenever $\xi\neq0$.  Alternatively, this
follows from the Routh--Hurwitz criterion, since $\tilde{\kappa}>0$, $\xi^2>0$ and  $\tilde{\kappa} \sigma^2\xi^2>0$, as well as $\tilde{\kappa}\xi^2-\tilde{\kappa}\sigma^2\xi^2 =\tilde{\kappa}(1-\sigma^2)\xi^2>0$.

By compactness of the annulus $K_{r_0,R_0} := \{\xi\in\mathbb R:r_0\leq|\xi|\leq R_0\}$, the spectrum is contained in $\{\operatorname{Re}\lambda\leq-2\theta_m\}$ for some $\theta_m>0$. This can be seen by a contradiction argument. Namely, if no uniform spectral gap existed, there would be sequences $\xi_n\in K_{r_0,R_0}$ and eigenvalues $\lambda_n$ of $M_Q(\xi_n)$ such that $\operatorname{Re}(\lambda_n) \to 0$.  Passing to a subsequence, we may assume $\xi_n\to\xi_*\in K_{r_0,R_0}$ and $\lambda_n\to\lambda_*$. Continuity of the characteristic polynomial then implies that $\lambda_*$ is again an eigenvalue of $M_Q(\xi_*)$ with $\operatorname{Re}\lambda_*=0$, in contradiction to the strict negativity established above.  

To obtain a uniform semigroup bound without assuming that the eigenvalues are simple,\footnote{The characteristic polynomial \eqref{eq:sec8-characteristic-polynomial} may have repeated roots, for example, for  $\sigma^2=1/9$ and $|\xi|=\tilde{\kappa}/\sqrt3$, one has $p(\lambda,\xi)=\left(\lambda+\frac{\tilde{\kappa}}{3}\right)^3$.}  choose a fixed (positively oriented) rectangular contour $\Gamma$ enclosing the spectra of all matrices $M_Q(\xi)$, for $\xi\in K_{r_0,R_0}$, and which is contained in
$$
\{z\in\mathbb C:\operatorname{Re}z\leq-\theta_m\}.
$$
Since $\Gamma$ does not meet any of the eigenvalues of $M_Q(\xi)$, for $\xi\in K_{r_0,R_0}$, the map
$$
(z,\xi)\longmapsto(zI-M_Q(\xi))^{-1}
$$
is continuous on the compact set $\Gamma\times K_{r_0,R_0}$, and hence uniformly bounded,
$$
C_R \equiv  \sup_{\substack{z\in\Gamma\\ \xi\in K_{r_0,R_0}}} \|(zI-M_Q(\xi))^{-1}\| <\infty.
$$
By the Dunford-Taylor formula (see (5.47) in \cite[Ch.~I]{Kato_Perturbation})--the matrix version of the Cauchy integral formula--we have 
$$
e^{tM_Q(\xi)} = \frac{1}{2\pi i} \int_\Gamma e^{tz}(zI-M_Q(\xi))^{-1}\,dz,
$$
which directly gives the estimate
\beq \nonumber 
\|e^{tM_Q(\xi)}\| \leq
\frac{|\Gamma|}{2\pi}
\sup_{z\in\Gamma}e^{t\operatorname{Re}z}
\sup_{\substack{z\in\Gamma\\
\xi\in K_{r_0,R_0}}}
\|(zI-M_Q(\xi))^{-1}\|
\leq
\frac{|\Gamma|C_R}{2\pi}e^{-\theta_m t},
\eeq 
where $|\Gamma|$ denotes the (fixed) arc-length of $\Gamma$. This completes the proof.
\QED

Let us remark that combining Lemmas \ref{lem:sec8-small-frequency} -- \ref{lem:sec8-middle-frequency}, choosing $r_0>0$ sufficiently small and $R_0>r_0$ sufficiently large, one obtains the spectral bound
\begin{equation}\nonumber
\max_{\lambda\in\operatorname{spec}M_V(\xi)} \operatorname{Re}\lambda \leq -c\frac{\xi^2}{1+\xi^2},\qquad \xi\neq0,
\end{equation}
where $M_Q$ is used for large and intermediate frequencies, and $M_V$ for small frequencies.

Based on Lemmas \ref{lem:sec8-small-frequency} -- \ref{lem:sec8-middle-frequency} we now establish the full linear decay estimates. We include this here for completeness. The proof of Proposition \ref{Prop_gl-ex_bootstrap} and Theorem \ref{Thm_gl-ex_W} only requires this linear decay estimate for small and intermediate frequencies.

\begin{Prop} \label{prop:sec8-linear-decay}
Let $s\geq 2$. There exist constants $\theta>0$ and $C>0$ such that every compatible solutions $Y$ of the linearized $Y$-equation $Y_t=LY$ with compatible initial data $Y_0\in H^s(\mathbb R)\cap L^1(\mathbb R)$, satisfies 
\begin{equation} \label{eq:sec8-linear-semigroup-Hs}
\|S(t)Y_0\|_{H^s} \leq C_s(1+t)^{-1/4}\|Y_0\|_{L^1}  +C_se^{-\theta t}\|Y_0\|_{H^s},
\end{equation}
where $S(t)=e^{tL}$ and $L$ is defined in \eqref{gl-ex_PDE_LY}.  
\end{Prop}

\Proof 
We first prove the more refined result that there exists a constant $\theta>0$ such that, for every integer $0\leq k \leq s$, there is a constant $C_k>0$ for which every compatible solution of the linearized $W$-equation \eqref{gl-ex_PDE_lin-W} satisfies
\begin{align}
\|\partial_x^kW(t)\|_{L^2} &\leq C_k(1+t)^{-k/2-1/4}\|W_0\|_{L^1}   +C_ke^{-\theta t}\|W_0\|_{H^k}, \label{eq:sec8-linear-decay-L1} 
\end{align}
where $W(t) = S_W(t) W_0$. By Lemma \ref{lem:sec8-exact-Y-equation}, the same estimates then hold with equivalent norms for the compatible solutions of the linearized $Y$-equation $Y_t=LY$. More precisely, on the compatible subspace equation \eqref{eq:sec8-linear-conjugacy} gives
$$
S(t)=P S_W(t)P^{-1}, \qquad \text{and} \qquad  S_W(t)=e^{tL_W}.
$$
where $P=D\mathcal T(0)$. Since $P$ is a constant invertible matrix, it induces equivalent $L^1$ and $H^k$ norms.  Therefore, once \eqref{eq:sec8-linear-decay-L1} is established for $W$, an analogous estimate will hold with modified constants for $S(t)Y_0$, and summing this estimate over $0\leq k\leq s$, using $(1+t)^{-k/2-1/4}\leq(1+t)^{-1/4}$, then gives the sought-after estimate \eqref{eq:sec8-linear-semigroup-Hs}. So it remains to prove \eqref{eq:sec8-linear-decay-L1} for the $W$-variables.

The linear flow $S_W(t)=e^{tL_W}$ of compatible solutions to the linear $W$-equation \eqref{gl-ex_PDE_lin-W}, after Fourier transformation, is equivalently given by
$$
\widehat V(t,\xi)=e^{tM_V(\xi)}\widehat V_0(\xi),
\qquad
V=(\varrho,w,u)^{\mathsf T},
$$
for $|\xi|\in [0,R_0)$, and for $|\xi| \geq R_0$ sufficiently large by
$$
\widehat Q(t,\xi)=e^{tM_Q(\xi)}\widehat Q_0(\xi),
\qquad
Q=(\varrho,w,z)^{\mathsf T}.
$$
Thus our results on the matrix exponentials $e^{tM_V(\xi)}$ and $e^{tM_Q(\xi)}$ in Lemmas \ref{lem:sec8-small-frequency} -- \ref{lem:sec8-middle-frequency} directly carry over to estimates on $S_W(t)$ on the subspace of compatible solutions.

Fix $r_0>0$ such that Lemma \eqref{lem:sec8-small-frequency} applies for all $|\xi|\leq r_0$. We then have for $|\xi|\leq r_0$ the spectral decomposition
$$
e^{tM_V(\xi)} =\sum_{j\in\{+,-,f\}}e^{t\lambda_j(\xi)}\Pi_j^V(\xi).
$$
By Lemma \ref{lem:sec8-small-frequency}, the spectral projectors $\Pi_j^V(\xi)$ are uniformly bounded, and the associated eigenvalues satisfy $\operatorname{Re}\lambda_\pm(\xi)\leq-c_0\xi^2$ and $\operatorname{Re}\lambda_f(\xi)\leq-\theta_0$.
As a consequence we obtain in the small frequency regime that 
\beq \label{gl-ex_Prop_linear-decay_eqn1}
\|e^{tM_V(\xi)}\| \leq C\bigl(e^{-c_0\xi^2t}+e^{-\theta_0t}\bigr), \qquad |\xi|\leq r_0.
\eeq 

Next fix $R_0>r_0$ so that Lemma \ref{lem:sec8-large-frequency} applies for $|\xi|\geq R_0$.  The uniform spectral projector bounds  of Lemma \ref{lem:sec8-large-frequency} together with the estimate $\operatorname{Re}\lambda_j(\xi)\leq-\theta_1$  give
$$
\|e^{tM_Q(\xi)}\|\leq Ce^{-\theta_1t}, \qquad |\xi|\geq R_0.
$$
Now, on the compact annulus $r_0\leq|\xi|\leq R_0$, Lemma \eqref{lem:sec8-middle-frequency} readily provides the estimate  $\|e^{tM_Q(\xi)}\|\leq Ce^{-\theta_m t}$.  In combination,  we obtain
\beq \label{gl-ex_Prop_linear-decay_eqn2}
\|e^{tM_Q(\xi)}\|\leq Ce^{-\theta t}, \qquad |\xi|\geq r_0,
\eeq 
where $\theta \equiv \min\{\theta_1,\theta_m\}$. 

\smallskip 
We are now prepared to establish \eqref{eq:sec8-linear-decay-L1}. By Plancherel's theorem,
\begin{align} \label{gl-ex_Prop_linear-decay_eqn3}
\|\partial_x^kW(t)\|_{L^2}^2
&=\int_{|\xi|\leq r_0}|\xi|^{2k}|\widehat W(t,\xi)|^2\,d\xi  +\int_{|\xi|\geq r_0}|\xi|^{2k}|\widehat W(t,\xi)|^2\,d\xi.
\end{align}
We now estimate the two terms on the right separately. On the low-frequency region, $|\xi|\leq r_0$, the equivalence between the $W$- and $V$-representations together with \eqref{gl-ex_Prop_linear-decay_eqn1} imply (after potentially decreasing $\theta>0$)
\begin{align} \label{gl-ex_Prop_linear-decay_eqn3b}
|\widehat W(t,\xi)| &\leq C|\widehat V(t,\xi)| \cr 
& \leq C \|e^{tM_V(\xi)}\| \; |\widehat V_0(\xi)| \cr 
& \leq C \; \bigl(e^{-c_0\xi^2t}+e^{-\theta t}\bigr) \; |\widehat V_0(\xi)| .
\end{align}
On the large and intermediate frequency region, where $|\xi|\geq r_0$, the equivalence between the $W$- and $Q$-representations together with \eqref{gl-ex_Prop_linear-decay_eqn2} imply
\begin{align*}
|\widehat W(t,\xi)| 
\leq C \|e^{tM_Q(\xi)}\| \; |\widehat Q_0(\xi)|  
\leq C \; e^{-\theta t} |\widehat Q_0(\xi)|.
\end{align*}
Substituting the above two estimates into \eqref{gl-ex_Prop_linear-decay_eqn3}, using that $|\widehat Q_0(\xi)| \leq |\widehat W_0(\xi)|$ and $|\widehat V_0(\xi)| \leq |\widehat W_0(\xi)|$, we obtain that
\begin{align} \label{gl-ex_Prop_linear-decay_eqn4}
\|\partial_x^kW(t)\|_{L^2}^2
& \leq C \int_{|\xi|\leq r_0} |\xi|^{2k} e^{-2c_0\xi^2t}|\widehat W_0(\xi)|^2\,d\xi 
+ C e^{-2\theta t} \|W_0\|_{H^k}^2 ,
\end{align}   
where we used Plancherel's Theorem to obtain the last term. 

To obtain \eqref{eq:sec8-linear-decay-L1}, we bound the low frequency integral on the right hand side of \eqref{gl-ex_Prop_linear-decay_eqn4} by
\begin{align} \label{gl-ex_Prop_linear-decay_eqn5}
\int_{|\xi|\leq r_0} |\xi|^{2k} e^{-2c_0\xi^2t}|\widehat W_0(\xi)|^2\,d\xi 
& \leq \|W_0\|_{L^1}^2 \int_{|\xi|\leq r_0} |\xi|^{2k} e^{-2c_0\xi^2t}\,d\xi \cr 
& \leq C_k \|W_0\|_{L^1}^2 (1+t)^{-k-1/2},
\end{align}
using a Gaussian moment estimate for the last step, (as can be computed by a change of variables to $\eta=\sqrt t\,\xi$, and absorbing integration over $0\leq t\leq1$ into $C_k$). Substituting \eqref{gl-ex_Prop_linear-decay_eqn5} back into \eqref{gl-ex_Prop_linear-decay_eqn4}, taking a square root, gives
$$
\|\partial_x^kW(t)\|_{L^2} \leq C_k(1+t)^{-k/2-1/4}\|W_0\|_{L^1} +C_ke^{-\theta t}\|W_0\|_{H^k},
$$
which is the sought-after estimate \eqref{eq:sec8-linear-decay-L1}. 
\QED

\subsection{Nonlinear estimates}\label{subsec:sec8-nonlinear-estimates}

We now derive the estimates on the non-linear Duhamel terms required for the proof of Proposition \ref{Prop_gl-ex_bootstrap}. For this, we first introduce a suitable frequency cut-off. Fix an even function $\varphi\in C^\infty(\mathbb{R})$ such that
$$
0\leq\varphi\leq1,\qquad
\varphi(\xi)=0\ \text{for }|\xi|\leq1,\qquad
\varphi(\xi)=1\ \text{for }|\xi|\geq2.
$$
For $R>1$, define the high-frequency multiplier $P_{\mathrm{hi}}$ in Fourier space by
$$
\widehat{P_{\mathrm{hi}}f}(\xi) := \varphi(\xi/R)\widehat f(\xi),
$$
and set $P_{\mathrm{lo}}=I-P_{\mathrm{hi}}$. Our first lemma estimates the integrand in Duhamel's formula, which produces an additional low-frequency decay factor.

\begin{Lemma} \label{lem:sec8-structured-linear}
Let $\phi=(\phi_1,\phi_2,\phi_3)^{\mathsf T}$ and let $\mathcal J$ be defined in
\eqref{eq:sec8-Phi-J}.  Then
\begin{align} \label{eq:sec8-structured-low} 
\|P_{\mathrm{lo}}S(t)\mathcal J\phi\|_{H^s}
&\leq C(1+t)^{-3/4}\|\phi\|_{L^1}  +Ce^{-\theta t}\|\phi\|_{L^2}, 
\end{align}
for every integer $s\geq0$, where $C>0$ is a constant independent of $\phi$ and $t$.
\end{Lemma}

\Proof 
Let $\chi_R(\xi):=1-\varphi(\xi/R)$ be the symbol of $P_{\mathrm{lo}}$. Since $\chi_R$ is supported in $|\xi|\leq2R$, Plancherel's theorem gives
\beq \label{gl-ex_non-lin_eqn0}
\|P_{\mathrm{lo}}S(t)\mathcal J\phi\|_{H^s}^2 =
\int_{|\xi|\leq2R} \langle\xi\rangle^{2s}|\chi_R(\xi)|^2 \bigl|\widehat{S(t)\mathcal J\phi}(\xi)\bigr|^2\,d\xi,
\eeq 
where $\langle\xi\rangle^{2} \equiv (1+\xi^2)$. To estimate $\widehat{S(t)\mathcal J\phi}$ on the right hand side we pass to the reduced $V$- and $Q$-representation of the linearization of \eqref{gl-ex_PDE_LY} and apply Lemmas \ref{lem:sec8-small-frequency} - \ref{lem:sec8-middle-frequency}.

\smallskip\noindent \textit{Claim:}
The flow in the $V$-representation satisfies
\begin{equation} \label{eq:sec8-L89-WV-representation}
\widehat{S(t)\mathcal J\phi}(\xi) = P\mathsf E_V(\xi)e^{tM_V(\xi)}\widehat F_\phi^V(\xi),
\end{equation}
where $P=D\mathcal T(0)$ and 
$$
\mathsf E_V(\xi) \equiv \begin{pmatrix} 1&0&0\\ 0&1&0\\ 0&0&i\xi\\ 0&0&1 \end{pmatrix}
\quad \text{and} \quad
\widehat F_\phi^V(\xi) = -i\xi \begin{pmatrix} \widehat\phi_1(\xi)\\[1mm] \epsilon^{-1}\widehat\phi_2(\xi)\\ 0 \end{pmatrix}
+ \widehat\phi_3(\xi) \begin{pmatrix} 0\\-\tilde{\kappa}\\1 \end{pmatrix}. 
$$ 

\smallskip \noindent \textit{Proof of Claim:} 
Note first that $\mathcal J\phi = (-\partial_x \phi_1,\; -\partial_x\phi_2,\; \partial_x\phi_3,\; \phi_3)^T$ is compatible, since its third component is the $x$-derivative of its fourth component. Hence, the conjugacy between the linearized $Y$- and $W$-flows in \eqref{eq:sec8-linear-conjugacy} applies and gives
\begin{equation} \label{gl-ex_non-lin_eqn1}
S(t)  [\mathcal J\phi]  =  P e^{tL_W}P^{-1}  [\mathcal J\phi].
\end{equation}
A direct computation gives
\begin{equation}\nonumber
P^{-1}\mathcal J\phi =
\begin{pmatrix}
-\partial_x\phi_1\\[1mm]
-\epsilon^{-1}\partial_x\phi_2-\tilde{\kappa}\phi_3\\[1mm]
\partial_x\phi_3\\
\phi_3
\end{pmatrix},
\qquad \text{since} \quad 
P= \begin{pmatrix} 1&0&0&0\\ 0&\epsilon&0&\kappa\\0&0&1&0\\ 0&0&0&1 \end{pmatrix}.
\end{equation}
The $V$-representation of the flow on $\mathcal J \phi$ thus acts on the reduced variables 
$$
F_\phi^V \equiv ( -\partial_x\phi_1,\; -\epsilon^{-1}\partial_x\phi_2-\tilde{\kappa}\phi_3,\; \phi_3 )^T,
$$ 
whose Fourier transform is given by $\widehat F_\phi^V(\xi)$. The four components of a compatible $W$ are reconstructed from $V$ by
\begin{equation}\nonumber
\widehat W(\xi)=\mathsf E_V(\xi)\widehat V(\xi).
\end{equation}
Finally, since the reduced linearized flow is
\[
\widehat V(t,\xi)=e^{tM_V(\xi)}\widehat V(0,\xi),
\]
we obtain the exact representation
\begin{equation} \nonumber
\widehat{S(t)\mathcal J\phi}(\xi) = P\mathsf E_V(\xi)e^{tM_V(\xi)}\widehat F_\phi^V(\xi),
\end{equation}
which is \eqref{eq:sec8-L89-WV-representation}. This proves the claim. \hfill $\Diamond$

\smallskip 
We first estimate \eqref{eq:sec8-L89-WV-representation} for $|\xi|\leq r_0$. For this recall the spectral decomposition 
\beq \label{gl-ex_non-lin_eqn2}
e^{tM_V(\xi)} = e^{t\lambda_+(\xi)}\Pi_+^V(\xi) +e^{t\lambda_-(\xi)}\Pi_-^V(\xi) +e^{t\lambda_f(\xi)}\Pi_f^V(\xi).
\eeq 
The fast projector is uniformly bounded, and thus
\begin{equation} \label{gl-ex_non-lin_eqn3a}
\bigl|\Pi_f^V(\xi)\widehat F_\phi^V(\xi)\bigr| \leq C|\widehat\phi(\xi)|, \qquad |\xi|\leq r_0.
\end{equation}
To estimate the slow projector terms, observe that the vector  $r_f(0)\equiv  (0,\; -\tilde{\kappa},\; 1)^T$, appearing in the second term in $\widehat F_\phi^V(\xi)$, is an eigenvector of $M_V(0)$ corresponding to the eigenvalue $\lambda_f(0)=-\tilde{\kappa}$. Hence $\Pi_\pm^V(0)r_f(0)=0$,\footnote{Although the eigenvectors are not orthogonal, the spectral projectors satisfy $\Pi_j^V(\xi)\Pi_k^V(\xi)= \delta_{jk}\Pi_j^V(\xi)$ for $0<|\xi|\leq r_0$, (cf. (5.21) in \cite[Ch.~I]{Kato_Perturbation}). Since the projectors extend continuously to $\xi=0$, these identities remain valid at $\xi=0$.  In particular, $\Pi_\pm^V(0) \Pi_f^V(0)=0$. Since $r_f(0)\in\operatorname{Ran}\Pi_f^V(0)$, it follows that $\Pi_\pm^V(0) r_f(0) = \Pi_\pm^V(0)\Pi_f^V(0)r_f(0)=0$.}             
and since by Lemma \ref{lem:sec8-small-frequency} these projectors extend smoothly to $\xi=0$, with invertible eigenvector matrix, we have
\begin{equation}\nonumber
\|\Pi_\pm^V(\xi)r_f(0)\| \leq C|\xi|, \qquad |\xi|\leq r_0.
\end{equation}
This, since the first term of $\widehat F_\phi^V$ carries an explicit factor $\xi$, gives us 
\begin{equation} \label{gl-ex_non-lin_eqn3b}
\bigl|\Pi_\pm^V(\xi)\widehat F_\phi^V(\xi)\bigr|
\leq C|\xi|\,|\widehat\phi(\xi)|,
\qquad |\xi|\leq r_0.
\end{equation}
Combining now \eqref{gl-ex_non-lin_eqn3a} and \eqref{gl-ex_non-lin_eqn3b} with the spectral decomposition \eqref{gl-ex_non-lin_eqn2} and the estimates of Lemma \ref{lem:sec8-small-frequency}, we conclude with the sough-after estimate 
\begin{equation} \nonumber 
\bigl|e^{tM_V(\xi)}\widehat F_\phi^V(\xi)\bigr|
\leq C\bigl(|\xi|e^{-c\xi^2t}+e^{-\theta t}\bigr) |\widehat\phi(\xi)|, \qquad |\xi|\leq r_0.
\end{equation}
Thus, using that $P$ and $E_V$ are uniformly bounded, implies 
\beq  \label{gl-ex_non-lin_eqn4}
|\widehat{S(t)\mathcal J\phi}(\xi) | \leq C\bigl(|\xi|e^{-c\xi^2t}+e^{-\theta t}\bigr) |\widehat\phi(\xi)|, \qquad |\xi|\leq r_0,
\eeq 
which is the sought-after estimate.

We next estimate \eqref{eq:sec8-L89-WV-representation} on the annulus $r_0\leq|\xi|\leq2R$. On this fixed annulus the transformations between the compatible
$V$- and $Q$-representations (as well as $P$, $P^{-1}$ and $E_V$) are uniformly bounded. Thus Lemma \ref{lem:sec8-middle-frequency}, with $R_0=2R$, directly implies
\beq \label{eq:sec8-L89-annulus-flow}
\bigl|\widehat{S(t)\mathcal J\phi}(\xi)\bigr| \leq Ce^{-\theta t}|\widehat\phi(\xi)|, \qquad r_0\leq|\xi|\leq2R,
\eeq 
which is the sough-after estimate on the annulus.

Substituting now \eqref{gl-ex_non-lin_eqn4} and \eqref{eq:sec8-L89-annulus-flow} into the Plancherel identity \eqref{gl-ex_non-lin_eqn0}, and absorbing the bounded factor $\langle\xi\rangle^{2s}|\chi_R(\xi)|^2$ into the constant, gives
\begin{align*}
\|P_{\mathrm{lo}}S(t)\mathcal J\phi\|_{H^s}^2
\leq  C\int_{|\xi|\leq r_0} |\xi|^2e^{-2c\xi^2t}|\widehat\phi(\xi)|^2\,d\xi
+ Ce^{-2\theta t} \int_{|\xi|\leq2R}|\widehat\phi(\xi)|^2\,d\xi.
\end{align*}
Plancherel's theorem bounds the second integral by $C\|\phi\|_{L^2}^2$. Since $|\widehat\phi(\xi)|\leq\|\phi\|_{L^1}$, the first integral is bounded by
$$
C\|\phi\|_{L^1}^2 \int_{|\xi|\leq r_0}|\xi|^2e^{-2c\xi^2t}\,d\xi  
\leq  C \; \|\phi\|_{L^1}^2\; (1+t)^{-3/2}.
$$
Combining this and taking square roots yields
$$
\|P_{\mathrm{lo}}S(t)\mathcal J\phi\|_{H^s} \leq C(1+t)^{-3/4}\|\phi\|_{L^1} +Ce^{-\theta t}\|\phi\|_{L^2},
$$
which proves the lemma.
\QED

For the next lemma and the proof of Proposition \eqref{Prop_gl-ex_bootstrap} it is convenient to use the following bootstrap quantity in $Y$,
\begin{equation}\label{eq:sec8-XK-Y}
\mathcal{X}^Y_s(T)=\sup_{0\leq t\leq T}(1+t)^{1/4}\|Y(t)\|_{H^s},
\end{equation}
which by Lemma \ref{lem:sec8-exact-Y-equation} can be used equivalently for our purposes here.

\begin{Lemma} \label{lem:sec8-Duhamel-nonlinear}
Let $Y$ be a compatible solution on $[0,T]$, contained in the neighborhood $\Omega_0$ of Lemma \ref{gl-ex_Lemma_diffeo}.  Then, for $0\leq t\leq T$, $s\geq2$,
\beq \label{eq:sec8-Duhamel-nonlinear-bound}
I_{\mathrm{lo}}(t) \equiv \left\|P_\text{lo} \int_0^tS(t-\tau)\mathcal J\Phi(Y(\tau))\,d\tau\right\|_{H^s} 
\leq C \; \mathcal{X}^Y_s(T)^2 \; (1+t)^{-1/4}.
\eeq 
\end{Lemma}

\begin{proof}
To begin, observe that by Morrey's inequality ($s\geq2$) we have for any $0\leq\tau\leq T$ the bound
\beq \label{gl-ex_non-lin_L2_eqn1}
\|Y(\tau)\|_{L^2}+\|Y(\tau)\|_{L^\infty}  \leq C \mathcal{X}^Y_s(T)(1+\tau)^{-1/4}.
\eeq 
Moreover, since $\Phi(0)=0$ and $D\Phi(0)=0$, Taylor's formula gives on a sufficiently small neighborhood $\Omega_0$ the point-wise estimate $|\Phi(Y)|\leq C|Y|^2$. Thus, H\"older's inequality and Sobolev embedding in combination with \eqref{gl-ex_non-lin_L2_eqn1} imply      
\begin{align*}
\|\Phi(Y(\tau))\|_{L^1}
&\leq C\|Y(\tau)\|_{L^2}^2
 \leq C \mathcal{X}^Y_s(T)^2(1+\tau)^{-1/2},\\
\|\Phi(Y(\tau))\|_{L^2}
&\leq C\|Y(\tau)\|_{L^\infty}\|Y(\tau)\|_{L^2}
 \leq C \mathcal{X}^Y_s(T)^2(1+\tau)^{-1/2}.
\end{align*}

From Lemma \ref{lem:sec8-structured-linear} we now obtain the estimate 
\begin{align*}
I_{\mathrm{lo}}&(t)
 \leq \int_0^t  \left\|P_{\mathrm{lo}} S(t-\tau)\mathcal J\Phi(Y(\tau))\right\|_{H^s} \,d\tau \\
&\leq C \int_0^t \left( (1+t-\tau)^{-3/4}\|\Phi(Y(\tau)\|_{L^1}  + e^{-\theta (t-\tau)}\|\Phi(Y(\tau)\|_{L^2}\right) \,d\tau,
\end{align*}
and applying the preceding bounds gives
\small 
\begin{align*}
I_{\mathrm{lo}}(t)
\leq C \mathcal{X}^Y_s(T)^2 \Big( \int_0^t (1+t-\tau)^{-3/4}(1+\tau)^{-1/2}\,d\tau  +\int_0^t e^{-\theta(t-\tau)}(1+\tau)^{-1/2}\,d\tau  \Big).  
\end{align*}
\normalsize 
Now, splitting the first integral at
$t/2$ gives
\begin{align*}
&\int_0^t(1+t-\tau)^{-3/4}(1+\tau)^{-1/2}\,d\tau\\
&\quad\leq C(1+t)^{-3/4}
   \int_0^{t/2}(1+\tau)^{-1/2}\,d\tau
 +C(1+t)^{-1/2}
   \int_{t/2}^t(1+t-\tau)^{-3/4}\,d\tau\\
&\quad\leq C(1+t)^{-1/4}.
\end{align*}
For the second integral we have
\[
\int_0^t e^{-\theta(t-\tau)}(1+\tau)^{-1/2}\,d\tau
 \leq C(1+t)^{-1/2},
\]
where we used on $[0,t/2]$ that $e^{-\theta(t-\tau)}\leq e^{-\theta t/2}$, and on $[t/2,t]$ that $(1+\tau)^{-1/2}\leq C(1+t)^{-1/2}$. Substituting these estimates into the preceding bound on $I_{\mathrm{lo}}(t)$, observing that $(1+t)^{-1/2}\leq(1+t)^{-1/4}$, we finally obtain
\[
I_{\mathrm{lo}}(t)\leq C \mathcal{X}^Y_s(T)^2(1+t)^{-1/4},
\]
which proves the lemma.
\end{proof}

\subsection{Proof of Proposition~\ref{Prop_gl-ex_bootstrap}} \label{Sec_gl-ex_bootstrapProp} \label{subsec:sec8-proof-proposition}

Recall that $X_s(T)  \equiv \sup\limits_{0\leq t\leq T}(1+t)^{1/4}\|W(t)\|_{H^s}$,  $s\geq2$. Proposition \ref{Prop_gl-ex_bootstrap} asserts that there exist constants $\eta>0$ and $C>0$, independent of $T$, such that every compatible solution of \eqref{gl-ex_PDE-W} on $[0,T]$ which satisfies $X_s(T)\leq\eta$, satisfies as well
\begin{equation}\label{eq:sec8-bootstrap-combined_recall}
X_s(T)\leq C\delta+CX_s(T)^2,
\end{equation}
where $\delta \equiv \|W_0\|_{L^1}+\|W_0\|_{H^{s}} \leq \delta_0$.

We first establish the analog bootstrap estimates in the \(Y\)-variables on 
$$
\mathcal{X}^Y_s(T)=\sup_{0\leq t\leq T}(1+t)^{1/4}\|Y(t)\|_{H^s},
$$ 
as defined in \eqref{eq:sec8-XK-Y}. We then deduce the estimate in the corresponding $W$ variables, $X_s(T)$ using the diffeomorphism estimates \eqref{gl-ex_eqn_diffeo-norm} of Lemma \ref{gl-ex_Lemma_diffeo}. In particular, \eqref{gl-ex_eqn_diffeo-norm} gives the norm equivalence 
$$
C^{-1}X_s(T) \leq \mathcal X_s^Y(T) \leq CX_s(T),
$$
(written as $\mathcal{X}^Y_s(T) \asymp  X_s(T)$), and implies on initial data by \eqref{gl-ex_diffeo_norms_data} that 
\begin{align} \label{eq:sec8-Y0-data-bound} 
\|Y_0\|_{L^1}+\|Y_0\|_{H^{s}} &\le C\left( \|W_0\|_{L^1} +\|W_0\|_{H^{s}} +\|W_0\|_{H^{s}}^2 \right) \cr 
&\leq C(\delta+\delta^2) \leq 2C \delta, 
\end{align} 
provided $\delta$ is sufficiently small.

To start, recall Duhamel's solution formula \eqref{eq:sec8-Duhamel},  
\begin{equation}\label{gl-ex_Duhamel_again}
Y(t)= S(t)Y_0+ \int_0^tS(t-\tau)\mathcal J\Phi(Y(\tau))\,d\tau.
\end{equation}
We estimate the linear term in \eqref{gl-ex_Duhamel_again} by use of Proposition~\ref{prop:sec8-linear-decay} together with \eqref{eq:sec8-Y0-data-bound}, which yields 
\[
\|S(t)Y_0\|_{H^s}     
\ \leq \ C(1+t)^{-1/4}\|Y_0\|_{L^1} +Ce^{-\theta t}\|Y_0\|_{H^s}
\ \leq \ C\delta(1+t)^{-1/4}.
\]
For the low frequency part of the nonlinear term in Duhamel's formula, we use \eqref{eq:sec8-Duhamel-nonlinear-bound} of Lemma~\ref{lem:sec8-Duhamel-nonlinear} to obtain for $0\leq t\leq T$
\begin{equation}  \nonumber 
\left\| P_\text{lo} \int_0^tS(t-\tau)\mathcal J\Phi(Y(\tau))\,d\tau\right\|_{H^s} 
\leq C \; \mathcal{X}^Y_s(T)^2\; (1+t)^{-1/4}.
\end{equation}
Combing the above two estimate then gives
\begin{align}\label{eq:sec8-XY-bound}
\sup_{0\leq t\leq T}(1+t)^{1/4}\|P_\text{lo}Y(t)\|_{H^s} 
\leq C \Big(\delta+ \mathcal{X}^Y_s(T)^2 \Big).
\end{align}

To control the high frequency part we use the energy estimates developed in Section \ref{Sec_gl-ex_energy} below, from which we obtain the following estimate for $\eta>0$ sufficiently small,
\begin{equation}\label{eq:sec8-energy-integrated_preview}
 \|P_{\mathrm{hi}}Y(t)\|_{H^{s}}^2 \leq Ce^{-ct}\|Y_0\|_{H^{s}}^2 +C\int_0^t e^{-c(t-\tau)}\|Y(\tau)\|_{H^s}^4\,d\tau.
\end{equation}
We assume for the moment estimate \eqref{eq:sec8-energy-integrated_preview} to be valid. The proof of \eqref{eq:sec8-energy-integrated_preview} is postponed to Lemma \ref{lem:sec8-high-energy} below.  From \eqref{eq:sec8-energy-integrated_preview} and the definition of $\mathcal{X}^Y_s(T)$, we obtain
\begin{align}
\|P_{\mathrm{hi}}Y(t)\|_{H^{s}}^2
&\leq Ce^{-ct}\|Y_0\|_{H^{s}}^2 +C\mathcal{X}^Y_s(T)^4 \int_0^te^{-c(t-\tau)}(1+\tau)^{-1}\,d\tau \nonumber\\
&\leq C e^{-ct} \delta^2 + C\mathcal{X}^Y_s(T)^4(1+t)^{-1},
\nonumber
\end{align}
and taking square roots gives
\begin{equation}\label{eq:sec8-KY-bound}
\|P_{\mathrm{hi}}Y(t)\|_{H^{s}} \leq C e^{-ct} \delta+C\mathcal{X}^Y_s(T)^2 (1+t)^{-1/2}.
\end{equation}

Combining \eqref{eq:sec8-KY-bound} with \eqref{eq:sec8-XY-bound} gives the desired bound in the $Y$ variables,
\begin{equation}\label{eq:sec8-XY-bound-2}
\mathcal{X}^Y_s(T) \leq C \Big( \delta + \mathcal{X}^Y_s(T)^2 \Big).
\end{equation}
Choosing now $\eta>0$ small enough for the norm equivalence $\mathcal{X}^Y_s(T) \asymp  X_s(T)$ of Lemma~\ref{gl-ex_Lemma_diffeo} to hold, and small enough to imply $\sup_{0\leq\tau\leq T}\|Y(\tau)\|_{H^s}\leq\eta_0$, (the incoming bound of the energy estimate \eqref{eq:sec8-energy-integrated_preview} of Lemma \ref{lem:sec8-high-energy}), we conclude with the sought after estimate \eqref{eq:sec8-bootstrap-combined_recall}. 
This completes the proof of Proposition \ref{Prop_gl-ex_bootstrap}, once the proof of our energy estimate \eqref{eq:sec8-energy-integrated_preview} is completed in Subsection \ref{Sec_gl-ex_energy} below.
\hfill $\Box$

\subsection{Kawashima type energy estimates}  \label{Sec_gl-ex_energy}

The goal of this subsection is the proof of estimate \eqref{eq:sec8-energy-integrated_preview}, in order to complete the proof of Proposition \ref{Prop_gl-ex_bootstrap}. We begin by constructing the algebraic high-frequency symmetrizer used in the nonlinear energy estimate.  In the exact $Y$ variables let $q=(e,m,z)^{\mathsf T}$ and define          
\begin{equation}\nonumber
\mathcal G(q,u)=
\begin{pmatrix}
T^{01}(\rho(q,u),u)-\epsilon v(u)z\\
T^{11}(\rho(q,u),u)-p(\bar{\rho})-\epsilon z\\
-\omega(q,u)
\end{pmatrix},
\qquad \rho(q,u)=\bar{\rho}+e-v(u)m.
\end{equation}
The first three equations of \eqref{gl-ex_PDE-Y} are
\begin{equation}\nonumber
q_t+\partial_x\mathcal G(q,u)=0.
\end{equation}
Since $u_x=z$, this is equivalent to
\begin{equation}\label{eq:sec8-q-quasilinear}
q_t+\mathcal A(q,u)q_x=\mathcal C(q,u)q,
\end{equation}
where
\begin{equation}\nonumber
\mathcal A(q,u)=D_q\mathcal G(q,u),
\qquad
\mathcal C(q,u)=-\mathcal G_u(q,u)e_3^{\mathsf T},
\qquad e_3=(0,0,1)^{\mathsf T}.
\end{equation}
Writing $a=T^{01}_{\rho}(\rho,u)$ and $b=T^{11}_{\rho}(\rho,u)$, one has the
explicit formula
\begin{equation}\nonumber
\mathcal A(q,u)=
\begin{pmatrix}
a&-va&-\epsilon v\\
b&-vb&-\epsilon\\
\epsilon^{-1}a&-\epsilon^{-1}(1+va)&0
\end{pmatrix}.
\end{equation}
A direct computation, using the identities $T^{01}_{\rho}-vT^{11}_{\rho}=v$ and $T^{00}_{\rho}T^{11}_{\rho}-(T^{01}_{\rho})^2=p'(\rho)$ as well as $vT^{01}_{\rho}-T^{00}_{\rho}=-1$, gives
\begin{equation}\nonumber
\det(\lambda I-\mathcal A(q,u)) =(\lambda-v(u))(\lambda^2-1).
\end{equation}
Thus the three principal speeds are $-1$, $v(u)$, and $1$ and remain uniformly
separated near the rest state.

\begin{Lemma} \label{lem:sec8-frozen-symmetrizer}
Let
\begin{equation}\nonumber
 g(q,u)\equiv\partial_u\mathcal G(q,u),
 \qquad
 \mathcal C(q,u)=-g(q,u)e_3^{\mathsf T},
 \qquad e_3=(0,0,1)^{\mathsf T}.
\end{equation}
There exist a neighborhood $\Omega_0$ of the equilibrium $(q,u)=0$ and a
smooth symmetric matrix $\mathsf S=\mathsf S(q,u)$ such that
\begin{equation}\label{eq:sec8-classical-S-coercive}
 c_0 I\leq \mathsf S(q,u)\leq C_0 I,
 \qquad (q,u)\in\Omega_0,
\end{equation}
and
\begin{equation}\label{eq:sec8-classical-S-identities}
 \mathsf S\mathcal A=(\mathsf S\mathcal A)^{\mathsf T},
 \qquad
 \mathsf Sg=e_3.
\end{equation}
The constants $c_0,C_0>0$ are uniform on $\Omega_0$.  After $\Omega_0$
has been fixed, the cutoff frequency $R$ in the definition of
$P_{\mathrm{hi}}$ may be chosen arbitrarily large.
\end{Lemma}

\begin{proof}
The characteristic polynomial computed above shows that the eigenvalues of
$\mathcal A(q,u)$ are
\begin{equation}\nonumber
 -1,\qquad v(u),\qquad 1.
\end{equation}
They are real and pairwise distinct in a sufficiently small neighborhood of
the rest state.  Hence there is a smooth real eigenvector matrix
\begin{equation}\nonumber
 \mathsf R(q,u)=\bigl(r_-(q,u),r_0(q,u),r_+(q,u)\bigr)
\end{equation}
satisfying
\begin{equation}\label{eq:sec8-classical-R-diagonalizes}
 \mathcal A\mathsf R
 =\mathsf R\Lambda,
 \qquad
 \Lambda=\operatorname{diag}(-1,v(u),1).
\end{equation}
At the equilibrium, the eigenvectors may be chosen as the columns of
\begin{equation}\label{eq:sec8-classical-R-zero}
 \mathsf R(0)=
 \begin{pmatrix}
  0&\epsilon\sigma^{-2}&0\\
  \epsilon&0&-\epsilon\\
  1&1&1
 \end{pmatrix}.
\end{equation}
In particular, the third component of each eigenvector is nonzero.  By
continuity, after shrinking the state-space neighborhood, the eigenvectors
can be normalized smoothly so that
\begin{equation}\label{eq:sec8-classical-normalization}
 e_3^{\mathsf T}r_j(q,u)=1,
 \qquad j\in\{-,0,+\},
\end{equation}
or equivalently, $\mathsf R(q,u)^{\mathsf T}e_3=(1,1,1)^{\mathsf T}$.

Now introduce the vector,
\begin{equation}\nonumber
 b(q,u)\equiv\mathsf R(q,u)^{-1}g(q,u) \equiv \begin{pmatrix} b_-(q,u)\\ b_0(q,u)\\ b_+(q,u) \end{pmatrix}, 
\end{equation}
and define on a sufficiently small neighborhood
\begin{equation}\nonumber
\mathsf D(q,u)  \equiv \operatorname{diag} \left( \frac{1}{b_-(q,u)}, \frac{1}{b_0(q,u)}, \frac{1}{b_+(q,u)} \right). 
\end{equation}
Indeed, the matrix $\mathsf D(q,u)$ is well-defined, since at the equilibrium
\begin{equation}\nonumber
 g(0)=
 \begin{pmatrix}
  \kappa\\0\\\kappa/\epsilon
 \end{pmatrix},
 \qquad \text{and} \qquad 
b(0)=\frac{\kappa}{\epsilon}
 \begin{pmatrix}
  (1-\sigma^2)/2\\[1mm]
  \sigma^2\\[1mm]
  (1-\sigma^2)/2
 \end{pmatrix},
\end{equation}
and the assumptions $\kappa \equiv \bar{\rho}+p(\bar{\rho})>0$ and $0<\sigma^2<1$ imply that all three entries of $b(0)$ are strictly positive. Thus, shrinking the neighborhood gives $b_j(q,u)\geq b_*>0$, for $j\in\{-,0,+\}$.

We now define  
\begin{equation}\label{eq:sec8-classical-S-def}
 \mathsf S(q,u)
 \equiv\mathsf R(q,u)^{-\mathsf T}\mathsf D(q,u)\mathsf R(q,u)^{-1}.
\end{equation}
The matrix $\mathsf S$ is smooth, symmetric, and positive definite.  Since
$\Omega_0$ may be taken pre-compact, its smallest and largest
eigenvalues are uniformly bounded away from zero and infinity on $\overline{\Omega_0}$, which proves \eqref{eq:sec8-classical-S-coercive}. Using \eqref{eq:sec8-classical-R-diagonalizes}, we obtain
\begin{equation}\nonumber
 \mathsf S\mathcal A
 =\mathsf R^{-\mathsf T}\mathsf D\Lambda\mathsf R^{-1}.
\end{equation}
The middle matrix $\mathsf D\Lambda$ is real diagonal, and hence
$\mathsf S\mathcal A$ is symmetric.  Finally,
\begin{align*}
 \mathsf Sg
 &=\mathsf R^{-\mathsf T}\mathsf D\mathsf R^{-1}g
  =\mathsf R^{-\mathsf T}\mathsf Db
  =\mathsf R^{-\mathsf T}(1,1,1)^{\mathsf T}
  =e_3,
\end{align*}
where the last identity follows from \eqref{eq:sec8-classical-normalization}.
This proves \eqref{eq:sec8-classical-S-identities}.
\end{proof}

In order to introduce an energy functional that controls the necessary higher order derivative terms in our proof of Proposition \ref{Prop_gl-ex_bootstrap} we introduce the variable
\begin{equation}\label{eq:sec8-Zj-def}
Q_s:=P_{\mathrm{hi}}\partial_x^s q,
\end{equation}
any integer $s\geq 2$. Applying $P_{\mathrm{hi}}\partial_x^s$ to
\eqref{eq:sec8-q-quasilinear}, Leibniz rule shows that $Q_s$ satisfies
\begin{equation}\label{eq:sec8-paralinearized}
(Q_s)_t+\mathcal A \partial_x Q_s =-g(Q_s)_3+\mathcal R_s,
\end{equation}
for $\mathcal A \equiv \mathcal A(q,u)$ and $g\equiv g(q,u)$, and some function $\mathcal{R}_s$ satisfying the following bound.

\begin{Lemma} \label{lem:sec8-tame-reduction}
Let $Y=(q,u)$ be a compatible $H^s$ solution and define $Q_s$ by \eqref{eq:sec8-Zj-def}, $s\geq 2$ an integer. Then $Q_s$ satisfies \eqref{eq:sec8-paralinearized} and there exists a constant $C>0$, independent of $t$ and the cut-off $R>1$, such that 
\begin{equation}\label{eq:sec8-RN-bound}
\|\mathcal R_s(t)\|_{L^2} \leq C\|Y(t)\|_{H^s}^2 .
\end{equation}
\end{Lemma}

\Proof 
Since $P_{\mathrm{hi}}$ is a fixed Fourier multiplier, it commutes with $\partial_t$ and $\partial_x$. Thus, applying $P_{\mathrm{hi}}\partial_x^s$ to \eqref{eq:sec8-q-quasilinear}, Leibniz rule shows that $Q_s$ satisfies \eqref{eq:sec8-paralinearized} for 
\begin{align} \label{eq:sec8-Rj-exact}
\mathcal R_s ={}&  -\sum_{\ell=1}^s\binom{s}{\ell}  P_{\mathrm{hi}}\!\Big((\partial_x^\ell\mathcal A)  \partial_x^{s-\ell+1}q\Big) -\sum_{\ell=1}^s\binom{s}{\ell}  P_{\mathrm{hi}}\!\Big((\partial_x^\ell g)  \partial_x^{s-\ell}z\Big) \nonumber\\
& -[P_{\mathrm{hi}},g]\partial_x^s z   -[P_{\mathrm{hi}},\mathcal A]\partial_x^{s+1}q, 
\end{align}
where commutators are defined by $[P_{\mathrm{hi}},f]h \equiv P_{\mathrm{hi}}(f h) - f(P_{\mathrm{hi}}h)$ for functions $f$ and $h$, and hence
$$
\widehat{[P_{\mathrm{hi}},f]h}(\xi) = \frac1{2\pi} \int_{\mathbb R} \bigl(\varphi_R(\xi)-\varphi_R(\eta)\bigr) \widehat f(\xi-\eta)\widehat h(\eta)\,d\eta.
$$ 
It remains to estimate the $L^2$ norms of each term in \eqref{eq:sec8-Rj-exact}.

The finite sums in \eqref{eq:sec8-Rj-exact} are controlled by the following Moser estimates
\begin{align}
\sum_{\ell=1}^s \| (\partial_x^\ell \mathcal A)\partial_x^{s-\ell+1}q\|_{L^2}
&\leq C\bigl( \|\mathcal A_x\|_{L^\infty}\|q\|_{H^s} +\|\mathcal A_x\|_{H^{s-1}}\|q_x\|_{L^\infty} \bigr), \label{eq:sec8-integer-Moser-transport}\\
\sum_{\ell=1}^s \| (\partial_x^\ell g)\partial_x^{s-\ell}z\|_{L^2}
&\leq C\bigl( \|g_x\|_{L^\infty}\|z\|_{H^{s-1}} +\|g_x\|_{H^{s-1}}\|z\|_{L^\infty} \bigr), \label{eq:sec8-integer-Moser-source}
\end{align}
which can be derived by Leibniz' formula and Sobolev embedding. These two estimates control the binomial sums in \eqref{eq:sec8-Rj-exact}. The second to last term in \eqref{eq:sec8-Rj-exact} is controlled directly by 
\begin{equation} \label{gl-ex_energy_comm-est1}
\|[P_{\mathrm{hi}},g]\partial_x^s z\|_{L^2}
=\|[P_{\mathrm{hi}},g-g(0)]\partial_x^s z\|_{L^2}
\leq C\|g-g(0)\|_{L^\infty}\|z\|_{H^s},
\end{equation}
since $P_{\mathrm{hi}}$ is bounded in $L^2$ and commutes with constants. 

To control the last commutator term, (containing one apparent excessive spatial derivative), we next prove the estimate     
\begin{equation} \label{gl-ex_energy_comm-est2}
\|[P_{\mathrm{hi}},\mathcal A]\partial_x^{s+1}q\|_{L^2} \leq C\|\mathcal A - \mathcal{A}(0)\|_{H^s}   \|q\|_{H^s},
\end{equation}
for some constant $C>0$ independent of the cut-off $R>1$. To prove \eqref{gl-ex_energy_comm-est2}, recall that $\widehat{P_{\mathrm{hi}}f}(\xi) =\varphi_R(\xi)\widehat f(\xi)$ for $\varphi_R(\xi)=\varphi(\xi/R)$, where $\varphi$ is the fixed smooth cutoff profile.  Based on the mean-value theorem together with the fact that $\varphi_R$ is constant outside the transition region $[R,2R]$, one can show for $\xi, \eta\in \R$ that  
$$
|\eta|\, |\varphi_R(\xi)-\varphi_R(\eta)| \leq C|\xi-\eta|,
$$
with a constant $C>0$, (depending on $\varphi'$, but independent of $R$).  Hence, setting $f=\mathcal A-\mathcal A(0)$ and $h=\partial_x^sq$,
\begin{align*}
\left| \widehat{[P_{\mathrm{hi}},f]\partial_xh}(\xi) \right| &=
\left|\frac{1}{2\pi}\int_{\mathbb R}\bigl(\varphi_R(\xi)-\varphi_R(\eta)\bigr) i\eta\,\widehat f(\xi-\eta)\widehat h(\eta)\,d\eta \right|\\
&\leq C\left(|\widehat{\partial_xf}|*|\widehat h|\right)(\xi).
\end{align*}
Young's convolution inequality and Plancherel's theorem therefore imply
$$
\|[P_{\mathrm{hi}},f]\partial_xh\|_{L^2} \leq C\|\widehat{\partial_xf}\|_{L^1}\|h\|_{L^2}, 
$$
and since $s\geq 2$ we have by the Cauchy-Schwartz inequality
$$
\|\widehat{\partial_xf}\|_{L^1} \leq \left(\int_{\mathbb R}\frac{\xi^2}{(1+\xi^2)^s}\,d\xi\right)^{1/2} \|f\|_{H^s} \leq C \|f\|_{H^s},
$$
and thus
$$
\|[P_{\mathrm{hi}},f]\partial_xh\|_{L^2} \leq C\|f\|_{H^s} \|h\|_{L^2}. 
$$
This, keeping in mind that the commutator vanishes on the constant $\A(0)$, gives the sought-after estimate \eqref{gl-ex_energy_comm-est2}. 

Finally, using that by Morrey's inequality $H^s(\mathbb R)\hookrightarrow W^{1,\infty}(\mathbb R)$, together with smoothness of components in $q$ and $z$, we obtain 
\begin{equation} \label{gl-ex_energy_Q-Lem_techeqn5}
\begin{aligned}
&\|\mathcal A-\mathcal A(0)\|_{H^s} +\|g-g(0)\|_{H^s} +\|\partial_x\mathcal A\|_{L^\infty}
 +\|\partial_xg\|_{L^\infty}        \leq C\|Y\|_{H^s},\\
&\|q\|_{H^s}+\|z\|_{H^s} +\|q_x\|_{L^\infty}+\|z\|_{L^\infty}    \leq C\|Y\|_{H^s}.
\end{aligned}
\end{equation}
Combining \eqref{gl-ex_energy_Q-Lem_techeqn5} with \eqref{eq:sec8-integer-Moser-transport} - \eqref{gl-ex_energy_comm-est2} proves the sought-after estimate \eqref{eq:sec8-RN-bound}.
\QED

We are now prepared to state our energy estimate, together with the high frequency estimate \eqref{eq:sec8-energy-integrated} which was crucial in our proof of Proposition \ref{Prop_gl-ex_bootstrap}.

\begin{Lemma} \label{lem:sec8-high-energy}
Let $s\geq2$ be an integer, and let $Y=(q,u)$ be a compatible $H^s$ solution, $q=(e,m,z)^{\mathsf T}$, contained within neighborhood fixed in Lemma \ref{lem:sec8-frozen-symmetrizer}.  There exist constants $\eta_0,c,C, C_E>0$ and a
high-frequency energy $\mathcal E_s^{\mathrm{hi}}(t)$ such that, whenever
\begin{equation}\nonumber
\sup_{0\leq\tau\leq T}\|Y(\tau)\|_{H^s}\leq\eta_0,
\end{equation}
one has
\begin{equation}\label{eq:sec8-energy-equivalence}
C_E^{-1}\|P_{\mathrm{hi}} Y(t)\|_{H^s}^2
\leq \mathcal E_s^{\mathrm{hi}}(t)
\leq C_E\|P_{\mathrm{hi}} Y(t)\|_{H^s}^2,
\end{equation}
and
\begin{equation}\label{eq:sec8-energy-differential}
\frac{d}{dt}\mathcal E_s^{\mathrm{hi}}(t) +c\mathcal E_s^{\mathrm{hi}}(t) \leq C\|Y(t)\|_{H^s}^4.
\end{equation}
Consequently, by integrating \eqref{eq:sec8-energy-differential}, we obtain the sought-after estimate
\begin{equation}\label{eq:sec8-energy-integrated}
\|P_{\mathrm{hi}} Y(t)\|_{H^s}^2
\leq Ce^{-ct}\|Y_0\|_{H^s}^2
+C\int_0^t e^{-c(t-\tau)}\|Y(\tau)\|_{H^s}^4\,d\tau.
\end{equation}
\end{Lemma}

\Proof
The following calculations should be understood to be carried out at the level of smooth Friedrichs mollifications, subsequently passing to the $H^s$ solution once the sought-after uniform estimates are obtained.

The uncompensated energy functional
\beq \label{gl-ex_energy_0} 
\mathcal E_{0,s}^{\mathrm{hi}}(t)
\equiv \frac12 \int_{\mathbb R}Q_s^{\mathsf T}S(q,u)Q_s\,dx
\eeq  
is insufficient to control the higher order derivatives, essentially since the right hand side in our four component balance law has rank one. To exploit that the damping effect of this component is propagated through the coupling of the equations, we introduce a Kawashima compensated energy functional. For this we first introduce
\begin{equation} \nonumber 
K_0\equiv 
\begin{pmatrix}
0&1&0\\ -1&0&-\epsilon\sigma^{-2}\\ 0&\epsilon\sigma^{-2}&0
\end{pmatrix}
\quad \text{and} \quad 
\mathcal A_0\equiv \mathcal A(0)=
\begin{pmatrix}
0&0&0\\ \sigma^2&0&-\epsilon\\ 0&-\epsilon^{-1}&0
\end{pmatrix},
\end{equation}
and observe that $K_0$ is skew-symmetric and satisfies
\begin{equation}\label{eq:sec8-Kawashima-algebra}
\frac12\bigl(K_0\mathcal A_0-
\mathcal A_0^{\mathsf T}K_0\bigr)
=D_0,
\qquad
D_0=\operatorname{diag}
\bigl(\sigma^2,\sigma^{-2},-\epsilon^2\sigma^{-2}\bigr).
\end{equation}
Now, choose
\begin{equation}\label{eq:sec8-vartheta-choice}
0<\vartheta<\frac{\sigma^2}{4\epsilon^2},
\end{equation}
and define our ``compensated'' energy as
\begin{equation}\label{eq:sec8-high-energy-def}
\mathcal E_s^{\mathrm{hi}}(t) 
\equiv \mathcal E_{0,s}^{\mathrm{hi}}(t) + \vartheta\; \mathcal I_s(t)
\end{equation}
where $\mathcal I_s(t)\equiv \langle K_0Q_{s-1}(t),Q_s(t)\rangle_{L^2}$, with $Q_s \equiv P_{\mathrm{hi}}\partial_x^s q$ and $Q_{s-1} \equiv P_{\mathrm{hi}}\partial_x^{s-1} q$ for $s\geq 2$.
\smallskip

We now prove the norm equivalence \eqref{eq:sec8-energy-equivalence}. For this, note that the support of $P_{\mathrm{hi}}$ is separated from the origin, which implies that
\begin{equation}\label{eq:sec8-top-high-equivalence}
\|P_{\mathrm{hi}} q\|_{H^s} \asymp \|Q_s\|_{L^2},
\qquad
\|Q_{s-1}\|_{L^2} \leq R^{-1}\|Q_s\|_{L^2}.
\end{equation}
Moreover, compatibility gives $\widehat z(\xi)=i\xi\widehat u(\xi)$, and therefore
\begin{equation}\label{eq:sec8-u-high-by-z}
\|P_{\mathrm{hi}} u\|_{H^s}  \leq C\|P_{\mathrm{hi}} z\|_{H^{s-1}}  \leq C\|Q_s\|_{L^2}.
\end{equation}
Together with \eqref{eq:sec8-top-high-equivalence}, this yields the norm equivalence
\begin{equation}\label{eq:sec8-high-norm-equivalence}
C^{-1}\|P_{\mathrm{hi}} Y\|_{H^s}  \leq \|Q_s\|_{L^2}  \leq C\|P_{\mathrm{hi}} Y\|_{H^s}.
\end{equation}
By Lemma \ref{lem:sec8-frozen-symmetrizer} and \eqref{eq:sec8-top-high-equivalence}, after possibly increasing the cutoff frequency $R$, we have
\begin{equation}\nonumber
|\mathcal I_s(t)| = |\langle K_0Q_{s-1},Q_s\rangle_{L^2}| \leq C R^{-1}\|Q_s\|_{L^2}^2
\end{equation}
and 
\begin{equation}\label{eq:sec8-Ehi-Z-equivalence}
\tilde{C}^{-1} \|Q_s\|_{L^2}^2 \leq \mathcal E_s^{\mathrm{hi}}(t) \leq \tilde{C}\|Q_s\|_{L^2}^2,
\end{equation}
form some constant $\tilde{C}>0$. Combining \eqref{eq:sec8-Ehi-Z-equivalence} with \eqref{eq:sec8-high-norm-equivalence} directly proves \eqref{eq:sec8-energy-equivalence}.

\smallskip
The main step remaining is to establish \eqref{eq:sec8-energy-differential}, the final estimate \eqref{eq:sec8-energy-integrated} then follows by integration. For this, we begin by differentiating the first term in \eqref{eq:sec8-high-energy-def} which, together with the equation for $Q_s$ in \eqref{eq:sec8-paralinearized}, gives
\begin{align}
\frac{d}{dt} \mathcal E_{0,s}^{\mathrm{hi}}(t) +\|(Q_s)_3\|_{L^2}^2 \ 
= \frac12\int_{\mathbb R}Q_s^{\mathsf T} \bigl(S_t+\partial_x(S\mathcal A)\bigr)Q_s\,dx +\int_{\mathbb R}Q_s^{\mathsf T}S\mathcal R_s\,dx, \label{eq:sec8-basic-energy-identity}
\end{align}
where we used symmetry of $S\mathcal A$ to integrate by parts, while $Sg=e_3$ of Lemma \ref{lem:sec8-frozen-symmetrizer} gives the exact dissipative term $-\|(Q_s)_3\|_{L^2}^2$. Since the balance law implies
\begin{equation}\nonumber
\|Y_t(t)\|_{H^{s-1}}
\leq C\|Y(t)\|_{H^s},
\end{equation}
and since $H^{s-1}(\mathbb R)\hookrightarrow L^\infty(\mathbb R)$ for $s\geq2$, we have
\begin{equation}\nonumber
\|S_t+\partial_x(S\mathcal A)\|_{L^\infty}
\leq C\|Y(t)\|_{H^s}.
\end{equation}
Using also \eqref{eq:sec8-RN-bound} in \eqref{eq:sec8-basic-energy-identity}, we obtain
\begin{equation}\label{eq:sec8-basic-energy-bound}
\frac{d}{dt} \mathcal E_{0,s}^{\mathrm{hi}}(t) +\|(Q_s)_3\|_{L^2}^2 
\leq C\|Y\|_{H^s}\|Q_s\|_{L^2}^2 +C\|Y\|_{H^s}^2\|Q_s\|_{L^2}.
\end{equation}
Equation \eqref{eq:sec8-basic-energy-bound} readily implies dissipation of $\|(Q_s)_3\|_{L^2}$, but to obtain decay of the remaining two components we need to turn to the full compensated energy functional \eqref{eq:sec8-high-energy-def}. 

\smallskip
To prove suitable damping for the first two components, we first rewrite our balance law \eqref{eq:sec8-paralinearized} as
\begin{equation} \label{gl-ex_energy-est_techeqn-9}
(Q_s)_t=-\mathcal A_0\partial_x Q_s +\mathcal F_s,
\end{equation}
for $\mathcal F_s  \equiv -(\mathcal A-\mathcal A_0)\partial_x Q_s -g(Q_s)_3 + \mathcal R_s$. Then, assuming for the moment $s\geq 3$, differentiation of $\mathcal I_s(t)= \langle K_0Q_{s-1}(t),Q_s(t)\rangle_{L^2}$, after substitution of \eqref{gl-ex_energy-est_techeqn-9} and subsequent  integration by parts, gives
\begin{align}
\frac{d}{dt} \mathcal I_s
&=-\langle K_0\mathcal A_0Q_s,Q_s\rangle_{L^2}  -\langle K_0Q_{s-1},\mathcal A_0\partial_xQ_s\rangle_{L^2} \nonumber\\
&\quad+\langle K_0\mathcal F_{s-1},Q_s\rangle_{L^2}   +\langle K_0Q_{s-1},\mathcal F_s\rangle_{L^2} \nonumber\\
&=-2\langle D_0Q_s,Q_s\rangle_{L^2}+\mathcal R_I, \label{eq:sec8-Iprime-adjacent}
\end{align}
where the remainder $\mathcal R_I$ is defined as the scalar function
$$
\mathcal R_I \equiv  \langle K_0\mathcal F_{s-1},Q_s\rangle_{L^2} +\langle K_0Q_{s-1},\mathcal F_s\rangle_{L^2}.
$$
We next estimate $\mathcal R_I$.  For the first term, using \eqref{eq:sec8-RN-bound} of Lemma \ref{lem:sec8-tame-reduction}, together with \eqref{eq:sec8-top-high-equivalence} and boundedness of $g$, we obtain 
\begin{equation}\nonumber
|\langle K_0\mathcal F_{s-1},Q_s\rangle_{L^2}| 
\leq C\bigl(R^{-1}+\|Y\|_{H^s}\bigr)\|Q_s\|_{L^2}^2 +C\|Y\|_{H^s}^2\|Q_s\|_{L^2}.
\end{equation}
For the second term, we control the leading order derivative term $(\mathcal A-\mathcal A_0)\partial_x Q_s$ in $\mathcal F_s$ by integration by parts as 
\begin{align*}
|\langle K_0Q_{s-1},& (\mathcal A-\mathcal A_0)\partial_xQ_s\rangle_{L^2}|  \cr 
&\leq  |\langle K_0Q_s,(\mathcal A-\mathcal A_0)Q_s\rangle_{L^2}| + |\langle K_0Q_{s-1},
(\partial_x\mathcal A)Q_s\rangle_{L^2}| \cr 
& \leq C\|Y\|_{H^s}\|Q_s\|_{L^2}^2,
\end{align*}
using \eqref{eq:sec8-top-high-equivalence} for the last inequality, while the lower order term $-g(Q_s)_3 + \mathcal R_s$ in $\mathcal{F}_s$ is controlled in a similar way as
\begin{align*}
|\langle K_0Q_{s-1}, \big(g(Q_s)_3 - \mathcal R_s\big)\rangle_{L^2}|  
 \leq C\bigl(R^{-1}+\|Y\|_{H^s}\bigr)\|Q_s\|_{L^2}^2 +C\|Y\|_{H^s}^2\|Q_s\|_{L^2}.
\end{align*}
Combining those estimates, we conclude with the sought-after bound
\begin{equation}\label{eq:sec8-RI-bound}
|\mathcal R_I| \leq C\bigl(R^{-1}+\|Y(t)\|_{H^s}\bigr)\|Q_s(t)\|_{L^2}^2 +C\|Y(t)\|_{H^s}^2\|Q_s(t)\|_{L^2},
\end{equation}
for $s\geq 3$. Now, when $s=2$, to make sense of $\mathcal F_{s-1} = \mathcal F_{1}$, define 
\begin{align} \nonumber  
\mathcal F_{1} &\equiv -(\mathcal A-\mathcal A_0)\partial_x Q_{1} -g(Q_{1})_3 + \mathcal R_1, \cr 
\mathcal R_1 &\equiv - P_{\mathrm{hi}}\big((\partial_x\mathcal A)  \partial_x q\big) -[P_{\mathrm{hi}},\mathcal A]\partial_x^{2}q -  P_{\mathrm{hi}}\big((\partial_x g) z\big) -[P_{\mathrm{hi}},g]\partial_x z,
\end{align} 
in analogy to \eqref{eq:sec8-Rj-exact}. Then, similar to Lemma \ref{lem:sec8-tame-reduction}, using the $L^2$-boundedness of $P_{\mathrm{hi}}$, and the preceding commutator estimates \eqref{gl-ex_energy_comm-est1} - \eqref{gl-ex_energy_comm-est2}, we obtain by  that 
\begin{align*}
\|\mathcal R_1\|_{L^2}
\leq C\Bigl(
&\|\partial_x\mathcal A\|_{L^\infty}\|\partial_xq\|_{L^2}
+\|\mathcal A-\mathcal A(0)\|_{H^2}\|\partial_xq\|_{L^2}\\
&+\|\partial_xg\|_{L^\infty}\|z\|_{L^2}
+\|g-g(0)\|_{H^2}\|z\|_{L^2}
\Bigr) ,
\end{align*}
which by $H^2(\mathbb R)\hookrightarrow W^{1,\infty}(\mathbb R)$, and smooth dependence of $\mathcal A$ and $g$ on $Y$, gives
$$
\|\mathcal R_1(t)\|_{L^2} \leq C\|Y(t)\|_{H^2}^2.
$$
From this, using that $\partial_t Q_1=-\mathcal A_0\partial_xQ_1+\mathcal F_1$, we finally obtain the sought-after bound \eqref{eq:sec8-RI-bound} for all integer $s\geq 2$.

\smallskip
We are now prepared to derive the sought-after energy estimate \eqref{eq:sec8-energy-differential}. For this, combining \eqref{eq:sec8-Iprime-adjacent} and \eqref{eq:sec8-RI-bound} with the basic energy estimate \eqref{eq:sec8-basic-energy-bound}, we obtain
\begin{align} \label{eq:sec8-pre-absorption-energy}
\frac{d}{dt}\mathcal E_s^{\mathrm{hi}} &+\|(Q_s)_3\|_{L^2}^2 +2\vartheta\langle D_0Q_s,Q_s\rangle_{L^2} \nonumber\\
&\leq C\bigl(R^{-1}+\|Y(t)\|_{H^s}\bigr)\|Q_s(t)\|_{L^2}^2 +C\|Y(t)\|_{H^s}^2\|Q_s(t)\|_{L^2}. 
\end{align}
From the explicit form of $D_0$ in \eqref{eq:sec8-Kawashima-algebra}, together with our choice $0<\vartheta<\frac{\sigma^2}{4\epsilon^2}$ in \eqref{eq:sec8-vartheta-choice}, we find that
\begin{align} \label{eq:sec8-high-dissipation-coercive}
\|(Q_s)_3\|_{L^2}^2 +& 2\vartheta\langle D_0Q_s,Q_s\rangle_{L^2} \cr 
=& 2\vartheta\sigma^2\|(Q_s)_1\|_{L^2}^2 +2\vartheta\sigma^{-2}\|(Q_s)_2\|_{L^2}^2 +\bigl(1-2\vartheta\epsilon^2\sigma^{-2}\bigr) \|(Q_s)_3\|_{L^2}^2 \cr 
\geq & c_0 \|Q_s\|_{L^2}^2,
\end{align}
for some constant $c_0>0$. From this in combination with \eqref{eq:sec8-pre-absorption-energy}, we obtain
\begin{align} \nonumber 
\frac{d}{dt}\mathcal E_s^{\mathrm{hi}}(t) +c_0 \|Q_s\|_{L^2}^2 
& \leq \frac{d}{dt}\mathcal E_s^{\mathrm{hi}} +\|(Q_s)_3\|_{L^2}^2 +2\vartheta\langle D_0Q_s,Q_s\rangle_{L^2} \cr 
&\leq C\bigl(R^{-1}+\|Y(t)\|_{H^s}\bigr)\|Q_s(t)\|_{L^2}^2 +C\|Y(t)\|_{H^s}^2\|Q_s(t)\|_{L^2},
\end{align}
and using the Young's inequality in the form
$$
C\; \|Y(t)\|_{H^s}^2 \|Q_s(t)\|_{L^2} \leq \frac{c_0}{4}\|Q_s(t)\|_{L^2}^2 + \frac{C^2}{c_0}\|Y(t)\|_{H^s}^4,
$$ 
we obtain by absorbing $\frac{c_0}{4}\|Q_s(t)\|_{L^2}^2$ into the left hand side that
\begin{align} \nonumber 
\frac{d}{dt}\mathcal E_s^{\mathrm{hi}}(t) + \frac{3}{4}c_0 \|Q_s\|_{L^2}^2 
&\leq C\Big(R^{-1}+\|Y(t)\|_{H^s} \Big)\|Q_s(t)\|_{L^2}^2 +C\|Y(t)\|_{H^s}^4.
\end{align}
Now, choosing $R$ sufficiently large and $\eta_0$ sufficiently
small, such that $C\bigl(R^{-1}+\|Y(t)\|_{H^s}\bigr)<c_0/4$, we absorb the first term on the right-hand side into the left hand side, and obtain
\begin{align} \nonumber 
\frac{d}{dt}\mathcal E_s^{\mathrm{hi}}(t) + \frac{c_0}{2} \|Q_s\|_{L^2}^2 
&\leq  C\|Y(t)\|_{H^s}^4.
\end{align}
Using finally $\mathcal E_s^{\mathrm{hi}} \leq \tilde{C} \|Q_s\|_{L^2}^2$ by the norm equivalence \eqref{eq:sec8-Ehi-Z-equivalence}, we obtain the sought-after energy estimate \eqref{eq:sec8-energy-differential} for a new constant $0<c\leq\frac{c_0}{2 \max\{1,\tilde{C}\}}$. Finally, multiplying \eqref{eq:sec8-energy-differential} by $e^{ct}$ and integrating from $0$ to $t$ gives
\begin{equation}\nonumber
\mathcal E_s^{\mathrm{hi}}(t)   \leq    e^{-ct}\mathcal E_s^{\mathrm{hi}}(0) + C\int_0^t e^{-c(t-\tau)}\|Y(\tau)\|_{H^s}^4 \,d\tau.
\end{equation}
The energy equivalence \eqref{eq:sec8-energy-equivalence} and $\mathcal E_s^{\mathrm{hi}}(0)\leq C\|Y_0\|_{H^s}^2$ now imply \eqref{eq:sec8-energy-integrated}.  This completes the proof.
\QED

\appendix

\section{Discussion of an extension to General Relativity}

We now give a brief discussion on an extension of the dissipative relativistic Euler equations to General Relativity. For this, define the ``viscosity tensor'' 
$$
T_\visc^{\mu\nu}  \equiv T^{\mu\nu} - \visc \nabla^\mu u^\nu,
$$ 
where $T^{\mu\nu}$ is again the energy momentum tensor of a perfect fluid \eqref{perfect_fluid} but with $\eta_{\mu\nu}$ replaced by a general Lorentzian metric $g_{\mu\nu}$ and $\nabla^\mu = g^{\mu\nu}\nabla_\nu$ denotes the covariant derivative (i.e., the Levi-Civita connection) of $g_{\mu\nu}$. The interaction of the matter fields $T_\visc^{\mu\nu}$ with the spacetime geometry--and hence the gravitational field--is governed by the Einstein equations coupled to the dissipative relativistic Euler equations,
\beq \label{Einstein-RAV-Euler}
G^{\mu\nu} = 8\pi T^{\mu\nu}_\visc 
\hspace{1cm} \text{and} \hspace{1cm}
\nabla_\nu T^{\mu\nu} = \visc\; \Box_g u^\mu, 
\eeq
where $G_{\mu\nu}$ denotes the Einstein curvature tensor of $g_{\mu\nu}$ and $\Box_g \equiv g^{\mu\nu} \nabla_\nu \nabla_\mu$ denotes the D'Alembert operator associated to $g_{\mu\nu}$. Since the Einstein curvature tensor is constructed to be divergence-free, the Einstein equations imply $T^{\mu\nu}_\visc$ to be divergence-free, and this implies the second equation in \eqref{Einstein-RAV-Euler} on the matter fields in $T^{\mu\nu}_\visc$, a covariant version of the dissipative relativistic Euler equations. Since the Einstein tensor is symmetric ($G_{\mu\nu}=G_{\nu\mu}$), the fluids described by \eqref{Einstein-RAV-Euler} are required to be irrotational in the sense that $\nabla_\mu u_\nu = \nabla_\nu u_\mu$, a additional constraint over $g(u,u)=-1$.

\section{Declarations}

\subsection*{Data Availability}

No data is associated to this manuscript.

\subsection*{Declaration of AI Use}

We used Chat GPT 5.6 Sol to draft, discuss and check mathematical results and methods of this paper, mostly in Section \ref{Sec_gl-ex}, and in particular for the energy estimates in Section \ref{Sec_gl-ex_energy}. The final version of the manuscript has been written and verified by the authors.  

\subsection*{Funding}
This research was supported by CityU Start-up Grant (7200748), CityU Strategic Research Grant (7005839), and by Hong Kong University Grants Council grants ECS-21306524 and GRF-11303326.

\subsection*{Acknowledgement}
We thank Blake Temple, Matthias Sroczinski and Heinrich Freist\"uhler for sharing their deep knowledge on relativistic dissipation.

\end{document}